\documentclass[12pt,a4paper,final]{amsart}

\input xy
\xyoption{all}

\DeclareMathAlphabet{\mathpzc}{OT1}{pzc}{m}{it}

\usepackage{amsfonts,amsmath,amssymb,indentfirst,mathrsfs,amscd}
\usepackage[mathscr]{eucal}
\usepackage{graphicx}
\usepackage{nicefrac}

\usepackage[utf8]{inputenc}
\usepackage[english]{babel}
\usepackage{booktabs}
\usepackage{multirow}
\usepackage{verbatim}

\usepackage[all]{xy}
\usepackage{graphicx}
\usepackage{stmaryrd}
\usepackage{enumitem}

\usepackage{empheq}
\usepackage{color}
\usepackage{xcolor}

\author[]{\'Angel Gonz\'alez-Prieto}
\address{Departamento de \'Algebra, Geometr\'ia y Topolog\'ia, Facultad de Ciencias Matem\'aticas, Universidad Complutense de Madrid, 28040 Madrid, Spain.}
\address{Instituto de Ciencias Matem\'aticas (CSIC-UAM-UC3M-UCM), C.\ Nicol\'as Cabrera 15, 28049 Madrid, Spain}
\email{angelgonzalezprieto@ucm.es}

\title[Pseudo-quotients and their application to character varieties]{Pseudo-quotients of algebraic actions\\and their application to character varieties}

\thanks{}


\DeclareMathOperator{\img}{Im \,}           %

\DeclareMathOperator{\Hom}{Hom\,}           %
\DeclareMathOperator{\tr}{tr\,}             %

\DeclareMathOperator{\Inn}{Inn}

\DeclareMathOperator{\Spec}{Spec}

\DeclareMathOperator{\Gr}{Gr}

\usepackage{amsmath}
\begin{document}

\newtheorem{thm}{Theorem}[section]
\newtheorem{prop}[thm]{Proposition}
\newtheorem{lem}[thm]{Lemma}
\newtheorem{cor}[thm]{Corollary}
\newtheorem{conjecture}{Conjecture}
\newtheorem*{theorem*}{Theorem}

\theoremstyle{definition}
\newtheorem{defn}[thm]{Definition}
\newtheorem{ex}[thm]{Example}
\newtheorem{as}{Assumption}

\theoremstyle{remark}
\newtheorem{rmk}[thm]{Remark}

\theoremstyle{remark}
\newtheorem*{prf}{Proof}

\newcommand{\iacute}{\'{\i}} %
\newcommand{\norm}[1]{\lVert#1\rVert} %

\newcommand{\lto}{\longrightarrow}
\newcommand{\hra}{\hookrightarrow}

\newcommand{\suchthat}{\;\;|\;\;}
\newcommand{\dbar}{\overline{\partial}}

\newcommand{\cA}{\mathcal{A}}
\newcommand{\cC}{\mathcal{C}}
\newcommand{\cD}{\mathcal{D}} %
\newcommand{\cF}{\mathcal{F}} %
\newcommand{\cG}{\mathcal{G}} %
\newcommand{\cO}{\mathcal{O}} %
\newcommand{\cL}{\mathcal{L}} %
\newcommand{\cM}{\mathcal{M}} %
\newcommand{\cN}{\mathcal{N}} %
\newcommand{\cP}{\mathcal{P}} %
\newcommand{\cS}{\mathcal{S}} %
\newcommand{\cU}{\mathcal{U}} %
\newcommand{\cX}{\mathcal{X}}
\newcommand{\cT}{\mathcal{T}}
\newcommand{\cV}{\mathcal{V}}
\newcommand{\cB}{\mathcal{B}}
\newcommand{\cR}{\mathcal{R}}
\newcommand{\cH}{\mathcal{H}}
\renewcommand{\k}{k}

\newcommand{\ext}{\mathrm{ext}} %
\newcommand{\x}{\times}

\newcommand{\mM}{\mathscr{M}} %

\renewcommand{\AA}{\mathbb{A}} %
\newcommand{\CC}{\mathbb{C}} %
\newcommand{\QQ}{\mathbb{Q}} %
\newcommand{\PP}{\mathbb{P}} %
\newcommand{\HH}{\mathbb{H}} %
\newcommand{\RR}{\mathbb{R}} %
\newcommand{\ZZ}{\mathbb{Z}} %
\newcommand{\NN}{\mathbb{N}} %
\newcommand{\TT}{\mathbb{T}} %

\newcommand{\acts}{\circlearrowright} %
\newcommand{\Stab}{\textrm{Stab}} %

\renewcommand{\lg}{\mathfrak{g}} %
\newcommand{\lh}{\mathfrak{h}} %
\newcommand{\lu}{\mathfrak{u}} %
\newcommand{\la}{\mathfrak{a}} %
\newcommand{\lb}{\mathfrak{b}} %
\newcommand{\lm}{\mathfrak{m}} %
\newcommand{\lgl}{\mathfrak{gl}} %
\newcommand{\too}{\longrightarrow}
\newcommand{\imat}{\sqrt{-1}} %

\newcommand\CBord[3]{\textbf{Bord}_{{#1 #3}}^{#2}}
\newcommand\CBordp[1]{\CBord{#1}{}{}}
\newcommand\CTub[3]{\textbf{Tub}_{{#1 #3}}^{#2}}
\newcommand\CTubp[1]{\CTub{#1}{}{}}
\newcommand\CClose[3]{\textbf{Cl}_{{#1 #3}}^{#2}}
\newcommand\CClosep[1]{\CClose{#1}{}{}}
\newcommand\Obj[1]{\textrm{Obj}(#1)}
\newcommand\Mor[1]{\textrm{Mor}(#1)}
\newcommand\Vect[1]{#1\textrm{-}\textbf{Vect}}
\newcommand\MHS[1]{\textbf{HS}}
\newcommand\HS{\textbf{HS}}
\newcommand\Var[1]{\textbf{Var}_{#1}}
\newcommand\CVar{\Var{\CC}}
\newcommand\PHM[2]{\cR_{#1}^p(#2)}
\newcommand\MHM[1]{\cR_{#1}}
\newcommand\K[1]{\textrm{K}#1}
\newcommand\Kh[1]{\hat{\textrm{K}}#1}

\newcommand\Ab{\textbf{Ab}}
\newcommand\GL[1]{\mathrm{GL}_{#1}}
\newcommand\SL[1]{\mathrm{SL}_{#1}}
\newcommand\PGL[1]{\mathrm{PGL}_{#1}}
\newcommand\Rep[1]{\mathfrak{X}_{#1}}

\newcommand\coh[1]{\left[#1\right]}
\newcommand\DHPol[1]{e\left(#1\right)}
\newcommand\eVect{\mathcal{E}}
\newcommand\e[1]{\eVect\left(#1\right)}
\newcommand\intMor[2]{\int_{#1}\,#2}

\newcommand\Dom[1]{\mathcal{D}_{#1}}

\newcommand\Xf[1]{{X}_{#1}}					%
\newcommand\Xs[1]{\mathfrak{X}_{#1}}							%

\newcommand\Xft[2]{\overline{{X}}_{#1, #2}} %
\newcommand\Xst[2]{\overline{\mathfrak{X}}_{#1, #2}}			%

\newcommand\Xfp[2]{X_{#1, #2}}			%
\newcommand\Xsp[2]{\mathfrak{X}_{#1, #2}}						%

\newcommand\Xfd[2]{X_{#1; #2}}			%
\newcommand\Xsd[2]{\mathfrak{X}_{#1; #2}}						%

\newcommand\Xfm[3]{\mathcal{X}_{#1, #2; #3}}		%
\newcommand\Xsm[3]{X_{#1, #2; #3}}					%

\newcommand\XD[1]{#1^{D}}
\newcommand\XDh[1]{#1^{\delta}}
\newcommand\XU[1]{#1^{UT}}
\newcommand\XP[1]{#1^{U}}
\newcommand\XPh[1]{#1^{\upsilon}}
\newcommand\XI[1]{#1^{\iota}}
\newcommand\XTilde[1]{#1^{\varrho}}
\newcommand\Xred[1]{#1^{r}}
\newcommand\Xirred[1]{#1^{ir}}
\newcommand\Id{\textrm{Id}}

\newcommand\RM[2]{R\left(\left.#1\right|#2\right)}
\newcommand\set[1]{\left\{#1\right\}}

\hyphenation{mul-ti-pli-ci-ty}

\hyphenation{mo-du-li}
\hyphenation{mo-no-dro-my}

\begin{abstract}
In this paper, we propose a weak version of quotient for the algebraic action of a group on a variety, which we shall call a pseudo-quotient. They arise when we focus on the purely topological properties of good GIT quotients regardless of their algebraic properties. The flexibility granted by their topological nature enables an easier identification in geometric constructions than classical GIT quotients. We obtain several results about the interplay between pseudo-quotients and good quotients. Additionally, we show that in characteristic zero pseudo-quotients are unique up to virtual class in the Grothendieck ring of algebraic varieties. As an application, we compute the virtual class of $\SL{2}(\k)$-character varieties for free groups and surface groups as well as their parabolic counterparts with punctures of Jordan type. 
\end{abstract}

\maketitle

\setcounter{section}{0}

\vspace{-0.5cm}
\normalsize

\section{Introduction}
\let\thefootnote\relax\footnotetext{\noindent \emph{2020 Mathematics Subject Classification}. Primary:
 14L24. %
 Secondary:
 14C25, %
 20G05. %
 
\emph{Keywords and phrases}: Geometric Invariant Theory, character varieties, Grothendieck ring, Hodge theory.}

Let $\Gamma$ be a finitely generated group and let $G$ be a reductive algebraic group over a field $k$. The set of representations $\rho: \Gamma \to G$ can be endowed with the structure of an algebraic variety, the so-called representation variety of $\Gamma$ into $G$ and denoted by $\Rep{G}(\Gamma)$. The group $G$ itself acts on $\Rep{G}(\Gamma)$ by conjugation so, taking the Geometric Invariant Theory (GIT) quotient of $\Rep{G}(\Gamma)$ by this action, we obtain the moduli space of representations $\cR_G(\Gamma) = \Rep{G}(\Gamma) \sslash G$, also known as the character variety.

These varieties have been studied for more than thirty years, in part due to their prominent role in non-abelian Hodge theory. In this context we take $k = \CC$, fix $\Gamma = \pi_1(\Sigma)$ where $\Sigma$ is a closed orientable surface, and set $\cR_G(\Sigma) = \cR_G(\pi_1(\Sigma))$. For $G = \GL{n}( \CC)$ (resp.\ $G = \SL{n}(\CC)$), the Riemann-Hilbert correspondence \cite{SimpsonI,SimpsonII} gives a real analytic correspondence between $\cR_G(\Sigma)$ and the moduli space of flat bundles of rank $n$ and degree $0$ (resp.\ with trivial determinant bundle). Analogously, the Hitchin-Kobayashi correspondence \cite{Simpson:1992,Corlette:1988} shows that $\cR_{\GL{n}(\CC)}(\Sigma)$ is real analytically equivalent to the moduli space of rank $n$ and degree $0$ Higgs bundles (see \cite{Hitchin}).

The situation is even richer in the so-called parabolic context. From the point of view of character varieties, a parabolic structure $Q$ on $\Sigma$ is a finite set of points $p_1, \ldots, p_s \in \Sigma$ and conjugacy classes $\lambda_1, \ldots, \lambda_s \subseteq G$. The $G$-representation variety with parabolic structure $Q$, denoted by $\Rep{G}(\Sigma, Q)$, is the space of representations $\rho: \pi_1(\Sigma - \left\{p_1, \ldots, p_s\right\}) \to G$ such that the image of the positively oriented loop around the puncture $p_i$ is forced to lie in the conjugacy class $\lambda_i$. Analogously, we can consider the associated parabolic character variety $\cR_G(\Sigma, Q)= \Rep{G}(\Sigma, Q) \sslash G$. In this context, the non-abelian Hodge correspondence can be extended to give new relations between the moduli spaces. For instance, if we remove a single point $p \in \Sigma$ and we take $\lambda = \{e^{2\pi i d/n}\Id\}$ for $d \in \NN$, we obtain the so-called twisted character variety. This variety turns out to be diffeomorphic to the moduli space of rank $n$ logarithmic flat bundles of degree $d$ with a pole at $p$ with residue $-\frac{d}{n}\Id$ and to the moduli space of degree $d$ Higgs bundles. 

Moreover, for an arbitrary number of punctures $p_1, \ldots, p_s$ and different semi-simple conjugacy classes of $G=\GL{n}(\CC)$, we obtain diffeomorphisms between moduli spaces of parabolic Higgs bundles with parabolic structures (with general weights) on $p_1, \ldots, p_s$ and with the moduli space of logarithmic flat connections with poles on the punctures \cite{Simpson:parabolic}. Additionally, other correspondences may also appear. For example, if $G=\SL{2}(\CC)$, the surface $\Sigma$ is an elliptic curve and $Q$ has two marked points with generic semi-simple conjugacy classes, then the variety $\cR_G(\Sigma, Q)$ is diffeomorphic to the moduli space of doubly periodic instantons through the Nahm transform \cite{Biquard-Jardim,Jardim}.

However, these correspondences are far from being algebraic so it is important to study the Hodge structures presented on these complex character varieties. For this purpose, it is customary to associate a polynomial {$e(X) \in  \ZZ[u, v]$} to each complex algebraic variety $X$ by taking the alternating sum of its compactly supported Hodge numbers. This polynomial is usually called the $E$-polynomial, or Deligne-Hodge polynomial, and should be seen as the Hodge-theoretic analog halfway between the Poincar\'e polynomial and the Euler characteristic.

Three different strategies exist in the literature to address the computation of $E$-polynomials of complex character varieties, usually called the arithmetic, geometric and quantum methods. The first strategy, the arithmetic method, was developed by \cite{Hausel-Rodriguez-Villegas:2008} by means of a theorem of Katz. Roughly speaking, this result states that if the number of points of a variety $X$, over the finite field of $q$ elements, is a polynomial in $q$, say $P(q) = |X(\mathbb{F}_q)|$, then the $E$-polynomial of $X(\CC)$ is $e(X(\CC)) = P(uv)$. Following this method, an expression for these polynomials for twisted character varieties was given in terms of generating functions for $G=\GL{n}(\CC)$ \cite{Hausel-Rodriguez-Villegas:2008} and for $G=\SL{n}(\CC)$ \cite{Mereb}. This arithmetic approach has also been applied for counting the number of stable Higgs bundles on a compact Riemann surface defined over a finite field \cite{Schiffmann2016}. However, in the general parabolic case much remains to be done. One of the most important advances was \cite{Hausel-Letellier-Villegas} where the $E$-polynomials were computed for $G=\GL{n}(\CC)$ and generic semi-simple marked points. 

The second approach is the geometric method, initiated by Logares, Mu\~noz and Newstead in \cite{LMN}. The strategy is to focus on the computation of the $E$-polynomial of $\Rep{G}(\Sigma)$ by decomposing it into simpler geometric strata for which the $E$-polynomial can be easily computed.
Following these ideas, in the case $G=\SL{2}(\CC)$, they computed the $E$-polynomials of the character variety for a single marked point and underlying surface of genus $g=1,2$. Later, these computations were extended for two marked points and $g=1$ \cite{LM} and for a marked point and $g=3$ \cite{MM:2016}. For arbitrary genus and at most one marked point, the case $G=\SL{2}(\CC)$ was accomplished in \cite{MM} and the case $G=\PGL{2}(\CC)$ in \cite{Martinez:2017}. Moreover, using a combination of both the arithmetic and the geometric approaches, explicit expressions of the $E$-polynomials for character varieties have been computed for orientable surfaces with $G=\GL{3}(\CC),\SL{3}(\CC)$ and for non-orientable surfaces with $G=\GL{2}(\CC),\SL{2}(\CC)$ \cite{Baraglia-Hekmati:2016}.
In order to compute the $E$-polynomial of the character variety from one of the representation varieties, the authors stratified the representation variety into its irreducible locus (where the action is essentially free) and its reducible locus (where the quotient reduces to the study of semi-simple reducible representations). Despite their achievements, the arguments in these papers are very handcrafted and it is not clear how to generalize these results to higher rank or to Hodge structures. Nevertheless, these ideas lie in the core of the present paper.

It is worth mentioning that these stratification methods have also been successfully applied in other contexts through the non-abelian Hodge correspondence. For instance, in \cite{Gonzalez-Martinez2012Jan} the author computed explicit formulas for the $E$-polynomial of the moduli space of rank $2$ stable vector bundles by constructing a suitable stratification of the semi-stable locus. 

The third approach, the quantum method, is based on a novel framework described in \cite{GPLM-2017}. In that paper it is proven that there exists a lax monoidal Topological Quantum Field Theory (TQFT) that computes the $E$-polynomial of the representation variety of any compact manifold (maybe equipped with a parabolic structure). Using this TQFT, in \cite{GP-2018} (see also \cite{Gonzalez-Prieto:Thesis}) explicit expressions were computed for the $E$-polynomials of the representation varieties for $\SL{2}(\CC)$ and parabolic structures with any number of punctures and Jordan type conjugacy classes. These cases cannot be addressed with the arithmetic and geometric approaches.

Furthermore, the quantum method is able to go forward and to compute a subtler invariant than the $E$-polynomial: the virtual class of the representation variety. Recall that we can form the Grothendieck ring of algebraic varieties over a ground field $k$, denoted by $\K{\Var{k}}$. This is the ring generated by isomorphism classes of algebraic varieties, denoted by $[X]$ and called virtual classes, modulo the `cut-and-paste' relations: $[X] = [Y] + [U]$ if $X = Y \sqcup U$ with $Y \subseteq X$ a closed subvariety and $U \subseteq X$ an open set. In the case $k = \CC$ there is a natural ring homomorphism {$\K{\Var{\CC}} \to \ZZ[u,v]$} sending $[X] \mapsto e(X)$, so in this sense virtual classes are stronger than $E$-polynomials. The crucial point here is that the TQFT developed in \cite{GP-2018} not only computes the $E$-polynomial but also the virtual classes of the representation varieties, that is $\left[\Rep{G}(\Sigma, Q)\right] \in \K{\Var{k}}$. This enables the calculation of the virtual classes of character varieties, a problem unapproachable with the arithmetic and geometric techniques.

In this paper, we compute the virtual class of $\SL{2}(k)$-character varieties of closed orientable surfaces, starting from those (recently obtained) of the corresponding $\SL{2}(k)$-representation varieties. We achieve this in two steps: firstly we use the method of stratification, mentioned above and introduced in \cite{LMN}, which consists in decomposing the representation variety into its irreducible and reducible loci; secondly we show how the virtual class of the GIT quotients of the irreducible and reducible loci can be computed from the (known) virtual classes of the irreducible and reducible loci. For this second step, we introduce the notion of `pseudo-quotients' which we present in the first part of this paper. We have chosen to present the theory of pseudo-quotients in greater generality than is required for its application to character varieties as we believe that it is interesting in its own right and moreover that it could be useful for studying other
types of character varieties (for example character varieties involving non-reductive groups) or even other moduli spaces.

The aim of pseudo-quotients is to provide a weak version of a classical GIT quotient with more flexibility. Suppose that $X$ is an algebraic variety with an action of an algebraic group $G$. A pseudo-quotient is a surjective $G$-invariant regular morphism
$$
	\pi: X \to \overline{X}
$$  
onto an algebraic variety $\overline{X}$ such that, for any disjoint $G$-invariant closed sets $W_1, W_2 \subseteq X$ we have that the Zariski closure of their images $\overline{\pi(W_1)}, \overline{\pi(W_2)}$ are also disjoint. It is worth noting that a good quotient (in the standard GIT sense) is then a pseudo-quotient with the additional algebraic property that for any open set $U \subseteq \overline{X}$, the map $\pi$ induces an isomorphism between regular functions on $U$ and $G$-invariant regular functions on $\pi^{-1}(U)$, i.e.\ $\pi^*: \cO(U) \cong \cO(\pi^{-1}(U))^G$. In this sense, the definition of pseudo-quotients only captures the topological properties of good quotients and gets rid of their algebraic part.

Thanks to its purely topological nature, pseudo-quotients behave better than good quotients with respect to several geometric constructions. The main results regarding pseudo-quotients (in characteristic zero) which we will prove in this paper are summarised in the theorem below.

\begin{theorem*}
Let $X$ be an algebraic variety with an action of an algebraic group $G$ (possibly non-reductive) and let $X \to \overline{X}$ be a pseudo-quotient for the action of $G$ on $X$. Then we have:
\begin{enumerate}
	\item If $X \to \overline{X}'$ is another pseudo-quotient for the action of G on X, then
$$
	[\overline{X}] = [\overline{X}']
$$
in the Grothendieck ring of algebraic varieties.
	\item More generally, if $X = Y  \sqcup U$ where $Y \subseteq X$ is closed and $U \subseteq X$ is open and saturated, then
$$
[\overline{X}] = [\overline{Y}] + [\overline{U}]
$$
for any pseudo-quotients $Y \to \overline{Y}$ and $U \to \overline{U}$ of $Y$ and $U$ respectively.
	\item If $Y \subseteq X$ is a subvariety, $H \leq G$ is a subgroup and $(Y, H)$ is a core for the action of $G$ on $X$ (see Proposition \ref{prop:core} for the definition), then the restriction of $X \to \overline{X}$ to $Y$ is a pseudo-quotient $Y \to \overline{X}$ for the action of $H$ on $Y$.
\end{enumerate}
\end{theorem*}

This notion of pseudo-quotient can be used to characterize quotients even if $G$ is not reductive. To be precise, suppose that there is an action of a non-reductive group $G$ on a complex algebraic variety $X$ and that $X$ can be stratified into subvarieties $X = X_1 \sqcup \ldots \sqcup X_n$ with each $X_i$ orbitwise-closed. Suppose that, by some topological argument, it is possible to prove that there exists a pseudo-quotient of $G$ on each $X_i$, call it $\overline{X}_i$. Then this completely characterizes the class of any quotient $\overline{X}$ of $X$ under $G$ in the Grothendieck ring of algebraic varieties as $[\overline{X}] = [\overline{X}_1] + \ldots + [\overline{X}_n]$. This approach is particularly useful for character varieties, as applied in \cite{GPLM-2020} or \cite{Vogel:2020}, since it allows us to give a meaning to the `virtual character variety' even if $G$ is not reductive.
Moreover, it is also interesting to combine pseudo-quotients with the techniques developed to construct non-reductive GIT quotients, as in \cite{Berczi-Dolan-Hawes-Kirwan:2016}, \cite{Dolan-Kirwan:2007} and \cite{Kirwan:2009}. Once the GIT quotient has been constructed by means of these powerful techniques, pseudo-quotients may be used to extract topological information from these sophisticated quotients, e.g.\ their virtual class or natural stratifications for the action.

As an application to character varieties, we have a decomposition of the representation variety as $\Rep{G}(\Gamma) = \Xred{\Rep{G}}(\Gamma) \sqcup \Xirred{\Rep{G}}(\Gamma)$, where $\Xred{\Rep{G}}(\Gamma)$ denotes the set of reducible representations and $\Xirred{\Rep{G}}(\Gamma)$ the set of irreducible ones. In that case, the results of Section \ref{sec:rep-var} will imply that
$$
	[\cR_G(\Gamma)] = [\Xred{\Rep{G}}(\Gamma) \sslash G] + [\Xirred{\Rep{G}}(\Gamma) \sslash G].
$$

At this point, each stratum can be analyzed separately. For the stratum $\Xirred{\Rep{G}}(\Gamma)$, the situation is simple since the action on it is closed and (essentially) free. Hence, the quotient map $\Xirred{\Rep{G}}(\Gamma) \to \Xirred{\Rep{G}}(\Gamma) \sslash G$ is a locally trivial fibration in the analytic topology with trivial monodromy and fiber the inner automorphism group $\Inn(G)$ of $G$.

For the stratum $\Xred{\Rep{G}}(\Gamma)$, the situation is a bit harder. The idea here is that the closures of the orbits of elements of $\Xred{\Rep{G}}(\Gamma)$ always intersect the subvariety of semi-simple reducible representations. This is precisely the setting of a core, as mentioned above. Hence, the variety $\Xred{\Rep{G}}(\Gamma) \sslash G$ can be understood as a pseudo-quotient for the action of the symmetric group on the set of semi-simple reducible representations under permutation of its irreducible components. In particular, this means that the virtual class $[\Xred{\Rep{G}}(\Gamma) \sslash G]$ can be computed by means of the analysis of the quotient of the subvariety of semi-simple representations by a finite group, a task that can be accomplished via equivariant methods. Observe that the use of pseudo-quotients is crucial at this point since the available topological information of the action only allows us to prove that the semi-simple reducible representations are a pseudo-quotient for $\Xred{\Rep{G}}(\Gamma)$ and not a good quotient.

Using these ideas, in Section \ref{subsec:reducible-rep} we shall compute the virtual classes of character varieties of free groups and surface groups from the virtual classes of the corresponding representation varieties, for $G = \SL{2}(k)$. In particular, taking the Deligne-Hodge homomorphism {$e: \K{\CVar} \to \ZZ[u, v]$}, this calculation reproves the results of \cite{Cavazos-Lawton:2014} and \cite{MM} respectively. Moreover in Section \ref{sec:parabolic-rep} we will explore the parabolic case and we will compute the virtual classes of $\SL{2}(k)$-parabolic character varieties of free and surface groups with any number { of} punctures with conjugacy classes of Jordan type. In the case of surface groups, the obtained result is the following.

\begin{theorem*} 
Fix $k$ an algebraically closed field of characteristic zero and $G = \SL{2}(k)$. Let $\Sigma_g$ be the compact orientable surface of genus $g$, let $Q = \left\{(p_1, \lambda_1), \ldots, (p_s, \lambda_s)\right\}$ be a parabolic structure where $\lambda_i$ are conjugacy classes of $-Id$ or of the Jordan matrices $J_+, J_-$. Let $r_+$ be the number of $J_+$, let $r_-$ be the number of $J_-$ and $t$ the number of $-\Id$ in $Q$. Set $r = r_+ + r_-$ and $\sigma = (-1)^{t + r_-}$.

Let $q = \coh{\AA_k^1} \in \Kh{\Var{k}}$ be the virtual class of the affine line {in the localization $\Kh{\Var{k}}$ of the Grothendieck ring of algebraic varieties by $q, q+1$ and $q-1$. Then we have that the virtual classes of the $\SL{2}(k)$-character varieties in $\Kh{\Var{k}}$ are the following:}
\begin{itemize}
	\item If $\sigma = 1$, then
	\begin{align*}
	\coh{\cR_{\SL{2}(k)}(\Sigma_g, Q)} =& \,{\left(q^2 - 1\right)}^{2g + r -
2} q^{2g - 2} +\left(-1\right)^{r} 2^{2g}  {\left(q - 1\right)} q^{2g - 2}
{\left({1-\left(1-q\right)}^{r - 1}\right)}\\
&+ \frac{1}{2}{\left(q - 1\right)}^{2g +r - 2} q^{2g - 2} \, {\left(2^{2g} + q - 3\right)}  \\
&+ \frac{1}{2} {\left(q + 1\right)}^{2g +
r - 2} q^{2g - 2}\,
\left(2^{2g} + q - 1\right).
	\end{align*}
	\item If $\sigma = -1$, then
	\begin{align*}
	\hspace{-1.75cm}\coh{\cR_{\SL{2}(k)}(\Sigma_g, Q)} =& \left(-1\right)^{r-1}2^{2g - 1}  {\left(q + 1\right)}^{2g + r - 2} q^{2g - 2} \\ &+ {\left(q - 1\right)}^{2g + r - 2} q^{2g - 2}\left( {\left(q + 1\right)}^{2g + r - 2} + 2^{2g - 1} - 1\right).
\end{align*} 
\end{itemize}
\end{theorem*}

The $E$-polynomials of these varieties were only known in the case of a single marked point \cite{MM:2016} or two marked points on a genus one surface \cite{LM}.

As will be clear from the calculations, the general case of any conjugacy class does not present more difficulties than the case of Jordan type classes. However, as in the previous cases, in order to compute the virtual class of the character variety we need to know the virtual class of the associated representation variety. The problem is that in the case of semi-simple conjugacy classes a new interference phenomenon arises. This makes the calculations for the corresponding representation variety {more} difficult. Since dealing with these interferences requires a subtle analysis and developing novel techniques at the level of TQFTs, we have postponed these results to the paper \cite{GP-2019}.

As a final remark regarding the higher rank case, the methods developed in this paper can be easily generalized. As happens for the rank $m=2$ case, on the irreducible locus of the $\SL{m}(k)$-representation variety the action is essentially free, whereas on the reducible locus the role of the core is played by the semi-simple reducible representations up to the action of the symmetric group by permutation. Thereby, we expect that the methods developed in this paper, together with the TQFT approach, will help to address the higher rank case in the future.

\subsection*{Acknowledgements}

The author wants to thank David Ben-Zvi, Alejandro Calleja, Carlos Florentino, Marton Hablicsek, Frances Kirwan, Javier Mart\'inez, Mart\'in Mereb, Peter Newstead, Jesse Vogel and Thomas Wasserman for very useful conversations, and specially Sean Lawton for very fruitful discussions regarding Proposition \ref{prop:multiplicativity}. Moreover, the author wants to thank the anonymous referee of this paper for their careful reading of the manuscript and for his/her insightful comments.

Finally, I would also like to express my highest gratitude to my PhD advisors Marina Logares and Vicente Mu\~noz for their invaluable help, support and encouragement throughout the development of this paper.

The author acknowledges the hospitality of the School of Computing, Electronics and Mathematics at the University of Plymouth, where this work was completed during a research visit.
The author has been partially supported by a "\!la Caixa" scholarship for PhD studies in Spanish Universities from "\!la Caixa" Foundation LCF/BQ/DE15/10360013, by MINECO (Spain) Project MTM2015--63612--P, and by the \textit{Madrid Government (Comunidad de Madrid -- Spain)} under the Multiannual Agreement with the Universidad Complutense de Madrid in the line Research Incentive for Young PhDs in the context of the V PRICIT (Regional Programme of Research and Technological Innovation) through the project PR27/21-029. Various calculations in the paper were assisted by the computer algebra system SageMath.

\section{Review of Geometric Invariant Theory}
\label{sec:review-GIT}

Throughout this section we will work over an arbitrary algebraically closed field $k$.
The aim of this section is to review some of the most important notions of Geometric Invariant Theory (GIT for short). The definitions provided here follow Newstead's book \cite{Newstead:1978}.

In order to fix some notation, given an algebraic group $G$ acting algebraically on a variety $X$ we will denote by $Gx$ or $[x]_G$ the orbit of $x \in X$, and by $\overline{Gx}$ or $\overline{[x]}_G$ its Zariski closure. In general, the space of orbits $X/G$ does not have the structure of an algebraic variety so we need to consider more sophisticated quotients.

Recall that a \emph{categorical quotient} of $X$ by $G$ is a $G$-invariant regular morphism $\pi: X \to Y$ onto some algebraic variety $Y$ such that, for any $G$-invariant regular morphism $f: X \to Z$ with $Z$ and algebraic variety, there exists a unique regular morphism $\tilde{f}: Y \to Z$ such that the following diagram commutes
\[
\begin{displaystyle}
   \xymatrix
   {
   	X \ar[r]^{f} \ar[d]_{\pi} & Z \\
   	Y \ar@{--{>}}[ru]_{\tilde{f}}&
  	}
\end{displaystyle}
\]
Using this universal property, it follows that the categorical quotient, if it exists, is unique up to regular isomorphism.

\begin{rmk}
In this paper we will always work with categorical quotients within the category of algebraic varieties and regular morphisms. However, sometimes in the literature larger categories are considered, like the category of schemes.
\end{rmk}

The definition of a categorical quotient does not say anything about the geometry of $Y$. To capture these geometric properties, a regular morphism $\pi: X \to Y$ is called a \emph{good quotient} if it satisfies the following properties:
\begin{enumerate}[label=$\roman*)$,ref=$\roman*)$]
	\item $\pi$ is $G$-invariant.
	\item $\pi$ is surjective.
	\item For any open set $U \subseteq Y$, the map $\pi$ induces an isomorphism
	$$
		\pi^*: \cO_Y(U) \stackrel{\cong}{\longrightarrow} \cO_X(\pi^{-1}(U))^G \subseteq \cO_X(\pi^{-1}(U)),
	$$
where $\cO_X$ is the sheaf of regular functions on $X$.
	\item If $W \subseteq X$ is a closed $G$-invariant set, then $\pi(W) \subseteq Y$ is closed.
	\item Given two closed $G$-invariant subsets $W_1, W_2\subseteq X$, we have that $W_1 \cap W_2 = \emptyset$ if and only if $\pi(W_1) \cap \pi(W_2) = \emptyset$.
\end{enumerate}

If $\pi: X \to Y$ is a good quotient, then it is a categorical quotient \cite[Corollary 3.5.1]{Newstead:1978}. 

Finally, the third type {of} quotient considered in the literature is the so-called \emph{geometric quotient}. It is a good quotient $\pi: X \to Y$ such that, for any $y \in Y$, we have $\pi^{-1}(y) = Gy$. This last condition is sometimes referred to as $Y$ being an \emph{orbit space}. Geometric quotients are categorical so they are unique. From the properties of good quotients we get that a geometric quotient is the same as a good quotient where the action of $G$ on $X$ is closed, i.e.\ $Gx$ is a closed subset of $X$ for all $x \in X$.

\begin{rmk}
Sometimes in the literature, good quotients are referred to as good categorical quotients and geometric quotients are called good geometric quotients. Moreover, some authors (in particular \cite{Newstead:1978}) add to the definition of good quotient the requirement that $\pi$ is affine. The inclusion of this hypothesis is purely technical and is only useful for proving some restrictions to the existence of geometric quotients, such as Proposition 3.24 and Remark 3.25 of \cite{Newstead:1978}.
Indeed, in general this condition can be automatically taken for granted \cite[Propositions 0.7 and 0.8]{MFK:1994}.
\end{rmk}

\begin{ex}\label{ex:GIT-affine}
Let $X = \Spec(R)$ be an affine algebraic variety, where $R$ is a finitely generated reduced $k$-algebra. By considering $R$ as the $k$-algebra of regular functions on $X$, the action of $G$ on $X$ induces an action on $R$. By Nagata's theorem (see \cite{Nagata:1963} or \cite[Theorem 3.4]{Newstead:1978}), if $G$ is reductive (i.e.\ if its radical is isomorphic to a torus group) then the $k$-algebra of $G$-invariant elements of $R$, denoted by $R^G$, is finitely generated. Therefore, the inclusion $R^G \hookrightarrow R$ induces a regular morphism $\pi: X \to \Spec(R^G)$ onto an affine variety. This map is actually a good quotient \cite[Theorem 3.5]{Newstead:1978} and thus defines the so-called \emph{affine GIT quotient}
$$
	\pi: X \to X \sslash G := \Spec(R^G).
$$
\end{ex}

The previous example shows that constructing good quotients for actions of reductive groups on affine varieties is quite easy thanks to Nagata's theorem. Nevertheless, for general algebraic varieties the problem is more involved and some extra pieces of data are needed, even though $G$ is reductive. 

Let us suppose that we are considering the algebraic action of a reductive group $G$ on an algebraic variety $X$. The first problem we must address is that now $\cO(X)$, the ring of regular functions on $X$, does not characterizes $X$ (e.g. $\cO(X) = k$ whenever $X$ is smooth projective). Thus it do<es not make sense to use it to construct the quotient. Nevertheless, if $X \subseteq {\mathbb{P}^N_k}$ then we have a natural substitute, namely considering an extension of the action of $G$ to the ambient space ${\mathbb{A}^{N+1}_k}$ where the affine cone of $X$ lies. In this vein, the action of $G$ on $X$ can be studied through its action on $\cO({\mathbb{A}^{N+1}_k}) = k[x_0, \ldots, x_N]$.

Furthermore, we may suppose that this action on ${\mathbb{A}^{N+1}_k}$ is actually through linear maps. In this situation, we get the notion of a {\emph{linearization}} of the action: a linear representation $G \to \GL{N+1}(k)$ {which}, when restricted to $X \subseteq {\mathbb{P}^N_k}$, gives the original action.
Notice that the same concept can be described intrinsically by means of line bundles. The action of $G$ is said to be \emph{linearizable} if there exists a line bundle $L \to X$ with a fiberwise linear action of $G$ compatible with the one on $X$. Observe that, if $L$ is ample, then some tensor power of $L$ gives an embedding $X \hookrightarrow {\mathbb{P}^N_k}$ for $N$ large enough, recovering the extrinsic notion of {linearization}.

If $L \to X$ is a fixed linearization of the action of $G$, then a point $x \in X$ is called \emph{semi-stable} if there exists a $G$-invariant section $f$ of $L^r$ for some $r>0$ {with} $f(x) \neq 0$ and $X_f = \left\{x \in X \,|\, f(x) \neq 0\right\} \subseteq X$ is affine. If the action of $G$ on $X_f$ is, in addition, closed and $\dim Gx = \dim G$, then $x$ is called \emph{stable}. The set of semi-stable and stable points are open subsets of $X$ and we will denote them by $X^{SS}$ and $X^S$, respectively.

With these definitions, the most important result of GIT about the existence of quotients states that if $G$ is a reductive group acting via a linearizable action on a (quasi-projective) variety $X$, then there exists a good quotient of $X^{SS}$ that restricts to a geometric quotient on $X^S$. It is customary to call this good quotient the \emph{GIT quotient} and denote it by $X^{SS} \to X \sslash G$, or $X^{SS} \sslash G$ when we want to emphasize that it is defined only on $X^{SS}$. The proof of this result is just an appropriate gluing of the good quotients constructed in Example \ref{ex:GIT-affine} for an affine covering of $X^{SS}$ (see \cite[Theorem 3.14]{Newstead:1978}). Despite the fact that, a priori, the result of this gluing is an algebraic scheme, if we push forward the ample line bundle $L|_{X^{SS}} \to X^{SS}$ to $X \sslash G$ we obtain an ample line bundle there that embeds it into projective space, turning the scheme into a variety.

\begin{rmk}\label{rmk:affine-git}
The affine GIT quotient, as described in Example \ref{ex:GIT-affine}, can be seen as a degenerate case of this result. In this situation, if $X \subseteq {\mathbb{A}^N_k}$, then $X$ is naturally embedded into ${\PP^{N}_k}$ through the slice $x_0=1$, where $(x_0: x_1: \ldots, x_N)$ are the homogeneous coordinates in ${\PP^{N}_k}$. The linearization bundle is $\cO(1)|_X \to X$ so that sections of $\cO(1)|_X^k$ are polynomials of degree $k$ (equivalently, homogeneous polynomials of degree $k$ restricted to $x_0=1$). Thus, the action of $G$ is canonically linearizable. Moreover, all the points of $X$ are semi-stable. Indeed, the constant polynomial $f(x) \equiv 1$ (equivalently, the restriction of the homogeneous polynomial $\bar{f}(x) = x_0$) is $G$-invariant, since the action of $G$ is trivial on $x_0$, and $X_f = X$ is affine. Hence, the good quotient exists on the whole $X$, which justifies the notation $X \sslash G$ for the affine GIT quotient.
\end{rmk}

\section{Definition of pseudo-quotients and the question of uniqueness}
\label{section:pseudo-quotients}

Let $X$ be an algebraic variety defined over $k$, which as we recall is an algebraically closed field (not necessarily of characteristic zero), and suppose that an algebraic group $G$ acts on $X$. In this section, we shall introduce a notion of weak quotient for the action of $G$ on $X$ that only captures the expected topological properties. We will call them `pseudo-quotients'. As we will see, since we are getting rid of the algebraic/categorical nature of the quotient, pseudo-quotients are no longer unique. Nevertheless we will show that, in characteristic zero, they are unique enough if we are only concerned with $K$-theory: pseudo-quotients are unique in the Grothendieck ring of algebraic varieties.

\begin{defn}
Let $X$ be an algebraic variety with an action of an algebraic group $G$. A \emph{pseudo-quotient} for the action of $G$ on $X$ is a surjective $G$-invariant regular morphism $\pi: X \to Y$ such that, for any disjoint $G$-invariant closed sets $W_1, W_2 \subseteq X$, we have that $\overline{\pi(W_1)} \cap \overline{\pi(W_2)} = \emptyset$.

\end{defn}

\begin{rmk}\label{rmk:prop-pseudo-quotients}
Suppose that $\pi: X \to Y$ is a pseudo-quotient. Directly from its definition, it satisfies the following:
\begin{enumerate}[label=$\roman*)$,ref=$\roman*)$]
	\item Let $x_1,x_2 \in X$. Since $\pi$ is $G$-invariant, the morphism $\pi$ maps every point of $\overline{Gx_i}$ into the same point of $Y$, namely $\pi(x_i)$. Therefore since the $\overline{Gx_i}$ are closed $G$-invariant sets, we have that $\overline{Gx_1} \cap \overline{Gx_2} = \emptyset$ if and only if $\pi(x_1) \neq \pi(x_2)$.
	\item \label{rmk:prop-pseudo-quotients-closed-sets} 
Let $W \subseteq X$ be a $G$-invariant closed set and suppose that $\pi(W)$ is not closed. Then, for any $y \in \overline{\pi(W)} - \pi(W)$, we would have that $\pi^{-1}(y)$ and $W$ are closed $G$-invariant sets so $\left\{y\right\} \cap \overline{\pi(W)} = \emptyset$, which is impossible. Thus the image of any $G$-invariant closed set is closed. In particular, good quotients are pseudo-quotients.	
	\item \label{rmk:prop-pseudo-quotients-good} Let $U \subseteq Y$ be an open set and let $\pi^*: \cO_Y(U) \to \cO_X(\pi^{-1}(U))$ be the induced ring morphism. Since $\pi$ is $G$-invariant, this morphism factorizes through the inclusion $\cO_X(\pi^{-1}(U))^G \subseteq \cO_X(\pi^{-1}(U))$ so it defines a ring morphism $\pi^*: \cO_Y(U) \to \cO_X(\pi^{-1}(U))^G$. However, in contrast with good quotients, we no longer require this morphism to be an isomorphism. For all these reasons, a pseudo-quotient is a regular map satisfying conditions $i), ii), iv)$ and $v)$ of a good quotient, but maybe failing $iii)$.
	\item \label{rmk:prop-pseudo-quotients-restrictions} If $\pi: X \to Y$ is a pseudo-quotient and $W \subseteq X$ is a closed $G$-invariant set, then the restriction $\pi: W \to \pi(W)$ is also a pseudo-quotient. For open sets, an easy adaptation of Lemma 3.6 of \cite{Newstead:1978} shows that, if $U \subseteq Y$ is an open set, then $\pi: \pi^{-1}(U) \to U$ is a pseudo-quotient.
\end{enumerate}
\end{rmk}

\begin{ex}\label{rmk:uniqueness-pseudo-quotients}
In contrast with categorical quotients, pseudo-quotients may not be unique. As an example, let us take $X = {\mathbb{A}^2_k}$ and $G = k^* = k - \left\{0\right\}$ acting by $\lambda \cdot (x,y)=(\lambda x, \lambda^{-1}y)$, for $\lambda \in k^*$ and $(x,y) \in {\mathbb{A}^2_k}$. Since the ring of $G$-invariant regular functions on $X$ is $\cO_X(X)^G = k[xy]$, standard GIT results show that the inclusion $k[xy] \hookrightarrow k[x,y] =  \cO_X(X)$ induces a good quotient $\pi: X \to {\mathbb{A}^1_k}$. Now, let $C = \left\{y^2 = x^3\right\} \subseteq {\mathbb{A}^2_k}$ be the standard nodal cubic curve. The map $\alpha: {\mathbb{A}^1_k} \to C$ given by $\alpha(t) = (t^2, t^3)$ is a regular bijective morphism, so $\alpha \circ \pi: X \to C$ is a pseudo-quotient for $X$. However, the curve $C$ is not isomorphic to ${\mathbb{A}^1_k}$ since $C$ is not normal.
\end{ex}

The previous example is general in the sense that, if $\pi: X \to Y$ is a pseudo-quotient and $\alpha: Y \to Y'$ is any regular bijective morphism, then $\alpha \circ \pi: X \to Y'$ is also a pseudo-quotient.
We will devote the rest of this section to {studying} the question of uniqueness of pseudo-quotients. We start with the next proposition, which shows that the origin of non-uniqueness of pseudo-quotients arises in essentially the same manner as in the previous example.

\begin{prop}\label{prop:pseudo-quotient-regular-bijective}
Let $X$ be an algebraic variety {acted on} by an algebraic group $G$. Suppose that $\pi: X \to Y$ and $\pi': X \to Y'$ are pseudo-quotients and that $\pi$ is a categorical quotient. Then there exists a regular bijective morphism $\alpha: Y \to Y'$.

\begin{proof}
By definition, the map $\pi': X \to Y'$ is $G$-invariant so, using the categorical property of $Y$, it defines a regular map $\alpha: Y \to Y'$ such that $\pi'=\alpha \circ \pi$. The surjectivity of $\alpha$ follows from {that} of $\pi'$. For the injectivity, suppose that $\alpha(y)=\alpha(y')$ for some $y,y' \in Y$. Then, there exists $x,x' \in X$ such that $y = \pi(x)$ and $y' = \pi(x')$ so $\pi'(x)=\alpha(\pi(x))=\alpha(\pi(x')) = \pi'(x')$. Thus, since $\pi'$ is a pseudo-quotient, we get that $\overline{Gx} \cap \overline{Gx'} \neq \emptyset$ so $y=\pi(x) = \pi(x') = y'$.
\end{proof}
\end{prop}

If there exists a pseudo-quotient $\pi: X \to Y$ that is also a categorical quotient, we will say that $X$ admits a \emph{categorical pseudo-quotient}.

In order to obtain other results about uniqueness we will also need to make, for the rest of the section, the additional assumption that the field $k$, fixed earlier, has characteristic zero. This is because we will require the following results
which are valid only in characteristic zero.

\begin{prop}\label{rmk:properties-regular-bijective:isomorphism}
 If $\alpha: X \to Y$ is a regular bijective morphism and $Y$ is normal, then $\alpha$ is an isomorphism.
\end{prop}

\begin{proof}
In characteristic zero, every dominant injective regular morphism is birational because every field extension is separable and thus the degree of $\alpha$ is the degree of the field extension $K(X)/K(Y)$. Hence, the morphism $\alpha$ is a birational equivalence. By Zariski's main theorem \cite[Proposition 8.12.3]{EGAIV}, the map $\alpha$ factorizes as $X \stackrel{\alpha_1}{\to} Z \stackrel{\alpha_2}{\to} Y$, with $\alpha_1$ an open immersion and $\alpha_2$ finite. Since $\alpha$ is birational, the morphism $\alpha_2$ must be an isomorphism so $X$ is isomorphic to an open subvariety $\alpha_1(X) \subseteq Z \cong Y$. But, since $\alpha$ is surjective, we have that $\alpha_1(X) = W$ and thus $\alpha$ is an isomorphism. 
\end{proof}

\begin{rmk}\label{rmk:arbitrary-field-ramified}
This result is not true if $k$ has positive characteristic, since a regular bijective map can be ramified and thus it may be not an isomorphism.
\end{rmk}

\begin{cor}
If $\pi: X \to Y'$ is a pseudo-quotient with $Y'$ normal and the action on $X$ admits a categorical pseudo-quotient $X \to Y$ then $Y \cong Y'$ and $\pi$ is categorical too.
\end{cor}

In the general case, it may happen that the categorical quotient does not exist. Nevertheless, in characteristic zero it is possible to compare pseudo-quotients. The key point is the following proposition that adapts Proposition 0.2 of \cite{MFK:1994} to the context of pseudo-quotients. Also, {we can} compare it with Proposition \ref{rmk:properties-regular-bijective:isomorphism} above.

\begin{prop}\label{prop:pseudo-quotient-is-good}
Let $\pi: X \to Y$ be a pseudo-quotient for the action of some algebraic group $G$ on $X$. If $X$ is irreducible and $Y$ is normal, then $\pi$ is a good quotient.

\begin{proof}
As mentioned in Remark \ref{rmk:prop-pseudo-quotients} \ref{rmk:prop-pseudo-quotients-good}, it is enough to prove that, for any open set $U \subseteq Y$, the induced ring morphism $\pi^*: \cO_Y(U) \to \cO_X(\pi^{-1}(U))^G \subseteq \cO_X(\pi^{-1}(U))$ is an isomorphism. Recall that, since $\pi$ is $G$-invariant, we automatically get that $\pi^*\left(\cO_Y(U)\right) \subseteq \cO_X(\pi^{-1}(U))^G$. Moreover $\pi^*$ is injective since $\pi$ is surjective.

In order to prove the surjectivity of the morphism $\pi^*$, let $f: \pi^{-1}(U) \to k={\mathbb{A}^1_k}$ be a regular $G$-invariant function. We want to show that $f$ descends to a regular morphism $\tilde{f}: U \to {\mathbb{A}^1_k}$ such that $f= \tilde{f} \circ \pi$. The only possible candidate is the one defined as $\tilde{f}(y)=f(x)$ for any $x \in \pi^{-1}(y)$, which is well-defined since $\pi$ is a pseudo-quotient.
Observe that $\tilde{f}$ is continuous in the Zariski topology. To check that, let $W \subseteq {\mathbb{A}^1_k}$ {be} a closed set. Then the set $f^{-1}(W)$ is a closed $G$-invariant set so, by Remark \ref{rmk:prop-pseudo-quotients} \ref{rmk:prop-pseudo-quotients-closed-sets}, the image $\pi(f^{-1}(W)) \subseteq U$ is closed. But $\pi(f^{-1}(W)) = \tilde{f}^{-1}(W)$ by construction, proving the continuity of $\tilde{f}$.

Therefore, it is enough to prove that $\tilde{f}: U \to {\mathbb{A}^1_k}$ is regular. For this purpose, let us construct the morphism $\pi'= f \times \pi : \pi^{-1}(U) \to {\mathbb{A}^1_k} \times U$ and let $U' \subseteq {\mathbb{A}^1_k} \times U$ be the closure of $\pi'(\pi^{-1}(U))$. Set $p_1: {\mathbb{A}^1_k} \times U \to {\mathbb{A}^1_k}$ and $p_2: {\mathbb{A}^1_k} \times U \to U$ for first and second projections respectively and let $\omega = p_2|_{U'}$. We have a commutative diagram
\[
\begin{displaystyle}
   \xymatrix
   {
  	 \pi^{-1}(U) \ar[rr]^{\pi'}\ar[rd]_{\pi} \ar@/^2.5pc/[rrrr]^{f} && U' \ar@{^{(}-{>}}[r] \ar[ld]^{\omega} & {\mathbb{A}^1_k} \times U \ar[r]^{p_1} & {\mathbb{A}^1_k} \\
  	 & U &&&
   }
\end{displaystyle}   
\]
To finish the proof, it is enough to prove that $\omega$ is an isomorphism since in that case $\tilde{f} = p_1 \circ \omega^{-1}$ would be regular. The map $\omega$ is surjective since $\pi$ is surjective. Moreover $\omega$ is injective on $\pi'(\pi^{-1}(U)) \subseteq U'$ since, for any $x,x' \in \pi^{-1}(U)$, if $\omega(\pi'(x)) = \omega(\pi'(x'))$ then $\pi(x) = \pi(x')$, which happens if and only if $\overline{Gx} \cap \overline{Gx'} \neq \emptyset$. In that case, since $f$ is $G$-invariant, we get that $f(x)=f(x')$ and thus $\pi'(x)=(f(x), \pi(x)) = (f(x'), \pi(x')) = \pi'(x')$.

Therefore, to conclude that $\omega$ is bijective, it is enough to prove that $\pi'$ is surjective. Let $A = U' - \pi'(\pi^{-1}(U))${,} which can be described as the set of points $z \in U'$ such that, for any $x \in \pi^{-1}(U)$ with $\pi(x)=\omega(z)$, {we have} $\pi'(x) \neq z$. Using the continuity of $\tilde{f}$,  we can rewrite $A$ as the open set
$$
	A = \left\{(l, y) \in U' \,|\, \tilde{f}(y) \neq l\right\}.
$$
If $A$ was non-empty then, since $\pi'(\pi^{-1}(U)) \subseteq U'$ is dense, we would get $\pi'(\pi^{-1}(U)) \cap A \neq \emptyset$, which is impossible.

Hence $\pi'$ is surjective and thus $\omega: U' \to U$ is a regular bijective morphism. But, since the characteristic of $k$ is zero, Proposition \ref{rmk:properties-regular-bijective:isomorphism} implies that $\omega$ is an isomorphism.
\end{proof}
\end{prop}

\section{The Grothendieck ring of varieties and quotients}
\label{sec:grothendieck-ring}

A very appropriate setting for unifying the previous results about the uniqueness of pseudo-quotients is the so-called Grothendieck ring of algebraic varieties. In this section we will review briefly its definition and will show how it can be applied to recast the uniqueness of pseudo-quotients. For further information, see \cite{Bass:1968}.

\begin{defn}\label{defn:grothendieck-ring}
Fix a ground field $k$. The \emph{Grothendieck ring of algebraic varieties}, denoted by $\K{\Var{k}}$, also known as the $K$-theory of algebraic varieties, is the { commutative} ring generated by the isomorphism classes of algebraic varieties $X$, denoted by $[X]$, { modulo the so-called `cut-and-paste' relations
$$
	[X] = [Y] + [X - Y],
$$
for any algebraic variety $X$ and any closed subvariety $Y \subseteq X$. The addition and multiplication in the ring $\K{\Var{k}}$ are given by disjoint union and cartesian product of varieties, respectively.}
Given an algebraic variety $X$, { the image of its} isomorphism class $[X] \in \K{\Var{k}}$ is called the \emph{virtual class} of $X$. The virtual class of the affine line will be denoted by $q = [\AA^1_k]$ and is called the \emph{Lefschetz motif}.
\end{defn}

\begin{rmk}
Beware of the similarity of the notation between virtual classes and orbits. The virtual class of $X$ is denoted by $[X]$ and the orbit of $X$ under the action of a group $G$ is denoted by $[X]_G$.
\end{rmk}

If we work in $\K{\Var{k}}$ the question of uniqueness of pseudo-quotients simplifies in many cases. As we have seen in Example \ref{rmk:uniqueness-pseudo-quotients}, pseudo-quotients are not unique in general. However, the following result shows that they are unique in $\K{\Var{k}}$ in characteristic zero, as byproduct of Theorem \ref{prop:pseudo-quotient-is-good}.

\begin{cor}\label{cor:equality-epol-pseudo-quotients}
Suppose that $k$ has characteristic zero. Let $X$ be an algebraic variety with an action of an algebraic group $G$. For any pseudo-quotients $\pi_1: X \to Y_1$ and $\pi_2: X \to Y_2$, the varieties $Y_1$ and $Y_2$ have the same virtual class in $\K{\Var{k}}$.
\begin{proof}
Restricting to the irreducible components of $X$ if necessary, we can suppose that $X$ is irreducible. Let $Y^{'}_i \subseteq Y_i$ be the open subset of normal points of $Y_i$ for $i = 1,2$. Since the preimages $\pi_i^{-1}(Y_i^{'}) \subseteq X$ are saturated open sets, we get that $U = \pi_1^{-1}(Y_1^{'}) \cap \pi_2^{-1}(Y_2^{'}) \subseteq X$ is a saturated open set. Thus the restrictions $\pi_i|_U: U \to \pi_i(U) \subseteq Y_i^{'}$ are pseudo-quotients onto normal varieties so, by Proposition \ref{prop:pseudo-quotient-is-good}, they are good quotients. In particular, they are categorical quotients so $\pi_1(U)$ is isomorphic to $\pi_2(U)$ by uniqueness.

Therefore, working inductively on $X-U$, we find a stratification $Y_1 = Z_1 \sqcup \ldots Z_s$ and $Y_2 = \hat{Z}_1 \sqcup \ldots \hat{Z}_s$ such that $Z_j$ is isomorphic to $\hat{Z}_j$ for $j = 1, \ldots, s$. Hence, they define the same object in the Grothendieck ring of algebraic varieties.
\end{proof}
\end{cor}

\begin{rmk}\label{rmk:properties-regular-bijective:polynomial}
If we assume in addition that a categorical pseudo-quotient exists for the action of $G$ on $X$, then a simpler proof of Corollary \ref{cor:equality-epol-pseudo-quotients} can be given, which does not use Proposition \ref{prop:pseudo-quotient-is-good}. The proof of this claim is straightforward using that in characteristic zero every dominant injective regular morphism is birational (c.f.\  Proposition \ref{rmk:properties-regular-bijective:isomorphism}) and arguing as in Corollary \ref{cor:equality-epol-pseudo-quotients}.
\end{rmk}

\begin{rmk}\label{rmk:arbitrary-field-uniqueness}
Opposed to Remark \ref{rmk:properties-regular-bijective:polynomial}, in positive characteristic it is still an open problem whether there exist algebraic varieties with a regular bijective morphism between them, but different virtual classes in $\K{\Var{k}}$. On the other hand, if we replace $\K{\Var{k}}$ with the Grothendieck semi-ring of algebraic varieties (for which equality is a cut-and-paste matter), such {an} example is known \cite[Theorem 3.2]{Beke2017Feb}.

Additionally, in this setting it is also known that if we replace $\K{\Var{k}}$ with the Grothendieck ring of constructible sets a regular bijective morphism does induce equality of the virtual classes \cite[Lemma 1.2]{Beke2017Feb}. In this way, in positive characteristic pseudo-quotients are unique up to virtual class as constructible sets. Thus, the results of this paper are valid in positive characteristic without further modifications if we work on the Grothendieck ring of constructible sets.
\end{rmk}

\subsection{Some calculations in $\K{\Var{k}}$}\label{sec:some-calculations}

By definition, on $\K{\Var{k}}$ we have that $\coh{X \times Y} = \coh{X} \cdot \coh{Y}$. However, this property also holds with more generality.

\begin{prop}\label{prop:principal-bundle}
Let $\pi: X \to B$ be a regular morphism which is a fiber bundle in the Zariski topology with fiber $F$. Then
$$
	\coh{X} = \coh{F} \cdot \coh{B}.
$$
\begin{proof}
Decompose $B = \bigcup_{i} U_i$ into a finite open cover so that $\pi|_{\pi^{-1}(U_i)}: \pi^{-1}(U_i) \to U_i$ is trivial. Make the cover disjoint by taking $W_1 = U_1$, $W_2 = U_2 - U_1$ and in general $W_i = U_i - \bigcup_{j < i} U_j$. In this way, we have decomposed $B = \bigsqcup_i W_i$, with $W_i$ locally closed sets of $B$ so that $[\pi^{-1}(W_i)] = [F]\cdot[W_i]$. Thus, summing all the elements of this cover we get
$$
	[X] = \sum_i \left[\pi^{-1}(W_i)\right] = \sum_i [F]\cdot[W_i] = [F] \sum_i [W_i] = [F]\cdot [B].
$$
\end{proof}
\end{prop}

{
\begin{ex}
As a direct application of Proposition \ref{prop:principal-bundle} and the cut-and-paste relations of Definition \ref{defn:grothendieck-ring}, we can perform the following computations. Recall that we set $q = [\mathbb{A}^1_k] = [k] \in \K{\Var{k}}$.
\begin{itemize}
	\item $[\PP^n_k] = q^n + q^{n-1} + \ldots + 1$. Observe that we have a natural stratification
	$$
		\PP^n_k = \mathbb{A}^n_k \sqcup \mathbb{A}^{n-1}_k \sqcup \ldots \sqcup \mathbb{A}^{1}_k \sqcup \star,
	$$
	so we have $[\PP^n_k] = [\mathbb{A}^1_k]^n + [\mathbb{A}^1_k]^{n-1} + \ldots [\mathbb{A}^{1}_k] + [\star] = q^n + q^{n-1} + \ldots + 1$.
	\item $[\GL{2}(k)] = q^4 - q^3 - q^2 + q$. To prove it, observe that we have a locally trivial fibration in the Zariski topology $\GL{2}(k) \to k^2 - \left\{(0,0)\right\}$ given by $A \mapsto Ae_1$ where $e_1$ is the first vector of the canonical basis of $k^2$. The fiber of this map is the set of vectors of $k^2$ that do not lie in the line spanned by $Ae_1$. In this way
	$$
		[\GL{2}(k)] = \left[k^2 - \left\{(0,0)\right\}\right]\left[k^2 - k\right] = (q^2-1)(q^2-q).
	$$
	Similar expressions for the virtual class of $[\GL{n}(k)]$ for $n \geq 2$ can be obtained with a similar argument.
	\item $[\SL{2}(k)] = q^3-q$. In this case, we still have the Zariski locally trivial fibration $\SL{2}(k) \to k^2 - \left\{(0,0)\right\}$, $A \mapsto Ae_1$, but now the fiber at $v = (v_1, v_2)$ of this map are the set of vectors $w = (w_1, w_2)$ such that $\det(v, w) = v_1w_2 - v_2w_1 = 1$. Since either $v_1$ or $v_2$ are non-zero, then this fiber is isomorphic to $k$. Thus, we get
	$$
		[\SL{2}(k)] = \left[k^2 - \left\{(0,0)\right\}\right]\left[k\right] = (q^3-q).
	$$
	
	\item $[\PGL{2}(k)] = q^3-q$. To show it, recall that a map in $\PGL{2}(k)$ is given by the choice of a projective reference of $\PP^1_k$, which in this case amounts to give three different points of $\PP^1_k$. Hence, if we denote by $\Delta_{i,j} = \{(x_1, x_2, x_3) \in {(\PP^1_k)^3} \,|\, x_i = x_j\} \cong (\PP^1_k)^2$, and by $\Delta_{1,2,3} = \{(x_1, x_2, x_3) \in {(\PP^1_k)^3} \,|\, x_1 = x_2=x_3\} \cong \PP^1_k$ we have the description
	\begin{align*}
		[\PGL{2}(k)] &= [\PP^1_k]^3 - \left[\Delta_{1,2} - \Delta_{1,2,3}\right] - \left[\Delta_{1,3} - \Delta_{1,2,3}\right] - \left[\Delta_{2,3} - \Delta_{1,2,3}\right] - \left[\Delta_{1,2,3}\right] \\
					&= [\PP^1_k]^3 - 3[\PP^1_k]^2 + 2 [\PP^1_k]= q^3-q.
	\end{align*}
\end{itemize}
\end{ex}
}

{ Finally, let us discuss the effect of finite group actions on virtual classes. Suppose that $X$ is an algebraic variety with an algebraic action of $\ZZ_2$ and set $\coh{X}^+ = \coh{X \sslash \ZZ_2}$ and $\coh{X}^- = \coh{X} - \coh{X \sslash \ZZ_2}$ so that $\coh{X} = \coh{X}^+ + \coh{X}^-$. Given two quasi-projective algebraic varieties $X_1$ and $X_2$ each of them endowed with an algebraic action of $\ZZ_2$, in \cite{GPHV} a formula is proved to compute the virtual class of the quotient $(X_1 \times X_2)/ \ZZ_2$ with respect to the diagonal action of $\ZZ_2$, i.e.\ $-1 \cdot (x_1, x_2) = (-1 \cdot x_1, -1 \cdot x_2)$ for $x_1 \in X_1$ and $x_2 \in X_2$. For this purpose, we shall say that a variety $Y$ is \textit{quasi-linear} if there exists a decomposition $Y = \PP^n_k - L_1 - \ldots - L_r$ for some linear subspaces $L_1, \ldots, L_r \subseteq \PP^n_k$. Similarly, an action of $\ZZ_2$ on $Y$ is said to be \textit{quasi-linear} if it is the restriction of a $\ZZ_2$-action on $\PP^n_k$ that preserves the arrangement $L_1 \cup \ldots \cup L_r$.

\begin{prop}[\cite{GPHV}]\label{prop:pm-formula}
Let $X_1$ and $X_2$ be algebraic varieties equipped with an algebraic action of $\ZZ_2$. Suppose that $X_1$ belong to the subring of $\K{\Var{k}}$ generated by quasi-linear varieties with quasi-linear actions. Then we have
\begin{align}\label{form:pm-formula}
    \coh{X_1 \times X_2}^+ = \coh{X_1}^+\coh{X_2}^+ + \coh{X_1}^-\coh{X_2}^-
\end{align}
for the diagonal action of $\ZZ_2$ on the product $X_1 \times X_2$.
\end{prop}

Two useful variants of formula (\ref{form:pm-formula}) can be easily obtained:
\begin{itemize}
	\item Proceeding recursively with (\ref{form:pm-formula}), we get that
\begin{align}\label{form:Z2-action}
\begin{split}
    \coh{X^n}^+ = \frac{1}{2}\left[\coh{X}^n + \left(\coh{X}^+ - \coh{X}^-\right)^n\right],\quad \coh{X^n}^- = \frac{1}{2}\left[\coh{X}^n - \left(\coh{X}^+ - \coh{X}^-\right)^n\right].
\end{split}
\end{align}
	\item If $X \to X_1$ is a Zariski locally-trivial fiber bundle with fiber $X_2$, all of them equipped with a $\ZZ_2$-action such that the trivializations are equivariant, then by arguing on each trivializating open set and summing up the contributions we get
\begin{align}\label{form:pm-formula-bundle}
    \coh{X}^+ = \coh{X_1}^+\coh{X_2}^+ + \coh{X_1}^-\coh{X_2}^-.
\end{align}
\end{itemize}}

{ In the complex case, { (\ref{form:pm-formula})} was proven for the $E$-polynomial of the variety (see Section \ref{sec:mhs}) in \cite{LMN}, and the general case of any finite group and $E$-polynomials is treated in \cite{Florentino-Silva}.}

{
As an application of formula (\ref{form:pm-formula-bundle}), we obtain the following result.

\begin{lem}\label{lem:PGL2-quotient}
Consider the action of $\ZZ_2$ on $\PGL{2}(k)$ by permutation of columns. Then
$$
	[\PGL{2}(k)]^+ = q^3-q, \qquad [\PGL{2}(k)]^- = 0.
$$
\begin{proof}
The proof is an adaptation of \cite[Proposition 3.3]{LMN} to work in the Grothendieck ring of varieties. Writing an element $A \in \PGL{2}(k)$ by columns as $A = [c_1 | c_2]$, let us consider the action of $k^* = k - \{0\}$ given by $\lambda \cdot [c_1 | c_2] = [c_1 | \lambda c_2]$ for $\lambda \in k^*$. This action defines a $k^*$-fibration
$$
	k^* \longrightarrow \PGL{2}(k) = \GL{2}(k) / k^* \stackrel{\pi}{\longrightarrow} \GL{2}(k) / (k^* \times k^*).
$$
The action of $k^* \times k^*$ on $\GL{2}(k)$ is given by $(\lambda, \mu) \cdot (c_1 | c_2) = (\mu c_1 | \lambda \mu c_2)$ so the quotient $\GL{2}(k) / (k^* \times k^*)$ is isomorphic to $\PP^1_k \times \PP^1_k - \Delta$, where $\Delta$ is the diagonal.

Moreover, this fibration is locally trivial in the Zariski topology. Indeed, if we consider $U_0 = \{[x_0: x_1] \in \PP^1_k \,|\, x_0 \neq 0\}$ and $U_1 = \{[x_0: x_1] \in \PP^1_k \,|\, x_1 \neq 0\}$ the usual affine charts, then $U_i \times U_j - \Delta$ are the trivializing open sets for $i,j \in \{0, 1\}$. For instance, in $U_0 \times U_0 - \Delta$ we have an isomorphism $(U_0 \times U_0 - \Delta) \times k^* \cong \pi^{-1}(U_0 \times U_0 - \Delta)$ given by
\begin{align}\label{eq:triviality-PGL2}
	([1 : x_1], [1 : y_1], \lambda) \mapsto \begin{bmatrix}1 & \lambda \\ x_1 & \lambda y_1\end{bmatrix},
\end{align}
where $\lambda \in k^*$ and $[1 : x_1], [1 : y_1] \in U_0$ are two different points.

Furthermore, if we consider the action of $\ZZ_2$ on $\PP^1_k \times \PP^1_k - \Delta$ by permutation of the factors and on $k^*$ by $\lambda \mapsto \lambda^{-1}$, the map (\ref{eq:triviality-PGL2}) is equivariant. Hence, applying formula (\ref{form:pm-formula-bundle}) we get
\begin{align}\label{form:PGL2-quotient}
	[\PGL{2}(k)]^+ = [\PP^1_k \times \PP^1_k - \Delta]^+[k^*]^+ + [\PP^1_k \times \PP^1_k - \Delta]^-[k^*]^-.
\end{align}
Now, observe that $(\PP^1_k \times \PP^1_k)/\ZZ_2 \cong \PP^2_k$ and that under this identification the diagonal $\Delta$ becomes a plane conic $X \subseteq \PP^2_k$, so we get 
$[\PP^1_k \times \PP^1_k - \Delta]^+ = [\PP^2_k] - [X] = q^2 + q + 1 - (q+1) =q^2$.
On the other hand, the map $k^* \to k$ given by $\lambda \mapsto \lambda + \lambda^{-1}$ provides an isomorphism $k^* / \ZZ_2 \cong k$. Since $[\PP^1_k \times \PP^1_k - \Delta] = (q+1)^2 - (q + 1)$ and $[k^*] = q-1$, plugging these values in (\ref{form:PGL2-quotient}) the result follows.
\end{proof}
\end{lem}
}

\subsection{Hodge theory}
\label{sec:mhs}

In the case that the ground field is $k = \CC$, additionally to the virtual class of a complex algebraic variety $X$, we can also consider other algebro-geometric invariants called Hodge structures. These are linear structures attached to $H^\bullet_c(X; \QQ)$, the rational compactly supported cohomology of $X$. Given a finite dimensional $\QQ$-vector space $H$, a \emph{mixed Hodge structure} on $H$ is a pair of filtrations $W_\bullet$ of $H$ (increasing, called the weight filtration) and $F^\bullet$ of $H_\CC =H \otimes_\QQ \CC$ (decreasing, called the Hodge filtration) satisfying the following additional property: for any $w \in \ZZ$, the induced filtration of $F^\bullet$ on the graded complex $(Gr^W_w\,H)_\CC$ is $w$-orthogonal i.e.\ $F^p \oplus \, \overline{F^{w-p+1}} = (Gr^W_w\,H)_\CC$ for all $p$. For further information, see \cite{Peters-Steenbrink:2008}.

In his celebrated works \cite{DeligneII:1971} and \cite{DeligneIII:1971}, Deligne proved that $H_c^k(X; \QQ)$ is naturally equipped with a mixed Hodge structure. In other words, the cohomology functor $H_c^k(-; \QQ)$ factorizes through the category $\HS$ of mixed Hodge structures
\[
\begin{displaystyle}
   \xymatrix
   {
   	\CVar \ar[rr]^{H_c^k(-;\QQ)} \ar@{--{>}}[d] && \Vect{\QQ} \\
   	\HS \ar[rru] &&
  	}
\end{displaystyle}
\]

\begin{defn}
Given a mixed Hodge structure $H$, we define its Hodge pieces by
$$
	H^{p,q} = \Gr_p^{F}\left(\Gr_W^{p+q}H\right)_\CC.
$$
In addition, if $H = H_c^k(X; \QQ)$ is the compactly supported cohomology of degree $k$ of a complex algebraic variety $X$, the non-negative integers
$$
	h_c^{k;p,q}(X) = \dim_\CC \left(H_c^k(X; \QQ)\right)^{p,q}
$$
are called the \emph{Hodge numbers} of $X$. 
\end{defn}

In particular, given a complex algebraic variety $X$, we can collect the Hodge numbers of its cohomology into a single three-variable polynomial called the (compactly supported) \emph{mixed Hodge polynomial}
$$
	\mu(X)(t, u,v) = \sum_{k,p,q} h_c^{k;p,q}(X)\;t^ku^pv^q \in \ZZ[t, u,v].
$$
{ Notice that, as suggested by the notation, only Hodge numbers $h_c^{k;p,q}(X)$ with positive indices $k, p, q \geq 0$ appear in these Hodge structures.}

We can {also} consider the specialization
$$
	e(X)(u,v) = \mu(X)(-1, u,v) = \sum_{k,p,q} (-1)^k h_c^{k;p,q}(X)\;u^pv^q \in { \ZZ[u,v]},
$$
called the \emph{$E$-polynomial} or the \emph{Deligne-Hodge polynomial} of $X$. Thanks to the long exact sequence in compactly supported cohomology induced by disjoint union and the K\"unneth isomorphism, the $E$-polynomial descends to a ring homomorphism
$$
	e: \K{\Var{\CC}} \to { \ZZ[u,v]}.
$$
In this sense { $e(X) \in \ZZ[u,v]$ }is a weaker invariant than the whole virtual class $[X] \in \K{\Var{\CC}}$. This former invariant is the one that has been studied in the literature for character varieties, for instance in \cite{Hausel-Rodriguez-Villegas:2008,LM,LMN,Munoz-Martinez:2015}.

\section{Properties of pseudo-quotients}
\label{section:properties-pseudo-quotients}

In this section, we will prove several results regarding the flexibility of pseudo-quotients with respect to geometric constructions. In particular, we will show two properties of pseudo-quotients concerning their behavior with respect to stratifications and to `cores', a notion which we will define in this section. This kind of results will be very useful in the upcoming computations regarding character varieties. 

\subsection{Stratifications of algebraic quotients}
\label{section:stratification}
Let us come back for a while to classical GIT and consider the following problem. Suppose that a reductive group $G$ acts on an algebraic variety $X$ and that the action is linearizable. For simplicity, let us also suppose that all the points of $X$ are semi-stable (otherwise, we focus on $X^{SS}$). In that case we have a good quotient $\pi: X \to X \sslash G$. Consider a decomposition $X =  Y \sqcup U$ into a closed set $Y \subseteq X$ and an open set $U \subseteq X$ both saturated for the action. The key issue is that, in general, the GIT quotient is not well-behaved under decompositions in the sense that
$$
	X \sslash G \neq (Y \sslash G) \sqcup (U \sslash G).
$$
The problem is not in $U$ since the restriction $\pi: U \to \pi(U)$ is again a good quotient and thus $\pi(U) = U \sslash G$ \cite[Theorem 3.10]{Newstead:1978}. However, for the closed part the restriction $\pi: Y \to \pi(Y)$ might not be a good quotient: there could exist $G$-invariant functions on $Y$ that do not factorize through $\pi|_{Y}$. This is a rather shocking phenomenon that cannot happen if $G$ is linearly reductive (in particular in characteristic zero) but it shows up in positive characteristic. However, the topological properties of the good quotient remain valid when restricted to $Y$. This suggests that pseudo-quotients are the right framework to study this {stratification} problem.

Except as otherwise stated, the arguments of this section are valid on any algebraically closed field $k$.

\begin{defn}\label{defn:orbitwise-closed}
Let $X$ be a variety and let $G$ be an algebraic group acting on $X$. A subset $A \subseteq X$ is said to be \emph{orbitwise-closed} if $\overline{Ga} \subseteq A$ for any $a \in A$. Analogously $A$ is said to be \emph{completely orbitwise-closed} if both $A$ and $X - A$ are orbitwise-closed.
\end{defn}

\begin{ex}\label{ex:orbitwise-closed} A closed invariant subset is orbitwise-closed. An open orbitwise-closed set is completely orbitwise-closed.
\end{ex}

\begin{lem}\label{lem:completely orbitwise-closed}
Let $G$ be an algebraic group acting on a variety $X$ and let $\pi: X \to Y$ be a pseudo-quotient.
\begin{enumerate}[label=$\roman*)$,ref=$\roman*)$]
	\item\label{lem:completely orbitwise-closed:enum:saturated} If $A \subseteq X$ is completely orbitwise-closed, then it is saturated for $\pi$, i.e.\ {$\pi^{-1}(\pi(A))=A$}.
	\item\label{lem:completely orbitwise-closed:enum:good} If $U$ is open and orbitwise-closed, then $\pi(U)$ is open. Moreover $\pi|_U: U \to \pi(U)$ is a pseudo-quotient.
\end{enumerate}
\begin{proof}
For \ref{lem:completely orbitwise-closed:enum:saturated}, trivially $A \subseteq \pi^{-1}(\pi(A))$. For the other inclusion, if $x \in \pi^{-1}(\pi(A))$ then $\pi(x) = \pi(a)$ for some $a \in A$. Hence, since $\pi$ is a pseudo-quotient, we have $\overline{Gx} \cap \overline{Ga} \neq \emptyset$ which implies that $x \in A$ since $A$ is orbitwise-closed. For \ref{lem:completely orbitwise-closed:enum:good}, by Remark \ref{rmk:prop-pseudo-quotients} \ref{rmk:prop-pseudo-quotients-closed-sets} if we set $W = X - U$, then $\pi(W)$ is closed so $Y - \pi(W)$ is open. But, since $U$ is completely orbitwise-closed and $\pi$ is surjective, we have that $Y = \pi(U) \sqcup \pi(W)$. Thus $\pi(U) = Y - \pi(W)$ is open. Finally, the fact that $\pi|_U: U \to \pi(U)$ is a pseudo-quotient follows from Remark \ref{rmk:prop-pseudo-quotients} \ref{rmk:prop-pseudo-quotients-restrictions} since $U$ is saturated.
\end{proof}
\end{lem}

\begin{rmk}
If $U$ is an open set, by Example \ref{ex:orbitwise-closed} and Lemma \ref{lem:completely orbitwise-closed}, we have that being orbitwise-closed is equivalent to the more familiar property of being saturated for a pseudo-quotient.
\end{rmk}

\begin{thm}\label{prop:decomposition-quotient}
Let $X$ be an algebraic variety with an action of an algebraic group $G$. Suppose that we have a decomposition $X = Y \sqcup U$ where $Y$ is a closed subvariety and $U$ is an open orbitwise-closed subvariety. Then, for any pseudo-quotient $\pi: X \to \overline{X}$, we have a decomposition
$$
	\overline{X} = \pi(Y) \sqcup \pi(U),
$$
where the map $\pi(U) \subseteq \overline{X}$ is open, the map $\pi(Y) \subseteq \overline{X}$ is closed and the morphisms $\pi|_U: U \to \pi(U)$ and $\pi|_Y: Y \to \pi(Y)$ are pseudo-quotients.

Furthermore, if $k$ has characteristic zero then, for any pseudo-quotients $Y \to {\overline{Y}}$ and $U \to {\overline{U}}$, we have that $[\overline{X}] = [\overline{Y}] + [\overline{U}]$ in the Grothendieck ring of algebraic varieties.

\begin{proof}
The decomposition $\overline{X} = \pi(Y) \sqcup \pi(U)$ and properties of $\pi(Y)$ and $\pi(U)$ follow immediately from the surjectivity of $\pi$ together with Remark \ref{rmk:prop-pseudo-quotients} \ref{rmk:prop-pseudo-quotients-restrictions} and Lemma \ref{lem:completely orbitwise-closed} \ref{lem:completely orbitwise-closed:enum:good}. For the last part, use Corollary \ref{cor:equality-epol-pseudo-quotients}.
\end{proof}
\end{thm}

\begin{rmk}\label{rmk:quotient-strat} With the notations of Theorem \ref{prop:decomposition-quotient}, suppose that
$\pi: X \to \overline{X}$ is good.
\begin{enumerate}
	\item For the open stratum, the restriction $\pi|_U: U \to \pi(U)$ is also good \cite[Proposition 3.10]{Newstead:1978}.
	\item A priori, the restriction may not be good for the closed stratum $\pi|_Y: Y \to \pi(Y)$. This holds if $G$ is a linearly reductive group \cite[Remark 3.4.3]{Newstead:1978}. Moreover, if $k$ has characteristic zero, any reductive group is linearly reductive \cite[Corollary 22.43]{milne:2017}, so the restriction of any good quotient to the closed stratum is also good.
	\item In particular, the previous observation provides a shorter proof of the last part of Theorem \ref{prop:decomposition-quotient} when the initial quotient is good. Nevertheless, in this paper we will need the slightly stronger version for pseudo-quotients as stated above. The key point is that the equality of virtual classes $[\overline{X}] = [\overline{Y}] + [\overline{U}]$ holds even if $Y \to \overline{Y}$ is not good and is only a pseudo-quotient.
	\item If $\pi(Y)$ is normal and $k$ has characteristic zero, then by Proposition \ref{prop:pseudo-quotient-is-good} the map $\pi|_Y: Y \to \pi(Y)$ is good. This provides an alternative proof to the fact mentioned in item $(2)$.
\end{enumerate}
\end{rmk}

\begin{ex}
The hypothesis that $\pi: X \to \overline{X}$ is a pseudo-quotient is required in Theorem \ref{prop:decomposition-quotient}, even if $k = \CC$. Consider $G = \CC$ and $X = {\mathbb{A}^2_k}$ with the action $\lambda \cdot (x,y)=(x, y + \lambda x)$, for $\lambda \in \CC$ and $(x,y) \in {\mathbb{A}^2_k}$. Recall that $G$ is not reductive, so classical GIT theory does not guarantee that a good quotient for the action exists. Actually, the map $\pi: X \to \overline{X} = {\mathbb{A}^1_k}$ given by $\pi(x,y)=x$ is a categorical quotient but is not a pseudo-quotient so it is not good.

Let us take $U = \left\{x \neq 0\right\} \subseteq X$ and $Y = \left\{x=0\right\} \subseteq X$. Observe that $U$ is orbitwise-closed (actually, the full action is closed), the restriction $\pi|_U: U \to \overline{U} = \pi(U) = {\mathbb{A}^1_k}-\left\{0\right\}$ is a good quotient and $\pi(Y) = \left\{0\right\}$. On the other hand, the identity map $Y \to Y$ is a categorical pseudo-quotient for the trivial action of $G$ on $Y$ but $[\overline{X}] \neq [Y] + [\overline{U}]$ (this can be checked using the $E$-polynomial, see Section \ref{sec:mhs}) so the conclusion of Theorem \ref{prop:decomposition-quotient} fails.
\end{ex}

{
We finish this section with a brief discussion on how virtual classes behave for geometric quotients. Recall from Proposition \ref{prop:principal-bundle} that if $X \to B$ is a Zariski locally trivial principal $G$-bundle, then $[X] = [G][B]$. Thanks to the celebrated Luna slice theorem  \cite{Luna:1973} (see also \cite[Proposition 5.7]{Drezet:2004}) stated below, geometric quotients of reductive groups are principal bundles in the étale topology, so some kind of multiplicative property might be expected for these quotients, at least for the action of some groups.

\begin{thm}[Luna slice theorem]\label{thm:Luna}
Let $X$ be an affine variety with an action of a reductive algebraic group $G$ on it. Let $X_0 \subseteq X$ be the set of points where the action of $G$ is free and closed (on $X$). Then $X_0$ is open and saturated for the GIT quotient $\pi: X \to X / G$ and the restriction $\pi|_{X_0}: X_0 \to \pi(X_0) \subseteq X / G$ is a principal $G$-bundle in the étale topology.
\end{thm}

This result can be used to provide a useful consequence for our upcoming computations. Recall that an algebraic group $G$ is said to be \emph{special} if any regular morphism $X \to B$ that is a principal $G$-bundle in the \emph{étale topology} is also a principal bundle in the Zariski topology. Hence, for special groups it holds that if $X \to B$ is locally trivial in the étale topology, then $[X] = [G]\cdot [B]$. The groups $\GL{n}(k)$ and $\SL{n}(k)$ are known to be special \cite[Théor\`{e}me 2 and Section 4.4b]{Serre}.

\begin{cor}\label{cor:luna-thm-epol}
With the notations and hypotheses of Theorem \ref{thm:Luna}, if $U \subseteq X_0$ is an orbitwise-closed open set, then $\pi|_U : U \to \pi(U) = U / G$ is a principal $G$-bundle in the étale topology. Moreover, if $G$ is special, we have $\coh{U} = \coh{U / G}\cdot \coh{G}$.

\begin{proof}
By Theorem \ref{thm:Luna}, the map $\pi|_{X_0}: X_0 \to \pi(X_0)$ is a principal $G$-bundle. Since $U$ is an orbitwise-closed open set, Lemma \ref{lem:completely orbitwise-closed} \ref{lem:completely orbitwise-closed:enum:saturated} implies that $U$ is saturated and $\pi(U)$ is open. In this setting, the restriction $\pi|_U: U \to \pi(U)$ is also a principal $G$-bundle. The multiplicativity of virtual classes follows directly from the fact that $G$ is special.
\end{proof}
\end{cor}
}

\subsection{Cores and pseudo-quotients}

Another useful application of pseudo-quotients is presented in Proposition \ref{prop:core} below. The idea is the following. Suppose that we are dealing with a $G$-action on an algebraic variety $X$ and that we manage to find a subvariety $Y \subseteq X$ such that $\overline{Gx} \cap Y \neq \emptyset$ for any $x \in X$. In this case, every point of a good quotient of $X$ by $G$ contains a {representative in} $Y$ so, in order to understand the quotient of $X$, we can just focus on $Y$. On the other hand $\overline{Gx} \cap Y$ might contain more than a single point, so $Y$ is not a slice: we need to quotient $Y$ by the induced action of the subgroup $H \subseteq G$ that keeps $Y$ invariant.

In this situation, it is natural to expect that good quotients for the action of $G$ on $X$ and for the action of $H$ on $Y$, if they exist, should be closely related. However, in general, these quotients might not be isomorphic as algebraic varieties. The following result shows that the right framework to state this problem is as pseudo-quotients and that, under this topological point of view, both pseudo-quotients agree.

\begin{defn}\label{defn:core}
Let $X$ be an algebraic variety with an action of an algebraic group $G$. A \emph{core} for this action is a pair $(Y, H)$ of a subvariety $Y \subseteq X$ and an algebraic subgroup $H \leq G$ such that:

\begin{enumerate}[label=\roman*),ref=$\roman*)$]
	\item\label{prop:core:enum:0} $Y$ is orbitwise-closed for the action of $H$.
	\item\label{prop:core:enum:1} For any $x \in X$, we have $\overline{Gx} \cap Y \neq \emptyset$.
	\item\label{prop:core:enum:2} For any two $W_1, W_2 \subseteq Y$ disjoint closed (in $Y$) $H$-invariant subsets, we have that $\overline{GW_1} \cap \overline{GW_2} = \emptyset$.
\end{enumerate}
\end{defn}

\begin{prop}\label{prop:core}
Let $X$ be an algebraic variety with an action of an algebraic group $G$ and let $(Y, H)$ be a core for this action. If there exists a pseudo-quotient $\pi: X \to \overline{X}$ for the action of $G$ on $X$, then the morphism $\pi$ restricts to a pseudo-quotient $\pi|_{Y}: Y \to \overline{X}$ for the action of $H$ on $Y$.
\begin{proof}
For the surjectivity of $\pi|_Y$, take $\overline{x} \in \overline{X}$. Since $\pi$ is surjective, we have that $\overline{x} = \pi(x)$ for some $x \in X$ and, by hypothesis \ref{prop:core:enum:1} of Definition \ref{defn:core}, there exists $y \in \overline{Gx} \cap Y$. Thus, we get $\pi(y) = \pi(x) = \overline{x}$. 

Now, let $W_1, W_2 \subseteq Y$ be two disjoint closed $H$-invariant subsets. If $\overline{\pi(W_1)} \cap \overline{\pi(W_2)} \neq \emptyset$, then $\pi^{-1}\left(\overline{\pi(W_1)}\right) \cap \pi^{-1}\left(\overline{\pi(W_2)}\right) \neq \emptyset$. But we claim that $\pi^{-1}\left(\overline{\pi(W_i)}\right) = \overline{GW_i}$, so this is impossible by hypothesis \ref{prop:core:enum:2} of Definition \ref{defn:core}. In order to check that, observe that the inclusion $\overline{GW_i} \subseteq \pi^{-1}\left(\overline{\pi(W_i)}\right)$ is trivial. For the other inclusion, since $\overline{GW_i}$ is a closed $G$-invariant set, then $\pi(\overline{GW_i})$ is a closed subset containing $\pi(W_i)$ and thus $\overline{\pi(W_i)} \subseteq \pi(\overline{GW_i})$ which implies that $\pi^{-1}\left(\overline{\pi(W_i)}\right) \subseteq \pi^{-1}\left(\pi(\overline{GW_i})\right) = \overline{GW_i}$.
\end{proof}
\end{prop}

\begin{rmk}
The idea of a core is implicitly presented in the differentiable setting, for instance in gauge theory. To endow the moduli space of connections on a fixed hermitian bundle with the structure of a Hilbert manifold, around a fixed connection $A_0$ it is customary to gauge connections $A$ near $A_0$ to a special `slice' of the action of the gauge group, the so-called Coulomb gauge \cite[Proposition 4.2.9]{donaldson1986geometry}. However, the {representative} of $A$ {in} this `slice' is not unique: the connections in {the Coulomb gauge} must be quotiented by the action of the stabilizer of $A_0$. In this sense, the notion of a core is introduced to capture this geometric situation in the algebraic setting. However, now it is not guaranteed that this restricted quotient is a categorical quotient, as in the differentiable setting, only a pseudo-quotient.
\end{rmk}

\begin{ex}\label{ex:cores-polystable}
Cores also appear naturally in classical GIT. Recall that a point $x \in X$ is said to be \emph{poly-stable} if its orbit is closed in $X^{SS}$. The set of poly-stable points of $X$ is denoted by $X^{PS} \subseteq X$. Since the action of $G$ on $X^{PS}$ is closed, the geometric quotient $X^{PS}/G$ can be formed. In can be proven \cite[Lemma 4.16]{Hoskins} that, for any $x \in X^{SS}$, the closure $\overline{Gx}$ contains a unique poly-stable orbit. For this reason, the natural inclusion $X^{PS} \hookrightarrow X^{SS}$ gives rise to a $G$-invariant map $\pi: X^{PS} \to X^{SS} \sslash G$ which can be {lifted} to a regular bijection (actually a homeomorphism \cite[Corollary 4.17]{Hoskins})
$$
	\alpha: X^{PS}/G {\longrightarrow} X^{SS} \sslash G.
$$
Nevertheless, in general the poly-stable quotient $X^{PS}/G$ is not isomorphic to the GIT quotient $X^{SS} \sslash G$. What is going on here is that $(X^{PS}, G)$ is a core for the action of $G$ on $X^{SS}$, the map $\pi$ is the induced pseudo-quotient map by restriction and $\alpha$ is a regular bijective morphism comparing the two pseudo-quotients (c.f.\ Proposition \ref{prop:pseudo-quotient-regular-bijective}). The aim of cores is to mimic this idea but allowing different groups acting on $X$ and on the subvariety $Y$.
\end{ex}

\section{Representation varieties}
\label{sec:rep-var}

For the sake of completeness, in this section we will review some well-known facts about the geometry of representation varieties. Throughout this section, we shall work in an arbitrary algebraically closed field $k$.

Given a finitely generated group $\Gamma$ and an algebraic group $G$, the set of representations $\Hom(\Gamma, G)$ can be naturally endowed with an algebraic structure as follows. Choose a presentation $\Gamma = \langle \gamma_1, \ldots, \gamma_r\;|\; R_\alpha(\gamma_1, \ldots, \gamma_r)=1\rangle$ of $\Gamma$ with finitely many generators, where $R_\alpha$ are the relations (possibly infinitely many). In that case, we define the injective map $\psi: \Hom(\Gamma, G) \to  G^r$ given by $\psi(\rho)=(\rho(\gamma_1), \ldots, \rho(\gamma_r))$. Moreover, the image of $\psi$ is the algebraic subvariety of $G^r$
$$
	\img{\psi} = \left\{(g_1, \ldots, g_r) \in G^r\;\right|R_\alpha(g_1, \ldots, g_r)=1\left.\right\}.
$$
Hence, we can impose an algebraic structure on $\Hom(\Gamma, G)$ by declaring that $\psi$ is a regular isomorphism over its image. This algebraic structure does not depend on the chosen presentation. With this structure $\Hom(\Gamma,G)$ is called the \emph{representation variety} of $\Gamma$ into $G$ and is denoted by $\Rep{G}(\Gamma)$.

The variety $\Rep{G}(\Gamma)$ has a natural action of $G$ by conjugation, i.e.\ $g \cdot \rho (\gamma) = g\rho(\gamma) g^{-1}$ for $g \in G$, $\rho \in \Rep{G}(\Gamma)$ and $\gamma \in \Gamma$. Recall that two representations $\rho, \rho'$ are said to be isomorphic if and only if $\rho' = g \cdot \rho$ for some $g \in G$. For this reason, if $G$ is reductive, it is interesting to consider the GIT quotient
$$
	\cR_G(\Gamma) = \Rep{G}(\Gamma) \sslash G,
$$
which is usually called the \emph{character variety}.

If $G$ is a reductive linear algebraic group, then as reviewed in Section \ref{sec:review-GIT} (particularly Example \ref{ex:GIT-affine}) the character variety coincides with the affine
variety $\Spec \cO(\Rep{G}(\Gamma))^G$, where $\cO(\Rep{G}(\Gamma))$ is the $k$-algebra of regular functions on $\Rep{G}(\Gamma)$. Note that $\Rep{G}(\Gamma)$ is an affine variety provided that $G$ is a linear algebraic group.
Moreover, for all $\gamma \in \Gamma$, the {character} maps $\chi_\gamma: \Rep{G}(\Gamma) \to k$ given by $\chi_\gamma(\rho) = \tr \rho(\gamma)$ belong to $\cO(\Rep{G}(\Gamma))^G$. Indeed, in \cite{Culler-Shalen} it was proven that if $G$ is reductive, then $\cO(\Rep{G}(\Gamma))^G$ is spanned by the characters, which justifies the name character variety.

Our aim in the remainder of this paper is to study the properties of this quotient and to show how to compute the virtual class of $\cR_G(\Gamma)$ from {that} of $\Rep{G}(\Gamma)$. In particular, in this paper we focus on the following cases:

\begin{itemize}
	\item $\Gamma = F_n$, the free group {on} $n$ generators. In that case, for short we will denote the associated representation variety by $\Xf{n}(G) = \Rep{G}(F_n) = G^n$ or, when the group $G$ is understood, just by $\Xf{n}$. The importance of this case comes from the fact that if $\Gamma$ is any finitely generated group with $n$ generators, then the epimorphism $F_n \to \Gamma$ gives an inclusion $\Rep{G}(\Gamma) \subseteq \Xf{n}(G)$.
	\item $\Gamma = \pi_1(\Sigma_g)$, the fundamental group of the compact orientable surface of genus $g$. In that case, we will denote the associated representation variety by $\Xs{g}(G) = \Rep{G}(\pi_1(\Sigma_g))$ or even just by $\Xs{g}$ when the group $G$ is understood from the context. The standard presentation of $\Gamma$ is
$$
	\pi_1(\Sigma_g) = \left\langle \alpha_1, \beta_1 \ldots, \alpha_{g}, \beta_g\;\;\left|\;\; \prod_{i=1}^g [\alpha_i, \beta_i] = 1 \right.\right\rangle,
$$
where $[\alpha_i,\beta_i] = \alpha_i\beta_i\alpha_i^{-1}\beta_i^{-1}$ is the group commutator. Hence, we have that $\Xs{g}(G) \subseteq \Xf{2g}(G)$. Actually $\Xs{g}(G)$ is given by tuples of $2g$ elements of $G$ satisfying the relation of $\pi_1(\Sigma_g)$, so it is a closed subvariety. The GIT quotient of this representation variety under the conjugacy action will be denoted by $\cR_g(G)$ or just by $\cR_g$.
\end{itemize}

\subsection{Stability of representation varieties}\label{sec:stability-rep-var}

In this section we summarize briefly some results about the stability conditions on the GIT quotient of the character variety $\cR_G(\Gamma) = \Rep{G}(\Gamma) \sslash G$. These are well-known results that can be found in the literature, for instance in \cite{Sikora2012}.

Recall that a linear representation $\rho: \Gamma \to \mathrm{GL}(V)$, where $V$ is a finite dimensional $k$-vector space, is said to be \emph{reducible} if there exists a proper $\Gamma$-invariant subspace of $V$. If $G$ is a linear group, then we can see $G \subseteq \GL{m}(k)$ for $m$ large enough so it also makes sense to talk about reducible $G$-representations. With this definition, we have a decomposition
$$
	\Rep{G}(\Gamma) = \Xred{\mathfrak{X}}_G(\Gamma) \sqcup \Xirred{\mathfrak{X}}_G(\Gamma),
$$
where $\Xred{\mathfrak{X}}_G(\Gamma) \subseteq \Rep{G}(\Gamma)$ is the closed subvariety of reducible representations and $\Xirred{\mathfrak{X}}_G(\Gamma)$ the open set of irreducible ones.

The elements in the center of $G$, denoted by $G^0 \subseteq G$, act trivially by conjugation on $\Rep{G}(\Gamma)$. We can get rid of these irrelevant elements by considering the action of the inner automorphism group $\Inn(G) = G/G^0$. The following result, proven by Sikora \cite[Theorem 30 and Corollary 32]{Sikora2012}, gives an algebraic identification of the stable locus of the representation variety as the open set of irreducible representations.

\begin{thm}\label{cor:stable-locus}
Let $\Gamma$ {be} a finitely generated group and $G$ a reductive linear group. Then the stable points for the action of $\Inn(G)$ on $\Rep{G}(\Gamma)$ by conjugation are precisely the irreducible representations, that is $\mathfrak{X}_{G}(\Gamma)^{S} = \Xirred{\mathfrak{X}}_{G}(\Gamma)$.
\end{thm}

\begin{rmk}
There is a small mismatch between our terminology and the one used in \cite{Sikora2012}. There, the term `properly stable' is used to refer to what we have called here `stable', a nomenclature tracing back to \cite{MFK:1994} (see also \cite{Newstead:1978}, particularly the remark above Lemma 3.13).
\end{rmk}

\begin{prop}\label{prop:action-free-on-irred}
Suppose that $k$ is an algebraically closed field. The action of $\Inn(G)=G/G^0$ on the irreducible representations $\Xirred{\mathfrak{X}_{G}}(\Gamma)$ is free. In addition, the map $\Xirred{\mathfrak{X}_{G}}(\Gamma) \to \Xirred{\mathfrak{X}_{G}}(\Gamma) / \Inn(G)$ is a principal $\Inn(G)$-bundle in the étale topology.
\begin{proof}
Take $m$ large enough so that $G \subseteq \GL{m}(k)$. Suppose that $P \in G \subseteq \GL{m}(k)$ fixes an irreducible representation $\rho \in \Xirred{\mathfrak{X}_{G}}(\Gamma)$. In that case, we have that $P \rho = \rho P$ so $P$ is a $\Gamma$-equivariant linear map. By Schur's lemma, this implies that $P = \lambda \Id$ for some $\lambda \in k^*$ and thus $P \in G^0$. The fact that the quotient map $\Xirred{\mathfrak{X}_{G}}(\Gamma) \to \Xirred{\mathfrak{X}_{G}}(\Gamma) / \Inn(G)$ is a principal bundle follows from Corollary \ref{cor:luna-thm-epol}.
\end{proof}
\end{prop}

\begin{rmk}\label{rmk:triang-irred}
\begin{itemize}
	\item If $G \subseteq \GL{2}(k)$, then a representation $\rho: \Gamma \to G$ is reducible if and only if all the elements of $\rho(\Gamma)$ share a common eigenvector.
	\item By \cite[Proposition 1.10]{Borel-1991}, every affine group is linear so the previous results automatically hold for affine reductive groups.
\end{itemize}
\end{rmk}

To finish this section, let {us} sketch how Theorem \ref{cor:stable-locus} can be proven with a hands-on approach following the spirit of the techniques we will use in Sections \ref{subsec:reducible-rep} and \ref{sec:parabolic-rep}. An important numerical criterion of stability is provided by the so-called Hilbert-Mumford criterion \cite[Theorem 4.9]{Newstead:1978} and \cite[Proposition 2.2]{MFK:1994}. This criterion states that $A \in \mathfrak{X}_{G}(\Gamma)^S$ (resp.\ $A \in \mathfrak{X}_{G}(\Gamma)^{SS}$) if and only if $\mu(A, \lambda) > 0$ (resp.\ $\geq 0$) for all $1$-parameter subgroups $\lambda: k^* \to G$, where $\mu(A, \lambda)$ is the minimum $\alpha \in \ZZ$ such that $\lim\limits_{t \to 0} t^\alpha \lambda(t) \cdot A$ exists.

However, the Hilbert-Mumford criterion only works for projective varieties. Hence, in order to apply it to representation varieties, strictly speaking we need to complete the representation variety to a projective variety. Let us focus on the case $G = \GL{m}(k)$ and consider the projective variety
$$
	\overline{\GL{m}(k)} = \left\{(a, b, A) \in {\PP^{m^2+1}_k}\,\left|\, \det(A)=ab^{m-1} \right.\right\},
$$
where $\det(A)$ denotes the usual homogeneous polynomial of a matrix of order $m$ and we see $b = 0$ as the hyperplane at infinity. Observe that $\overline{\GL{m}(k)} \cap \left\{a,b, \neq 0\right\} = \GL{m}(k)$. Using this variety and a presentation $\Gamma = \langle \gamma_1, \ldots, \gamma_r \,|\, R_\alpha(\gamma_1, \ldots, \gamma_r) = 1\rangle$, we can also consider the completion of the representation variety
$$
	\overline{\mathfrak{X}}_{\GL{m}(k)}(\Gamma) = \left\{(a_1, b_1, A_1, \ldots, a_r, b_r, A_r) \in \left(\overline{\GL{m}(k)}\right)^r\,\left|\, \begin{matrix}R_\alpha(A_1, \ldots, A_r) = \Id \end{matrix} \right.\right\},
$$
where the inverse matrices appearing in $R_\alpha(A_1, \ldots, A_r)$ are computed as transpose {adjoint} matrices and the denominators with the determinants are put in the right-hand side of the equation $R_\alpha(A_1, \ldots, A_r) = \Id$.

The action of $\GL{m}(k)$ on $\mathfrak{X}_{\GL{m}(k)}(\Gamma)$ by conjugation extends naturally to a linear action in $k^{r(m^2+2)}$ by fixing the auxiliary variables $a_i, b_i$. In this way, we can apply the Hilbert-Mumford criterion on $\overline{\mathfrak{X}}_{\GL{m}(k)}(\Gamma)$ as above to detect the stable and semi-stable locus. Since the determinant is invariant under the action, we get that $\mathfrak{X}_{\GL{m}(k)}(\Gamma)^S = \overline{\mathfrak{X}}_{\GL{m}(k)}(\Gamma)^S \cap \mathfrak{X}_{\GL{m}(k)}(\Gamma)$ and $\mathfrak{X}_{\GL{m}(k)}(\Gamma)^{SS} = \overline{\mathfrak{X}}_{\GL{m}(k)}(\Gamma)^{SS} \cap \mathfrak{X}_{\GL{m}(k)}(\Gamma)$.

In this vein, the Hilbert-Mumford criterion can be applied directly to $\mathfrak{X}_{\GL{m}(k)}(\Gamma)$ even though this is not a projective variety.
Let us sketch the calculation. By \cite[Theorem 4.11]{Newstead:1978} up to determinant, which is irrelevant for the conjugacy action, all the $1$-parameter subgroups are conjugate to one of the form
$$
	\lambda(t) = \begin{pmatrix} t^{s_1} & 0 & \cdots & 0 \\ 0 & t^{s_2} & \cdots & 0 \\ \vdots & & \ddots &  \vdots \\ 0 & \ldots && t^{s_m}\end{pmatrix} 
$$
for some $s_1 \geq s_2 \ldots \geq s_m \in \RR$ such that $\sum_i s_i = 0$. The action of these subgroups is as follows. Let $A=(A_1, \ldots, A_r) \in \mathfrak{X}_{\GL{m}(k)}(\Gamma)$ with entries $A_l = \left(a_{i,j}^l\right)$. A straightforward computation shows that the $(i,j)$-entry of $\lambda(t)  A_l\lambda(t)^{-1}$ is $a_{i,j}^l t^{s_i - s_j}$.

Now, suppose that $A$ is reducible with a proper invariant subspace of dimension $k < m$. In that case, maybe after conjugation, the element $A$ can be written as an array of $2 \times 2$-block matrices of sizes $k \times k$ (upper-left block), $k \times (m-k)$ (upper-right), $(m-k) \times k$ (bottom-left) and $(m-k) \times (m-k)$ (bottom-right), where the bottom-left block vanishes. 

In this situation, taking $\lambda(t)$ with $s_i = 1/k$ for $i \leq k$, and $s_i = -1 /(m-k)$ for $i > k$, we have that all the non-zero entries of $\lambda(t) \cdot A = \lambda(t) A \lambda(t)^{-1}$ are multiplied by $t^\beta$ with $\beta \geq 0$. Hence $\lim\limits_{t \to 0} \lambda(t) \cdot A$ exists so $\mu(A, \lambda)=0$ and $A$ cannot be stable. Reciprocally, if $A \in \mathfrak{X}_{\GL{m}(k)}(\Gamma)$ is not stable, the existence of a $1$-parameter subgroup such that $\lim\limits_{t \to 0}\lambda(t) \cdot A$ exists implies that $A$ must be conjugate {to} an array of block matrices with vanishing bottom-left block. This evidences a proper invariant subspace of the representation $A$.
 
It is worth mentioning that the same computation shows that all representations are semi-stable, in perfect agreement with the interpretation of affine GIT quotients as degenerate cases of quasi-projective quotients (Remark \ref{rmk:affine-git}). Indeed, all the elements in the diagonal of $\lambda(t) \cdot A$ or below it are multiplied by $t^\beta$ with $\beta \leq 0$. Since at least one of them {has} to be non-zero, we have that $\lim\limits_{t \to 0} t^{\alpha} \lambda(t) \cdot A$ cannot exist if $\alpha < 0$, so $\mu(A, \lambda) \geq 0$ for all $1$-parameter subgroups.

\section{The $\SL{2}(\k)$-representation variety}
\label{subsec:reducible-rep}

For the rest of the paper, we shall fix { an algebraically closed ground field} $k$ of characteristic zero and we will focus on the algebraic group $G = \SL{2}(k)$. To shorten notation, we will abbreviate $\SL{2} = \SL{2}(k)$ and $\PGL{2}=\PGL{2}(k)$.

The aim of this section is to compute the virtual classes of $\SL{2}$-character varieties for free and surface groups. These classes lie in the subring of $\K{\Var{k}}$ generated by $q = \coh{\mathbb{A}_k^1} \in \K{\Var{k}}$, the virtual class of the affine line (the Lefschetz motif). Therefore, in the case $k = \CC$ these results automatically compute the $E$-polynomials of the character varieties by taking the homomorphism {$e: \K{\Var{\CC}} \to \ZZ[u,v]$} of Section \ref{sec:mhs} which sends $e(q)=uv$. In this sense, these results improve { those} of \cite{Hausel-Rodriguez-Villegas:2008,LMN,MM}.

\subsection{Conjugacy classes in $\SL{2}$-character varieties}

The center of $\SL{2}$ is formed by the matrices $\pm \Id$ so the action of $\SL{2}$ by conjugation descends to an action of $\Inn(\SL{2})=\SL{2}/\left\{\pm \Id\right\} = \PGL{2}$. Given a finitely generated group $\Gamma$, we have a decomposition $\mathfrak{X}_{\SL{2}}(\Gamma) = \Xred{\mathfrak{X}_{\SL{2}}}(\Gamma) \sqcup \Xirred{\mathfrak{X}_{\SL{2}}}(\Gamma)$.
The first stratum of this decomposition is a closed set and the second stratum is an open set with a geometric free action of $\PGL{2}$.

\begin{prop}\label{prop:isotropy-of-matrices}
Let $\Gamma$ be a finitely generated group. Suppose that $A = (A_1, \ldots, A_r) \in \mathfrak{X}_{\SL{2}}(\Gamma)$ has non-trivial stabilizer for the action of $\PGL{2}$. Then $A$ is conjugate to an element of the following two types:
$$
\left(\begin{pmatrix} \pm 1 & \alpha_1 \\ 0 & \pm 1\end{pmatrix}, \ldots, \begin{pmatrix}\pm 1 & \alpha_r \\ 0 & \pm 1\end{pmatrix}\right) \;\;\;\txt{or}\;\;\; \left(\begin{pmatrix}\lambda_1 & 0 \\ 0 & \lambda_1^{-1}\end{pmatrix}, \ldots, \begin{pmatrix}\lambda_r & 0 \\ 0 & \lambda_r^{-1}\end{pmatrix}\right),
$$
for $\alpha_i \in k$ and $\lambda_i \in k^* = k - \left\{0\right\}$.
\begin{proof}
	By Proposition \ref{prop:action-free-on-irred}, the representation $A$ must be reducible. Moreover, by Remark \ref{rmk:triang-irred} $A$ is conjugate to a tuple of upper triangular matrices, namely
	$$
		A = (A_1, \ldots, A_r) = \left(\begin{pmatrix} \lambda_1 & \alpha_1 \\ 0 & \lambda_1^{-1}\end{pmatrix}, \ldots, \begin{pmatrix}\lambda_r & \alpha_r \\ 0 & \lambda_r^{-1}\end{pmatrix}\right).
	$$
If all the eigenvalues $\lambda_i = \pm 1$, then $A$ belongs to the first class. Let us suppose that $\lambda_j \neq \pm 1$ for some $1 \leq j \leq r$. In that case, after conjugating $A$ with
$$
	P_0 = \begin{pmatrix}1 & \frac{\alpha_j}{\lambda_j - \lambda_j^{-1}} \\ 0 & 1\end{pmatrix}
$$ we obtain that the upper right entry of the resulting matrix $P_0 A_j P_0^{-1}$ vanishes and is thus a diagonal matrix. We claim that the remaining matrices also become upper triangular upon conjugation so $A$ belongs to the second class. Otherwise, since $P_0$ preserves upper triangular matrices, there must be another upper triangular non-diagonal matrix $P_0A_kP_0^{-1}$ for $k \neq j$. In that case, if $P \in \SL{2}$ fixes $P_0AP_0^{-1}$, we obtain that $P$ stabilizes both the diagonal matrix $P_0 A_j P_0^{-1}$ and the non-diagonal upper triangular matrix $P_0A_kP_0^{-1}$. But this is only possible if $P = \pm \Id$.
\end{proof}
\end{prop}

Let us denote by $\XP{\mathfrak{X}_{\SL{2}}}(\Gamma)$ the representations that belong to the first type of elements of Proposition \ref{prop:isotropy-of-matrices} and do not belong to the second type and let $\XPh{\mathfrak{X}_{\SL{2}}}(\Gamma) = \SL{2} \cdot \XP{\mathfrak{X}_{\SL{2}}}(\Gamma)$ be their orbits. Analogously, we will denote by $\XD{\mathfrak{X}_{\SL{2}}}(\Gamma)$ the representations that only belong to the second class and by $\XDh{\mathfrak{X}_{\SL{2}}}(\Gamma) = \SL{2} \cdot \XD{\mathfrak{X}_{\SL{2}}}(\Gamma)$ the set of their orbits. Finally, the elements of $\XI{\mathfrak{X}_{\SL{2}}}(\Gamma)$ will be the representations that belong to both classes, that is, those of the form $A = (\pm \Id, \ldots, \pm \Id)$. Observe that the action of $\SL{2}$ on this later stratum is trivial.

We set $\XTilde{\mathfrak{X}_{\SL{2}}}(\Gamma) = \Xred{\mathfrak{X}_{\SL{2}}}(\Gamma) - {\XPh{\mathfrak{X}_{\SL{2}}}}(\Gamma) - {\XDh{\mathfrak{X}_{\SL{2}}}}(\Gamma) - \XI{\mathfrak{X}_{\SL{2}}}(\Gamma)$ for the remaining reducible representations, that is, those with a free action of $\PGL{2}$. This is the set of reducible representations, not completely reducible, such that some matrix {does not have} repeated eigenvalues. Using these pieces, we obtain a stratification
\begin{equation}\label{eq:strat-rep-var}
	\mathfrak{X}_{\SL{2}}(\Gamma) = \Xirred{\mathfrak{X}_{\SL{2}}}(\Gamma) \sqcup {\XPh{\mathfrak{X}_{\SL{2}}}}(\Gamma) \sqcup {\XDh{\mathfrak{X}_{\SL{2}}}}(\Gamma) \sqcup \XI{\mathfrak{X}_{\SL{2}}}(\Gamma) \sqcup \XTilde{\mathfrak{X}_{\SL{2}}}(\Gamma).
\end{equation}
This stratification also decomposes $\mathfrak{X}_{\SL{2}}(\Gamma)$ according to its stabilizers for the action of $\SL{2}$ by conjugation. To be precise, the stabilizers of the points of $\XPh{\mathfrak{X}_{\SL{2}}}(\Gamma)$ and $\XDh{\mathfrak{X}_{\SL{2}}}(\Gamma)$ are respectively conjugate to the subgroups of $\SL{2}$
$$
		\Stab\,J_+ = \Stab\,J_- = \left\{\begin{pmatrix} \pm 1 & \beta \\ 0 & \pm 1\end{pmatrix},\, \beta \in \k\right\}, \hspace{0.8cm} \Stab\,D_\lambda = \left\{\begin{pmatrix} \mu & 0\\ 0 & \mu^{-1} \end{pmatrix},\, \mu \in \k^*\right\}.
$$
Here $J_\pm$ are the Jordan type matrices with eigenvalues $\pm 1$ and $D_\lambda$ is the diagonal matrix with eigenvalues $\lambda^{\pm 1}$, that is
$$
	J_+ = \begin{pmatrix}
	1 & 1\\
	0 & 1\\
\end{pmatrix}, \hspace{1cm} J_- = \begin{pmatrix}
	-1 & 1\\
	0 & -1\\
\end{pmatrix},
\hspace{1cm} D_\lambda = \begin{pmatrix}
	\lambda & 0\\
	0 & \lambda^{-1}\\
\end{pmatrix}.
$$
For the stratum $\XI{\mathfrak{X}_{\SL{2}}}(\Gamma)$, the action is trivial and for the strata $\Xirred{\mathfrak{X}_{\SL{2}}}(\Gamma)$ and $\XTilde{\mathfrak{X}_{\SL{2}}}(\Gamma)$ the action of $\PGL{2}$ is free.

\begin{rmk}\label{rmk-luna-stratification}
The stratification considered in this paper is related to the so-called Luna stratification for ${\mathfrak{X}_{\SL{2}}}(\Gamma)$ (see \cite{Luna:1972,Luna:1973}, also \cite[Section 6.9]{Popov-Vinberg:1989}). This is a natural stratification of an algebraic variety under the action of a reductive group according to the conjugacy class of the stabilizer of each point. The Luna stratification plays an important role in the study of the monodromy of the action, similarly to the approach of \cite{LMN}, since on each stratum of the Luna stratification the quotient map is a locally trivial fibration in the \'etale topology.

In our case, the Luna strata for the action of $\SL{2}$ on the representation variety are $\XPh{\mathfrak{X}_{\SL{2}}}(\Gamma)$, determined by having stabilizer $\Stab\,J_+$; $\XDh{\mathfrak{X}_{\SL{2}}}(\Gamma)$, with stabilizer $\Stab\,D_\lambda$; $\XI{\mathfrak{X}_{\SL{2}}}(\Gamma)$, with stabilizer $\SL{2}$; and $\Xirred{\mathfrak{X}_{\SL{2}}}(\Gamma) \cup \XTilde{\mathfrak{X}_{\SL{2}}}(\Gamma)$ with stabilizer $\left\{\pm \Id\right\}$.
\end{rmk}

Before starting with the calculations, we shall need a final result regarding the quotient of the irreducible locus. Recall from Corollary \ref{cor:luna-thm-epol} that we have an étale principal bundle
$$
	\PGL{2} \to \Xirred{\mathfrak{X}_{\SL{2}}}(\Gamma) \to \Xirred{\mathfrak{X}_{\SL{2}}}(\Gamma) / \SL{2}.
$$

Since $\PGL{2}$ is {not} special this is not enough to assure that in $\K{\Var{k}}$ the virtual class of the total space can be written as $\coh{\Xirred{\mathfrak{X}_{\SL{2}}}(\Gamma)} = \coh{\PGL{2}} \coh{\Xirred{\mathfrak{X}_{\SL{2}}}(\Gamma) / \SL{2}}$ (c.f.\ Section \ref{sec:some-calculations}). Nevertheless, at least in this rank $2$ case the multiplicativity of the virtual class can be obtained by means of an auxiliary argument.

\begin{prop}\label{prop:multiplicativity}
For any finitely generated group $\Gamma$, we have
$$
	\coh{\Xirred{\mathfrak{X}_{\SL{2}}}(\Gamma)} = \coh{\PGL{2}} \coh{\Xirred{\mathfrak{X}_{\SL{2}}}(\Gamma) / \SL{2}}.
$$

\begin{proof}
Virtual classes are additive for stratifications so we can address this problem locally in the Zariski topology. {For fixed $i$ ($1 \leq i \leq r$)}, consider the subvariety $X_i \subseteq \Xirred{\mathfrak{X}_{\SL{2}}}(\Gamma)$ of tuples $(A_1, \ldots, A_r) \in \Xirred{\mathfrak{X}_{\SL{2}}}(\Gamma)$ such that the trace $\tr A_i \neq \pm 2$. We form the auxiliary variety
$$
	\tilde{X}_i = \left\{(A_1, \ldots, A_r, \lambda) \in X_i \times \left(k^*-\{\pm 1\}\right) \,|\, \tr A_i = \lambda + \lambda^{-1}\right\}.
$$

In other words, the set $\tilde{X}_i$ is the variety whose points are elements of $X_i$ together with a chosen eigenvalue of $A_i$. This actually allows us to explicitly identify both eigenvalues of $A_i$ since one is $\lambda$ and the other one is $\lambda^{-1}$ (and they are different since $\tr A_i \neq \pm 2$). The group $\PGL{2}$ also acts on $\tilde{X}_i$ by conjugation on the first component, so we get a $\PGL{2}$-principal bundle
\begin{equation}\label{eq:fibration-Xtilde}
	\PGL{2} \to \tilde{X}_i \to \tilde{X}_i/\PGL{2}.
\end{equation}

Let us show that (\ref{eq:fibration-Xtilde}) is actually locally trivial in the Zariski topology. Given a point $A^0 = (A_1^0, \ldots, A_r^0, \lambda^0) \in \tilde{X}_i$, there exists a Zariski open neighborhood $\tilde{U} \subseteq \tilde{X}_i$ of $A^0$ where we can compute algebraically the two (different) eigenspaces of $A_i$ for any $(A_1, \ldots, A_r, \lambda) \in \tilde{U}$. Let us denote these eigenspaces by $V_\lambda(A_i)$ and $V_{\lambda^{-1}}(A_i)$. The reason is that, once we know the eigenvalues of $A_i$, we can obtain the eigenspaces as $V_\lambda(A_i) = \ker(A_i-\lambda I)$ and $V_{\lambda^{-1}}(A_i)= \ker (A_i-\lambda^{-1} I)$.
This amounts to solve two linear systems of equations, a problem that can be solved on a Zariski open set explicitly and algebraically. { In other words, this provides two algebraic line bundles $V_\lambda \to \tilde{U}$ and $V_{\lambda^{-1}} \to \tilde{U}$. Shrinking $\tilde{U}$ even more if needed, we can find non-vanishing sections $v_\lambda: \tilde{U} \to V_{\lambda}$ and $v_{\lambda^{-1}}: \tilde{U} \to V_{\lambda^{-1}}$ of these line bundles, i.e. an algebraic choice of eigenvectors $v_\lambda(A_i)$ and $v_{\lambda^{-1}}(A_i)$ for any $(A_1, \ldots, A_r, \lambda) \in \tilde{U}$.}

Using these data, we can trivialize (\ref{eq:fibration-Xtilde}) on $\tilde{U}$. Given two linearly independent vectors $v_1, v_2 \in k^2$, denote by $M(v_1, v_2) \in \PGL{2}$ the unique projective mapping induced by the linear map sending the standard basis of $k^2$ into $v_1$ and $v_2$. Then, we get a regular isomorphism
$$
	\tilde{\varphi}: \tilde{U} \to \tilde{U}/\PGL{2} \times \PGL{2}, \quad \tilde{\varphi}(A, \lambda) = ((\PGL{2} \cdot A, \lambda), M(v_\lambda(A_i), v_{\lambda^{-1}}(A_i))),
$$
{ where $A = (A_1, \ldots, A_r) \in X_i$.}
Now, we can recover $X_i$ from $\tilde{X}_i$ as follows. There is an action of $\ZZ_2$ on $\tilde{X}_i$ by $-1 \cdot (A, \lambda) = (A,\lambda^{-1})$ and taking the quotient by the action we get $X_i = \tilde{X}_i / \ZZ_2$. Given $\tilde{U}$ trivializing (\ref{eq:fibration-Xtilde}), let $U = \tilde{U}/\ZZ_2  \subseteq X_i$. {Under the isomorphism $\tilde{\varphi}: \tilde{U} \to \tilde{U}/\PGL{2} \times \PGL{2}$, the group $\ZZ_2$ acts on the first component by $-1 \cdot (\PGL{2} \cdot A, \lambda) = (\PGL{2} \cdot A, \lambda^{-1})$ and on the second component by permutation of columns. Hence, using that $[\PGL{2}]^+ = [\PGL{2} / \ZZ_2] = [\PGL{2}]$ and $[\PGL{2}]^- = 0$ by Lemma \ref{lem:PGL2-quotient}, we get that
$$
	[U] = [\tilde{U}]^+ = \left[\tilde{U}/ \PGL{2} \times \PGL{2}\right]^+ = \left[\tilde{U}/ \PGL{2}\right]^+ \left[\PGL{2}\right]^+ =  \left[U/\PGL{2}\right] \left[\PGL{2}\right].
$$
Furthermore, adding up the contributions from all the trivializing open sets $U$, this also proves the equality $[X_i] = [X_i/\PGL{2}][\PGL{2}]$.}

Working on the different open sets $X_i \subseteq  \Xirred{\mathfrak{X}_{\SL{2}}}(\Gamma)$ we can prove the multiplicativity of virtual classes on the open set of irreducible representations with, at least, a matrix with trace different to $\pm 2$. For the remaining closed set $Y = \Xirred{\mathfrak{X}_{\SL{2}}}(\Gamma) - \bigcup_i X_i$, notice that any representation $A = (A_1, \ldots, A_r) \in Y$ must have at least two Jordan type matrices, say $A_i$ and $A_j$ (otherwise, the representation is reducible since all the matrices share an eigenvector). Therefore, the argument above can be repeated on $Y$ but instead of considering $V_\lambda(A_i), V_{\lambda^{-1}}(A_i)$ by considering the two eigenspaces of $A_i$ and $A_j$ (which again must be different since $A$ is irreducible). Putting together the multiplicativity property of virtual classes for each stratum, the results follows.
\end{proof}
\end{prop}

\begin{rmk}
For the last part of the previous argument, in which we choose two Jordan type matrices, we need that $r \geq 2$, i.e.\ $\Gamma$ needs to be generated by at least two elements. But if $\Gamma$ is cyclic, then any representation is automatically reducible so $\Xirred{\mathfrak{X}_{\SL{2}}}(\Gamma) = \emptyset$.
\end{rmk}

Working {in} rank $2$ was crucial in the previous proof since we only needed to identify two eigenspaces. To our knowledge, nothing is known about this issue for $\SL{n}$-character varieties with $n \geq 3$. It is an interesting problem to determine whether this multiplicative property still holds in the higher rank setting.

 Notice that, due to Proposition \ref{prop:multiplicativity}, in order to compute the virtual class $\coh{\Xirred{\mathfrak{X}_{\SL{2}}}(\Gamma) / \SL{2}}$ we will need to artificially invert the element $[\PGL{2}] = q(q-1)(q+1) \in \K{\Var{k}}$. For this purpose, let $S \subseteq \K{\Var{k}}$ be the multiplicative system generated by $q, q-1, q+1$, where recall that $q = \coh{\mathbb{A}_k^1} \in \K{\Var{k}}$ is the Lefschetz motif. We then consider the localization
$$
	\Kh{\Var{k}} = S^{-1}\K{\Var{k}}.
$$
Observe that, since $q$ is a zero divisor in $\K{\Var{k}}$, the natural map $\K{\Var{k}} \to \Kh{\Var{k}}$ is no longer injective: the kernel {is made up of} the annihilators of $q, q-1$ or $q+1$ (or a multiplicative combination of them). However, notice that this localization is not very restrictive. In \cite{behrend512640motive}, Behrend and Dhillon used a completion of this ring to formulate their conjectural motivic formula for the virtual class of the moduli stack of principal bundles over a curve. Furthermore, in \cite[Theorem 1.2]{ekedahl2009grothendieck}, Ekedahl proves that the localization of $\K{\Var{k}}$ by the multiplicative set generated by $q$ and $q^n-1$ for $n \geq 1$, which in particular is a further localization of $\Kh{\Var{k}}$, is isomorphic to the Grothendieck ring of algebraic stacks. In this sense, $\Kh{\Var{k}}$ is a natural ring to consider for moduli problems.

\begin{rmk}
When $k = \CC$, we still have the $E$-polynomial homomorphism $e: \Kh{\Var{\CC}} \to \ZZ[u,v]$. Moreover, since $\ZZ[u,v]$ is an integral domain, the $E$-polynomial of any annihilator $x$ of $q, q+1$ or $q-1$ is $e(x) = 0$. Indeed, the classical examples (e.g.\ \cite{Borisov}) of an annihilator of $q$ are constructed from non-isomorphic varieties $X$ and $Y$ such that $X \times \mathbb{A}^1_k \cong Y \times \mathbb{A}^1_k$. In this case, even though $X$ is not isomorphic to $Y$, their $E$-polynomials must agree.
\end{rmk}

From now on, we shall always work in the localization $\Kh{\Var{k}}$, so we are allowed to divide by $q$, $q-1$ and $q+1$.

\color{black}

\subsection{Free groups}\label{sec:free-groups}

Let us fix $\Gamma = F_n$, the free group of $n$ generators, and recall that set $\Xf{n} = \mathfrak{X}_{\SL{2}}(F_n)$. In this case, each of the pieces of the stratification (\ref{eq:strat-rep-var}) are as follows.

\begin{itemize}
	\item $\XPh{\Xf{n}}$. Given $A = (A_1, \ldots, A_n) \in \XPh{\Xf{n}}$, let
$$
    \left(\begin{pmatrix} \epsilon_1 & \alpha_1 \\ 0 & \epsilon_1\end{pmatrix}, \ldots, \begin{pmatrix}\epsilon_n & \alpha_n \\ 0 & \epsilon_n\end{pmatrix}\right)
$$ be the element of $\XP{\Xf{n}}$ conjugate to $A$ with $\epsilon_i = \pm 1$ and $\alpha_i \in \k$ not all vanishing. This special element of $\XP{\Xf{n}}$ is unique up to simultaneous rescaling of the upper triangular components $\alpha_i$. Thus the $\SL{2}$-orbit of $A$, denoted by $[A]_{\SL{2}}$, is the set of reducible representations $(B_1, \ldots, B_n) \in \Xf{n}$ such that in their upper triangular form
$$
    \left(\begin{pmatrix} \epsilon_1 & \beta_1 \\ 0 & \epsilon_1\end{pmatrix}, \ldots, \begin{pmatrix} \epsilon_n & \beta_n \\ 0 & \epsilon_n\end{pmatrix}\right)
$$
there exists $\lambda \neq 0$ satisfying $(\alpha_1, \ldots, \alpha_n) = \lambda(\beta_1, \ldots, \beta_n)$. Then, taking $\lambda \to 0$ we find that the closure of the orbit $\overline{[A]}_{\SL{2}}$ is precisely the set of reducible representations with double eigenvalues such that their upper triangular components satisfy $(\alpha_1, \ldots, \alpha_n) = \lambda(\beta_1, \ldots, \beta_n)$ for some $\lambda \in \k$. In particular, for $\lambda = 0$ we get that $(\epsilon_1\Id, \ldots, \epsilon_n\Id) \in \overline{[A]}_{\SL{2}}$.

	The virtual class of this stratum can be computed as follows. The tuple $(\alpha_1, \ldots, \alpha_n) \in \k^n-\left\{0\right\}$ determines the upper triangular form up to scaling. Hence, we obtain a regular fibration
	$$
		\k^* \longrightarrow \SL{2}/\Stab\,J_+ \times \left\{\pm 1\right\}^n \times \left(\k^n-\left\{0\right\}\right) \longrightarrow \XPh{\Xf{n}}
	$$
given by
$$
	\left(P, (\epsilon_1, \ldots, \epsilon_n), (\alpha_1, \ldots, \alpha_n)\right) \mapsto \left(P\begin{pmatrix} \epsilon_1 & \alpha_1 \\ 0 & \epsilon_1\end{pmatrix}P^{-1}, \ldots, P\begin{pmatrix}\epsilon_n & \alpha_n \\ 0 & \epsilon_n\end{pmatrix}P^{-1}\right).
$$

This fibration is locally trivial in the Zariski topology so, by Proposition \ref{prop:multiplicativity}, we get that
\begin{align*}
    \coh{\XPh{\Xf{n}}} &= \coh{\left\{\pm 1\right\}^n} \coh{\PP^{n-1}_k} \coh{\SL{2}/\Stab\,} \\&= 2^n(q^2-1)\frac{q^{n} -1}{q-1}.
\end{align*}
Recall that $q = \coh{\AA_\k^1} = \coh{k}$ denotes the virtual class of the affine line.
	\item $\XDh{\Xf{n}}$. Given $A = (A_1, \ldots, A_n) \in \XDh{\Xf{n}}$, let
$$
    \left(\begin{pmatrix} \lambda_1 & 0 \\ 0 & \lambda_1^{-1}\end{pmatrix}, \ldots, \begin{pmatrix}\lambda_n & 0 \\ 0 & \lambda_n^{-1}\end{pmatrix}\right)
$$ be the element of $\XD{\Xf{n}}$ conjugate to $A$ with $\lambda_i \in \k^*$ and not all equal to $\pm 1$.
	The diagonal form of an element of $\XDh{\Xf{n}}$ is determined up to permutation of the columns, so we have a double covering
	$$
		\SL{2}/\Stab\,D_\lambda \times \left((\k^*)^n-\left\{(\pm 1, \ldots, \pm 1)\right\}\right) \longrightarrow \XDh{\Xf{n}}.
	$$
Therefore, we obtain that
$$
	\XDh{\Xf{n}} = \frac{\SL{2}/\Stab\,D_\lambda \times \left((\k^*)^n-\left\{(\pm 1, \ldots, \pm 1)\right\}\right)}{\ZZ_2}.
$$
Using { formula (\ref{form:Z2-action})} and the fact that $[\SL{2}/\Stab\,D_\lambda]^+ = {q^2}$ (c.f.\ \cite[Proposition 3.2]{LMN}), its image in the Grothendieck ring of algebraic varieties is
$$
	\coh{\XDh{\Xf{n}}} = \frac{q^3-q}{2} \left((q-1)^{n-1} + (q+1)^{n-1}\right) - 2^nq^2.
$$

	\item $\XI{\Xf{n}}$. Here the action of $\SL{2}$ is trivial so in particular the action is closed.
	This stratum consists of $2^n$ points, so $\coh{\XI{\Xf{n}}} = 2^n$.
	\item $\XTilde{\Xf{n}}$. This is the set of reducible representations, not completely reducible, such that there is a matrix which does not have any repeated
eigenvalues. Given $A \in \XTilde{\Xf{n}}$, the tuple of matrices $A$ is conjugate to an element of the form
	$$
		\left(\begin{pmatrix} \lambda_1 & \alpha_1 \\ 0 & \lambda_1^{-1}\end{pmatrix}, \ldots, \begin{pmatrix}\lambda_n & \alpha_n \\ 0 & \lambda_n^{-1}\end{pmatrix}\right),
	$$
where the vector $(\alpha_1, \ldots, \alpha_n)$ is determined up to the action of $U = \k^* \times \k$ by
$$
	(\mu, a)\cdot (\alpha_1, \ldots, \alpha_n) = \left(\alpha_1\mu^2 - a\mu\left(\lambda_1-\lambda_1^{-1}\right), \ldots, \alpha_n\mu^2 - a\mu\left(\lambda_n-\lambda_n^{-1}\right)\right).
$$
Since $A \not\in \Xf{n}^D$, we must have $(\mu, a)\cdot(\alpha_1, \ldots, \alpha_n) \neq (0, \ldots, 0)$ for all $(\mu,a) \in U$. This vanishing takes place if and only if $(\alpha_1, \ldots, \alpha_n) = \lambda(\lambda_1 - \lambda_1^{-1}, \ldots, \lambda_n - \lambda_n^{-1})$ for some $\lambda \in \k$. Therefore, the allowed values for $(\alpha_1, \ldots, \alpha_n)$ lie in $\k^n - \ell$, where $\ell \subseteq \k^n$ is the line spanned by $(\lambda_1 - \lambda_1^{-1}, \ldots, \lambda_n - \lambda_n^{-1})$. Thus, we have regular fibration
$$
	U \longrightarrow \PGL{2} \times \Omega
 \stackrel{\pi}{\longrightarrow} \XTilde{\Xf{n}},
$$ 
where {$\Omega \cong \left((\k^*)^n - \left\{\left(\pm 1, \ldots, \pm 1\right)\right\}\right) \times \left(\k^n - \ell\right)$} is the set of allowed values for the eigenvalues and the anti-diagonal components.
Hence, the virtual class of this stratum is
$$
	\coh{\XTilde{\Xf{n}}} = \frac{q^3-q}{(q-1)q} \left((q-1)^n - 2^n\right)\left(q^n-q\right).
$$
\end{itemize}

From this analysis it is possible to describe the virtual class of the GIT quotient $\Xf{n} \sslash \SL{2}$ explicitly. First, we have a decomposition $\Xf{n} = \Xred{\Xf{n}} \sqcup \Xirred{\Xf{n}}$ with $\Xirred{\Xf{n}}$ an open subvariety. Observe that $\SL{2}$ is a reductive affine group, and its action on $\Xirred{\Xf{n}}$ is closed since this is the stable locus (Theorem \ref{cor:stable-locus}). Hence, by the results of Section \ref{section:stratification}, we have
$$
	 \coh{\Xf{n} \sslash \SL{2}} = \coh{\Xred{\Xf{n}} \sslash \SL{2}} + \coh{\Xirred{\Xf{n}} \sslash \SL{2}}.
$$
In addition, the action of $\PGL{2}$ on $\Xirred{\Xf{n}}$ is closed and free so by Corollary \ref{cor:luna-thm-epol} the map $\Xirred{\Xf{n}} \to \Xirred{\Xf{n}} \sslash \PGL{2}$ is a principal $\PGL{2}$-bundle in the étale topology. 
Moreover, by Proposition \ref{prop:multiplicativity} we get that
	$$
		\coh{\Xirred{\Xf{n}} \sslash \SL{2}} = \frac{\coh{\Xirred{\Xf{n}}}}{\coh{\PGL{2}}} = \frac{\coh{\Xirred{\Xf{n}}}}{q^3-q}.
	$$
The computation of $\coh{\Xirred{\Xf{n}}}$ can be done using the analysis above since $\Xirred{\Xf{n}} = \Xf{n} - \XPh{\Xf{n}} - \XDh{\Xf{n}} - \XI{\Xf{n}} - \XTilde{\Xf{n}}$. Using that the virtual class of the whole representation variety is $\coh{\Xf{n}} = \coh{\SL{2}}^n = (q^3-q)^n$, we find that
\begin{align*}
	\coh{\Xirred{\Xf{n}}} =&\, 2^{n} q^{2} - \frac{1}{2} \, {\left(q^{3} - q\right)} {\left({\left(q + 1\right)}^{n - 1} + {\left(q - 1\right)}^{n - 1}\right)} \\
	&- {\left(2^{n} q + {\left(q - 1\right)}^{n} q^{n} - {\left(q - 1\right)}^{n} q - 2^{n}\right)} {\left(q + 1\right)} - 2^{n} + {\left(q^{3} - q\right)}^{n}.
\end{align*}

For $\Xred{\Xf{n}}$ the situation is more involved since the action is very far from being closed. Actually, we are going to show that $(\XD{\Xf{n}}, \ZZ_2)$ is a core for the action of $\SL{2}$ on $\Xred{\Xf{n}}$. For that, we need a preliminary result about the behavior of the orbits.

\begin{prop}\label{prop:intersection-closure-core}
Let $W \subseteq \XD{\Xf{n}}(\SL{2})$ be a closed set. For any $A \in \overline{[W]}_{\SL{2}}$ we have that $\overline{[A]}_{\SL{2}} \cap W \neq \emptyset$.
\begin{proof}

Since $W$ is closed, it is the zero set of some regular functions $f_1, \ldots, f_r: \XD{\Xf{n}} \cong (\k^*)^n \to \k$. Since $\cO_{(\k^*)^n}((\k^*)^n) = \k[\lambda_1, \lambda_1^{-1}, \ldots, \lambda_n, \lambda_n^{-1}]$, after clearing denominators we can suppose that $f_i \in \k[\lambda_1, \ldots, \lambda_n]$. Now, consider the algebraic variety
$$
    \Omega = \left\{(A_i, \lambda_i, v) \in \Xf{n} \times (\k^*)^n \times \left(\k^2-\left\{0\right\}\right)\,\left|\\\;\begin{matrix}A_iv = \lambda_i v\\f_i(\lambda_1, \ldots, \lambda_n) = 0\end{matrix}\right.\right\} \subseteq \Xf{n} \times \k^{n + 2}.
$$
Let $p: \Xf{n} \times \k^{n + 2} \to \Xf{n}$ be the projection onto the first factor and set $W' = p(W)$. If $A = (A_1, \ldots, A_n) \in W'$, then the matrices $A_i$ have a common eigenvector so $A$ is reducible. Hence,
$$
	QAQ^{-1} = \left(\begin{pmatrix}
				\lambda_1 & a_1 \\
				0 & \lambda_1^{-1}
				\end{pmatrix},
				\ldots,
				\begin{pmatrix}
				\lambda_n & a_n \\
				0 & \lambda_n^{-1}
				\end{pmatrix}
				\right)
$$
for some $Q \in \SL{2}$ satisfying $f_i(\lambda_1, \ldots, \lambda_n)=0$. Therefore, taking $P_m =  \begin{pmatrix}
	m^{-1} & 0\\
	0 & m\\
\end{pmatrix}$ we have
$$
	(P_mQ) \cdot A \to \left(\begin{pmatrix}
	\lambda_1 & 0\\
	0 & \lambda_1^{-1}\\
\end{pmatrix}, \ldots, \begin{pmatrix}
	\lambda_n & 0\\
	0 & \lambda_n^{-1}\\
\end{pmatrix}
\right) \in \XD{\Xf{n}},
$$
for $m \to \infty$. This diagonal element is thus in the Zariski closure of $[A]_{\SL{2}}$. Moreover $f_i(\lambda_1, \ldots, \lambda_n)=0$ so this element belongs to $W$ and hence $\overline{[A]}_{\SL{2}} \cap W \neq \emptyset$. Therefore, in order to finish the proof it is enough to show that $\overline{[W]}_{\SL{2}} \subseteq W'$. Trivially we have that $[W]_{\SL{2}} \subseteq W'$ so it is enough to prove that $W'$ is closed.

To do so, let us consider the projectivization of $\Omega$ as the closed projective set
$$
    \tilde{\Omega} = \left\{\left(A_i, [\mu_0:\mu_1: \ldots:\mu_n], \bar{v}\right) \in \Xf{n} \times {\PP^n_k \times \PP^1_k}\,\left|\\\;\begin{matrix} A_i\mu_0v = \mu_i v\\ \tilde{f}_i(\mu_0, \ldots, \mu_n) = 0\end{matrix}\right.\right\} \subseteq \Xf{n} \times {\PP^{n}_k \times \PP^1_k},
$$
where $\tilde{f}_i$ denotes the homogenization of the polynomial $f_i$ with projective coordinates $[\mu_0: \ldots: \mu_n]$. Any element $\left(A_i, [\mu_0:\mu_1: \ldots:\mu_n], \bar{v}\right) \in \tilde{\Omega}$ must have $\mu_0 \neq 0$ so there are no points at infinity. Hence, we can also write $W' = \rho(\tilde{\Omega})$ where  $\rho: \Xf{n} \times {\PP^{n}_k \times \PP^1_k} \to \Xf{n}$ is the first projection. Projective spaces are universally closed, so $\rho$ is closed and thus $W'=\rho(\tilde{\Omega})$ is closed as we wanted to prove.
\end{proof}
\end{prop}

\begin{cor}\label{cor:core-diagonal}
Consider the action of $\ZZ_2$ on $\XD{\Xf{n}}$ by permutation of the eigenvalues. Then $(\XD{\Xf{n}}, \ZZ_2)$ is a core for the action of $\SL{2}$ on $\Xred{\Xf{n}}$. 
\begin{proof}
The orbit of an element of $\XD{\Xf{n}}$ by the action of $\ZZ_2$ is finite, so the action is automatically closed and, thus, part \ref{prop:core:enum:0} of Definition \ref{defn:core} holds. For part \ref{prop:core:enum:1}, given $A \in \Xred{\Xf{n}}$, taking $Q$ and $P_m$ as in the proof of Proposition \ref{prop:intersection-closure-core}, we get that $\overline{[A]}_{\SL{2}} \cap \XD{\Xf{n}} \neq \emptyset$. This intersection consists of at most two points uniquely determined by the eigenvalues of $A$ so in particular they are $\ZZ_2$-equivalent.

In order to prove condition \ref{prop:core:enum:2} of Definition \ref{defn:core}, let $W_1, W_2 \subseteq \XD{\Xf{n}}$ be two $\ZZ_2$-invariant disjoint closed subsets and suppose that $A \in \overline{[W_1]}_{\SL{2}} \cap \overline{[W_2]}_{\SL{2}}$. By Proposition \ref{prop:intersection-closure-core}, we have $\overline{[A]}_{\SL{2}} \cap W_1 \neq \emptyset$ and $\overline{[A]}_{\SL{2}} \cap W_2 \neq \emptyset$. However, this cannot hold since the elements of $\overline{[A]}_{\SL{2}} \cap \XD{\Xf{n}}$ are $\ZZ_2$-equivalent and $W_1, W_2$ are $\ZZ_2$-invariant and disjoint. Thus $\overline{[W_1]}_{\SL{2}} \cap \overline{[W_2]}_{\SL{2}} = \emptyset$, proving condition \ref{prop:core:enum:2}. 
\end{proof}
\end{cor}

\begin{rmk}
In analogy with the results of Section \ref{sec:stability-rep-var}, it can be shown that the poly-stable points of the representation variety for the conjugacy action are exactly the semi-simple representations (i.e.\ direct sums of irreducible representations). In our rank $2$ case, these semi-simple representations are either irreducible representations or diagonal representations (corresponding to the sum of two $1$-dimensional irreducible representations). In our notation, this means $\Xf{n}^{PS} = \Xirred{\Xf{n}} \sqcup \XDh{\Xf{n}}$. In particular, if we focus on the poly-stable points of the reducible locus, we have that $\left(\Xred{\Xf{n}}\right)^{PS} = \XDh{\Xf{n}}$. As mentioned in Example \ref{ex:cores-polystable}, we get a regular bijection
$$
	\left(\Xred{\Xf{n}}\right)^{PS}/\SL{2} = \XDh{\Xf{n}}/\SL{2} \stackrel{\cong}{\longrightarrow} \Xred{\Xf{n}} \sslash \SL{2},
$$
so in particular $\coh{\XDh{\Xf{n}}/\SL{2}} = \coh{\Xred{\Xf{n}} \sslash \SL{2}}$. Under this interpretation, the fact that $(\XD{\Xf{n}}, \ZZ_2)$ is a core can be also understood by observing that $\XD{\Xf{n}} \subseteq \XDh{\Xf{n}}$ plays the role of a `twofold slice' for the $\SL{2}$-orbit space, and the redundancy is killed by the $\ZZ_2$ action.
\end{rmk}

\begin{cor}\label{cor:hodge-red-quot} The virtual class of the GIT quotient of the reducible stratum is
$$
	\coh{\Xred{\Xf{n}} \sslash \SL{2}} = \frac{1}{2} \left((q-1)^{n} + (q+1)^{n}\right).
$$
\begin{proof}
The GIT quotient of $\XD{\Xf{n}}$ by $\ZZ_2$ exists and is a pseudo-quotient. Hence, by Proposition \ref{prop:core} we have that $\coh{\Xred{\Xf{n}} \sslash \SL{2}} = \coh{\XD{\Xf{n}} \sslash \ZZ_2}$. For computing this last term, notice that $\XD{\Xf{n}} = (\k^*)^n$ {and use formula (\ref{form:Z2-action}) (see also \cite{GPHV})}.
\end{proof}
\end{cor}

Once we have computed the virtual classes of the two strata so, summing the contributions, we obtain that
\begin{empheq}{align*}
	\coh{\Xf{n} \sslash \SL{2}} = \frac{1}{2} \, {\left(q + 1\right)}^{n - 1} q + \frac{1}{2} \, {\left(q
- 1\right)}^{n - 1} q - {\left(q - 1\right)}^{n - 1} q^{n - 1} +
{\left(q^{3} - q\right)}^{n - 1}.
\end{empheq}

\begin{rmk}
If we take the homomorphism {$e: \Kh{\Var{\CC}} \to \ZZ[u, v]$}, this result agrees with the computations of \cite{Lawton-Munoz:2016} and \cite{Cavazos-Lawton:2014}. In fact, our computation follows the lines of the calculations of the former paper, generalizing them to virtual classes. In particular, we used the same stratification as in \cite{Lawton-Munoz:2016}. This is not a coincidence. As mentioned in Remark \ref{rmk-luna-stratification}, this stratification is the Luna stratification of $\Xred{\Xf{n}}$, which is the natural one to consider if we want to stratify it according to the stabilizers of the conjugacy action. The notion of a core, introduced in this paper, plays a crucial role to extend the existing results regarding $E$-polynomials to virtual classes.
\end{rmk}

\begin{rmk}
Using the techniques of this section it is also possible to study the case $\Gamma = \ZZ^n$, the free abelian group with $n$ generators, and $G = \SL{m}$. Observe that, in this case, the variety $\Rep{\SL{m}}(\ZZ^n)$ is the set of tuples of $n$ pairwise commuting matrices of $\SL{m}$. Since commuting matrices share a common eigenvector, all the representations of $\Rep{\SL{m}}(\ZZ^n)$ are reducible so $\Rep{\SL{m}}(\ZZ^n) = \Xred{\Rep{\SL{m}}}(\ZZ^n)$. Hence, analogously to Corollary \ref{cor:core-diagonal}, we have that $(\XD{\Rep{\SL{m}}}(\ZZ^n), S_m)$ is a core for $\Xred{\Rep{\SL{m}}}(\ZZ^n)$, where $S_m$ acts on $\XD{\Rep{\SL{m}}}(\ZZ^n) = (\k^*)^{n(m-1)}$ by permutation of the eigenvalues. Therefore, we obtain that
$$
	\coh{\Rep{\SL{m}}(\ZZ^n) \sslash \SL{m}} = \coh{(\k^*)^{n(m-1)} \sslash S_m}.
$$
Analogously, for $G = \GL{m}(\k)$, we obtain the equality of virtual classes $\coh{\Rep{\GL{m}(\k)}(\ZZ^n) \sslash \GL{m}(\k)} = \coh{(\k^*)^{nm} \sslash S_m}$.
This reproves Theorem 5.1 of \cite{Florentino-Silva}.
In the case $m=2$, the virtual classes of these character varieties can be computed by means { of formula (\ref{form:Z2-action})}. In this way, we obtain that
\begin{align*}
	\coh{\Rep{\SL{2}}(\ZZ^n) \sslash \SL{2}} &= \coh{(\k^*)^n \sslash \ZZ_2} = \frac{1}{2} \left((q-1)^{n} + (q+1)^{n}\right),\\
	\coh{\Rep{\GL{2}}(\ZZ^n) \sslash \GL{2}} &= \coh{(\k^*)^{2n} \sslash \ZZ_2} = \frac{1}{2} \left((q-1)^{2n} + (q+1)^{2n}\right).
\end{align*}
In the higher rank case, we need to use stronger results about equivariant cohomology in order to compute the corresponding quotients under $S_m$. This is the strategy accomplished in \cite{Florentino-Silva}.
\end{rmk}

\subsection{Surface groups}

In this section we shall consider the case $\Gamma = \pi_1(\Sigma_g)$, where $\Sigma_g$ is the closed orientable surface of genus $g \geq 1$. Recall that, for short, we denote the associated $\SL{2}$-representation variety by $\Xs{g} = \mathfrak{X}_{\SL{2}}(\pi_1(\Sigma_g))$ and its GIT quotient by $\cR_g = \Xs{g} \sslash \SL{2}$. 

We have that $\Xs{g} \subseteq \Xf{2g}$ as a closed subvariety.
To identify the elements of $\Xs{g}$, let us denote the set of upper triangular matrices of $\Xf{n}$ by $\XU{\Xf{n}} \cong (\k^*)^n \times \k^n$ and set $\XU{\Xs{g}} = \XU{\Xf{2g}} \cap \Xs{g}$. Given $A \in \XU{\Xf{2g}}$, say
$$
	A= \left(
	\begin{pmatrix}
		\lambda_1 & \alpha_1 \\
		0 & \lambda_1^{-1}
	\end{pmatrix},
	\begin{pmatrix}
		\mu_1 & \beta_1 \\
		0 & \mu_1^{-1}
	\end{pmatrix}, \ldots,
	\begin{pmatrix}
		\lambda_{g} & \alpha_{g} \\
		0 & \lambda_{g}^{-1}
	\end{pmatrix},
	\begin{pmatrix}
		\mu_{g} & \beta_{g} \\
		0 & \mu_{g}^{-1}
	\end{pmatrix}
	\right)
$$
with $\lambda_i, \mu_i \in \k^*$ and $\alpha_i, \beta_i \in \k$, a straightforward computation shows that $A \in \XU{\Xs{g}}$ if and only if
\begin{equation*}\label{eq:cond-upper}
	{\sum_{i=1}^g \lambda_i\mu_i \left(\left(\lambda_i-\lambda_i^{-1}\right)\beta_i - \left(\mu_i-\mu_i^{-1}\right)\alpha_i\right)= 0.}
\end{equation*}
This implies that for the strata given by Proposition \ref{prop:isotropy-of-matrices} we get $\XP{\Xs{g}} = \XP{\Xf{2g}}$, $\XD{\Xs{g}} = \XD{\Xf{2g}}$ and $\XI{\Xs{g}} = \XI{\Xf{2g}}$. Let $\pi \subseteq \k^{2g}$ be the $(\alpha_i, \beta_i)$-plane defined by the previous equation for fixed $(\lambda_i, \mu_i)$. The analysis of the stratification is as follows.

\begin{itemize}
	\item For $\Xred{\Xs{g}}$, using the equality $\XD{\Xs{g}} = \XD{\Xf{2g}}$ and Corollary \ref{cor:core-diagonal}, we obtain that $(\XD{\Xs{g}}, \ZZ_2)$ is a core for the action. Therefore, since $\XD{\Xf{2g}} = (\k^*)^n$, we have
$$
	\coh{\Xred{\Xs{g}} \sslash \SL{2}} = \coh{(\k^*)^{2g} \sslash \ZZ_2} = \frac{1}{2} \left((q-1)^{2g} + (q+1)^{2g}\right).
$$
	\item For $\Xirred{\Xs{g}}$, since $\Xirred{\Xf{2g}}$ is an open set of $\Xf{2g}$ where the action of $\PGL{2}$ is closed and free and $\Xs{g} \subseteq \Xf{2g}$ is closed, then $\Xirred{\Xs{g}} = \Xirred{\Xf{2g}} \cap \Xs{g}$ is an open subset of $\Xs{g}$ with a closed and free action. Therefore, Proposition \ref{prop:multiplicativity} shows that
	$$
		\coh{\Xirred{\Xs{g}} \sslash \SL{2}} = \frac{\coh{\Xirred{\Xs{g}}}}{\coh{\PGL{2}}} = \frac{\coh{\Xirred{\Xs{g}}}}{q^3-q}.
	$$
In order to complete the calculation, it is enough to compute the virtual class $\coh{\Xirred{\Xs{g}}}$. For this purpose, we use that $\Xirred{\Xs{g}} = \Xs{g} - \XPh{\Xs{g}} - \XDh{\Xs{g}} - \XI{\Xs{g}} - \XTilde{\Xs{g}}$. Let us compute the virtual class of each stratum as follows.

\begin{itemize}
	\item $\XPh{\Xs{g}}$. In this case, since $\XP{\Xs{g}} = \XP{\Xf{2g}}$ and $\Xs{g}$ is $\SL{2}$-invariant, we have that $\XPh{\Xs{g}} = \XPh{\Xf{2g}}$. Therefore,
$$
    \coh{\XPh{\Xs{g}}} = \coh{\XPh{\Xf{2g}}} = 2^{2g}(q^2-1)\frac{q^{2g} -1}{q-1}.
$$
	\item $\XDh{\Xs{g}}$. Again $\XD{\Xs{g}} = \XD{\Xf{2g}}$ and thus $\XDh{\Xs{g}} = \XDh{\Xf{2g}}$. Therefore,
$$
	\coh{\XDh{\Xs{g}}} = \coh{\XDh{\Xf{2g}}} = \frac{q^3-q}{2} \left((q-1)^{2g-1} + (q+1)^{2g-1}\right) - 2^{2g}q^2.
$$

	\item $\XI{\Xs{g}}$. Again $\XI{\Xs{g}} = \XI{\Xf{2g}}$, which are $2^{2g}$ elements so $\coh{\XI{\Xs{g}}} = 2^{2g}$.
	\item $\XTilde{\Xs{g}}$. In this case, any element is conjugate to one of the form
$$
	\left(
	\begin{pmatrix}
		\lambda_1 & \alpha_1 \\
		0 & \lambda_1^{-1}
	\end{pmatrix},
	\begin{pmatrix}
		\mu_1 & \beta_1 \\
		0 & \mu_1^{-1}
	\end{pmatrix}, \ldots,
	\begin{pmatrix}
		\lambda_{g} & \alpha_{g} \\
		0 & \lambda_{g}^{-1}
	\end{pmatrix},
	\begin{pmatrix}
		\mu_{g} & \beta_{g} \\
		0 & \mu_{g}^{-1}
	\end{pmatrix}
	\right),
$$
with $(\lambda_1, \mu_1, \ldots, \lambda_g, \mu_g) \in (\k^*)^{2g} - \left\{(\pm 1, \ldots, \pm 1)\right\}$ and $(\alpha_1, \beta_1, \ldots, \alpha_{g}, \beta_g) \in \pi - \ell$ where $\ell$ is the line spanned by $(\lambda_1 - \lambda_1^{-1}, \mu_1 - \mu_1^{-1}, \ldots, \lambda_g - \lambda_g^{-1}, \mu_g - \mu_g^{-1})$. Thus, we have a fibration
$$
	U \longrightarrow \PGL{2} \times \Omega \longrightarrow \XTilde{\Xs{g}},
$$
where $U = \k \times \k^*$ and $\Omega$ is a Zariski locally trivial fibration $(\pi - \ell) \to \Omega \to (\k^*)^{2g} - \left\{(\pm 1, \ldots, \pm 1)\right\}$. Using that $\coh{\pi} = q^{2g-1}$ and $\coh{\ell} = q$, the virtual class is
$$
	\coh{\XTilde{\Xs{g}}} = \frac{q^3-q}{(q-1)q} \left((q-1)^{2g} - 2^{2g}\right)\left(q^{2g-1}-q\right).
$$
\end{itemize}
Putting all these data together, we get that
$$
	\coh{\Xred{\Xs{g}}} = (q+1) (q-1)^{2g}\left(q^{2g-1}-q\right) + \frac{q^3-q}{2} \left((q-1)^{2g-1} + (q+1)^{2g-1}\right) - 2^{2g}(q^2-1).
$$
\end{itemize}
From \cite[Section 5.4]{GP-2018} (see also \cite[Proposition 11]{MM} for the formula with $E$-polynomials), we know that the virtual class of the total space is
\begin{align*}
	\coh{\Xs{g}} =&\, 2^{2g - 1} {\left(q - 1\right)}^{2g - 1} {\left(q + 1\right)}
q^{2g - 1} + 2^{2g - 1} {\left(q + 1\right)}^{2g - 1}
{\left(q - 1\right)} q^{2g - 1} \\
&+ \frac{1}{2} \, {\left(q +
1\right)}^{2g - 1} {\left(q - 1\right)}^{2} q^{2g - 1} +
\frac{1}{2} \, {\left(q - 1\right)}^{2g - 1} {\left(q + 1\right)}
{\left(q - 3\right)} q^{2g - 1} \\&+ {{\left(q + q^{2 \,
g-1}\right)} {\left(q^2 - 1\right)}^{2 \,
g - 1}}.
\end{align*}
Hence, we finally obtain that
\begin{empheq}{align*}
	\coh{\cR_g} =&\, \frac{1}{2} \, {\left({\left(2^{2g} + 2 \, {\left(q - 1\right)}^{2g - 2} + q - 1\right)} q^{2g - 2} + q^{2} + 2 \, {\left(q -1\right)}^{2g - 2} + q\right)} {\left(q + 1\right)}^{2g - 2} \\
	&+\frac{1}{2} \, {\left({\left(2^{2g} - 1\right)} {\left(q -
1\right)}^{2g - 2} - {\left(q - 1\right)}^{2g - 2} q - 2^{2g
+ 1}\right)} q^{2g - 2} + \frac{1}{2} \, {\left(q - 1\right)}^{2 \,
g - 1} q.
\end{empheq}

\begin{rmk}
Taking {$e: \Kh{\Var{\CC}} \to \ZZ[u, v]$}, this result recovers the calculations of \cite[Theorem 14]{MM} and \cite[Theorem 1.3]{Baraglia-Hekmati:2016} for the $E$-polynomial.
\end{rmk}

\section{Parabolic representation varieties}
\label{sec:parabolic-rep}

In this section we will discuss a more general setting for representation varieties by considering parabolic structures. Adding parabolic information endows the character variety with an extra structure that allows us to extend the non-abelian Hodge correspondence. For instance, character varieties with generic parabolic semi-simple structures are diffeomorphic to moduli spaces of logarithmic flat connections and of parabolic Higgs bundles \cite{Simpson:parabolic}. Additionally, in the rank $2$ case and two punctures they are diffeomorphic to the moduli space of doubly periodic instantons \cite{Biquard-Jardim, Jardim}.

The definition of a parabolic structure that we will consider here is slightly more general than the one used in the literature, which typically focuses on surface groups. Let $\Gamma$ be a finitely generated group and let $G$ be an algebraic group. A \emph{parabolic structure} $Q$ is a finite set of pairs $(\gamma, \lambda)$ where $\gamma \in \Gamma$ and $\lambda \subseteq G$ is a locally closed subset which is closed under conjugation. Given a parabolic structure $Q$, we define the \emph{parabolic representation variety} $\Rep{G}(\Gamma, Q)$ as the subset of $\Rep{G}(\Gamma)$
$$
	\Rep{G}(\Gamma, Q) = \left\{\rho \in \Rep{G}(\Gamma)\,|\, \rho(\gamma) \in \lambda\,\,\textrm{for all } (\gamma, \lambda) \in Q\right\}.
$$
As in the non-parabolic case, the set $\Rep{G}(\Gamma, Q)$ can be endowed with the structure of an algebraic variety as follows. Suppose that $Q = \left\{(\gamma_1, \lambda_1), \ldots, (\gamma_s, \lambda_s)\right\}$. Choose a finite set of generators $S$ of $\Gamma$ that contains all the $\gamma_i$, namely $S = \left\{\eta_1, \ldots, \eta_r, \gamma_1, \ldots, \gamma_s\right\}$. In that case, using $S$ to identify $\Rep{G}(\Gamma)$ with a closed subvariety of $G^{r+s}$ (see Section \ref{sec:rep-var}), then we also have a natural identification $\Rep{G}(\Gamma, Q) = \Rep{G}(\Gamma) \cap \left(G^r \times \lambda_1 \times \ldots \times \lambda_s\right)$. We impose {on} $\Rep{G}(\Gamma, Q)$ the algebraic structure inherited from this identification.

The adjoint action of $G$ on $\Rep{G}(\Gamma)$ restricts to an action on $\Rep{G}(\Gamma, Q)$ since the subsets $\lambda_i$ are closed under conjugation. Moreover, as a subvariety of $\Rep{G}(\Gamma)$, all the previous results about GIT stability automatically hold in the parabolic setting. Thereby, all the points of $\Rep{G}(\Gamma, Q)$ are semi-stable for this action and the stable locus is $\Xirred{\Rep{G}}(\Gamma, Q) = \Xirred{\Rep{G}}(\Gamma) \cap\Rep{G}(\Gamma, Q)$, where the action of $\Inn(G)=G/G^0$ is closed and free.

In this paper, we are going to focus on the following cases:
\begin{itemize}
	\item Fix some elements $h_1, \ldots, h_s \in G$ and let $[h_i]_G$ be their conjugacy classes. Then, we take $\Gamma = F_{n+s}$ and $Q = \left\{(\gamma_{1}, [h_1]_G), \ldots, (\gamma_s, [h_s]_G)\right\}$, where $\gamma_1, \ldots, \gamma_s \in F_{r+s}$ is an independent set. Observe that $[h_i]_G \subseteq G$ is locally closed since, by \cite{Newstead:1978} Lemma 3.7, it is an open subset of {$\overline{[h_i]}_G$}.
	\item Let $\Sigma = \Sigma_g - \left\{p_1, \ldots, p_s\right\}$ with $p_i \in \Sigma_g$ distinct points, usually referred to as the punctures or the marked points. In this case we have a presentation of the fundamental group of $\Sigma$ given by
$$
	\pi_1(\Sigma) = \left\langle \alpha_1, \beta_1 \ldots, \alpha_{g}, \beta_g, \gamma_1, \ldots, \gamma_s\;\;\left|\;\; \prod_{i=1}^g [\alpha_i, \beta_i]\prod_{j=1}^s \gamma_s = 1 \right.\right\rangle,
$$
where the $\gamma_i$ are the positive oriented simple loops around the punctures. As parabolic structure, we take $Q = \left\{(\gamma_1, [h_1]_G), \ldots, (\gamma_s, [h_s]_G)\right\}$ for some fixed elements $h_i \in G$. Observe that the epimorphism $F_{2g +s} \to \pi_1(\Sigma)$ induces an inclusion $\Rep{G}(\pi_1(\Sigma), Q) \subseteq \Rep{G}(F_{2g+s}, Q)$.
\end{itemize}

As in the previous section, we will focus on the case that $k$ is an algebraically closed field of characteristic zero and $G = \SL{2}(k)$. Again, to shorten notation we will write $\SL{2} = \SL{2}(k)$. In this section, we shall compute the virtual classes of some parabolic $\SL{2}$-character varieties. Recall that, in the case $k = \CC$, these results automatically compute the $E$-polynomials of the character varieties by taking the homomorphism {$e: \Kh{\Var{\k}} \to \ZZ[u,v]$} of Section \ref{sec:mhs} that sends $e(q)=uv$. In this sense, these results extend and improve the ones of \cite{LM} (which are only provided for genus $g=1,2$).

In particular, in this section we are going to study the parabolic character varieties whose parabolic data lie in the conjugacy classes of the matrices
$$
	J_+ = \begin{pmatrix}
	1 & 1\\
	0 & 1\\
\end{pmatrix}, \hspace{1cm} J_- = \begin{pmatrix}
	-1 & 1\\
	0 & -1\\
\end{pmatrix},
\hspace{1cm} -\Id = \begin{pmatrix}
	-1 & 0\\
	0 & -1\\
\end{pmatrix}.
$$
As we will see, the most important case to be considered is to take the parabolic structure $Q_s^+ = \left\{(\gamma_1, [J_+]_{\SL{2}}), \ldots, (\gamma_s, [J_+]_{\SL{2}})\right\}$. The remaining cases can be easily obtained from this one, as shown in Section \ref{subsec:parabolic-data-gen}.

\subsection{Free group}
Set $\Xfp{n}{s} = \mathfrak{X}_{\SL{2}}(F_{n+s}, Q_s^+)$. We have a stratification $\Xfp{n}{s} = \Xirred{\Xfp{n}{s}} \sqcup \Xred{\Xfp{n}{s}}$. Hence, by the results of Section \ref{section:stratification}, we have that
$$
	\coh{\Xfp{n}{s} \sslash \SL{2}} = \coh{\Xirred{\Xfp{n}{s}} \sslash \SL{2}} + \coh{\Xred{\Xfp{n}{s}} \sslash \SL{2}}.
$$
Let us analyze each stratum separately.

\begin{itemize}
	\item $\Xirred{\Xfp{n}{s}}$ is an open subvariety and the action of $\PGL{2}$ induces multiplicativity of the virtual classes for the quotient by a proof analogous to Proposition \ref{prop:multiplicativity}. Thus
	$$
		\coh{\Xirred{\Xfp{n}{s}} \sslash \SL{2}} = \frac{\coh{\Xirred{\Xfp{n}{s}}}}{\coh{\PGL{2}}} = \frac{\coh{\Xirred{\Xfp{n}{s}}}}{q^3-q},
	$$
where {$q = \coh{\AA_k^1} = \coh{k} \in \Kh{\Var{k}}$ is the virtual class of the affine line in the Grothendieck ring of algebraic varieties, localized by $q, q+1,q-1$}. In order to compute $\coh{\Xirred{\Xfp{n}{s}}}$ we count:

\begin{itemize}
	\item $\XPh{\Xfp{n}{s}}$. In this case, the variety $\XP{\Xfp{n}{s}} = \left\{\pm 1\right\}^n \times \left(\k^{n} \times (\k^*)^s\right)$ and the action of $k^*$ restricts to {scaling} the second factor. Therefore,
$$
    \coh{\XPh{\Xfp{n}{s}}} = 2^{n}(q^2-1)\frac{q^{n}(q-1)^s}{q-1}.
$$
	
	\item $\XDh{\Xfp{n}{s}}$. In this case $\XD{\Xfp{n}{s}} = \emptyset$ since $J_+$ is not diagonalizable, so $\XDh{\Xfp{n}{s}} = \emptyset$ and this stratum makes no contribution.

	\item $\XI{\Xfp{n}{s}}$. Again, this stratum is empty so it makes no contribution.
	\item $\XTilde{\Xfp{n}{s}}$. In this case, any element is conjugate to one of the form
$$
	\left(
	\begin{pmatrix}
		\lambda_1 & \alpha_1 \\
		0 & \lambda_1^{-1}
	\end{pmatrix}, \ldots,
	\begin{pmatrix}
		\lambda_n & \alpha_n \\
		0 & \lambda_n^{-1}
	\end{pmatrix},
	\begin{pmatrix}
		1 & c_{1} \\
		0 & 1
	\end{pmatrix}, \ldots,
	\begin{pmatrix}
		1 & c_{s} \\
		0 & 1
	\end{pmatrix}
	\right)
$$
with $(\lambda_1, \ldots, \lambda_n) \in (\k^*)^n - \left\{(\pm 1, \ldots, \pm 1)\right\}$ and $(\alpha_1, \ldots, \alpha_{n}, c_1, \ldots, c_s) \in \k^{n} \times (\k^*)^s$. Thus, we have a fibration
$$
	U \longrightarrow \PGL{2} \times \Omega \longrightarrow \XTilde{\Xfp{n}{s}},
$$
where $\Omega = \left((\k^*)^n - \left\{(\pm 1, \ldots, \pm 1)\right\}\right) \times \left(\k^{n} \times (\k^*)^s\right)$. Observe that we do not need to remove any anti-diagonal value as we did in Section \ref{sec:free-groups}, since the intersection of the line spanned by $(\lambda_1-\lambda_1^{-1}, \ldots, \lambda_n - \lambda_n^{-1}, 0, \ldots, 0)$ with $\k^{n} \times (\k^*)^s$ is empty.
Hence, its virtual class is
$$
	\coh{\XTilde{\Xfp{n}{s}}} = \frac{q^3-q}{(q-1)q} \left((q-1)^n - 2^n\right)q^{n}(q-1)^s.
$$
\end{itemize}

Summing these results we get that
$$
	\coh{\Xred{\Xfp{n}{s}}} = {\left(q - 1\right)}^{n} {\left(q - 1\right)}^{s} {\left(q + 1\right)}
q^{n}.
$$
Using that $\Xfp{n}{s} = \SL{2}^n \times [J_+]_{\SL{2}}^s$ and that $\coh{[J_+]_{\SL{2}}} = (q^2-1)$, we obtain that $\coh{\Xfp{n}{s}} = (q^3-q)^n(q^2-1)^s$. Therefore, we get
$$
	\coh{\Xirred{\Xfp{n}{s}}} =  (q-1)^n(q-1)^s\left((q^2+q)^n(q+1)^s -(q+1)q^n\right).
$$
	\item For $\Xred{\Xfp{n}{s}}$ the situation is completely different from the previous ones. As we have shown, the stratum $\XD{\Xfp{n}{s}} = \emptyset$ so it can no longer be a core for the action. The key point now is that, precisely for this reason, the action of $\SL{2}$ on $\Xred{\Xfp{n}{s}}$ is closed. To be precise, take $A \in \Xred{\Xfp{n}{s}} \subseteq \Xred{\Xf{n + s}}$ and let $\overline{[A]}_{\SL{2}}$ be the closure of its orbit in $\Xf{n + s}$. The difference $\overline{[A]}_{\SL{2}} - [A]_{\SL{2}}$ lies in $\XD{\Xf{n + s}}$ so, since $\XD{\Xf{n + s}} \cap \Xfp{n}{s} = \XD{\Xfp{n}{s}} = \emptyset$, the orbit of $A$ is closed in $\Xred{\Xfp{n}{s}}$. However, the action of $\PGL{2}$ on $\Xred{\Xfp{n}{s}}$ is not free everywhere so we have to distinguish between two strata:
	\begin{itemize}
		\item $\XTilde{\Xfp{n}{s}}$. By { an analogous argument to Proposition \ref{prop:multiplicativity}} we have
		$$
			\coh{\XTilde{\Xfp{n}{s}} \sslash \SL{2}} = \frac{\coh{\XTilde{\Xfp{n}{s}}}}{q^3-q} = \left((q-1)^n - 2^n\right)q^{n-1}(q-1)^{s-1}.
		$$
		\item $\XPh{\Xfp{n}{s}}$. Here, the action of $\PGL{2}$ is not free but it has stabilizer isomorphic to $\Stab\,J_+ \cong \k$. The fact that $[J_+]_{\SL{2}} \cong \SL{2}/\Stab\,J_+$ implies that the GIT quotient
		$
			\XPh{\Xfp{n}{s}} \to \XPh{\Xfp{n}{s}} \sslash \SL{2}
		$
is a locally trivial fibration in the Zariski topology with fiber $\SL{2}/\Stab\,J_+$. Hence, since $\coh{\SL{2}/\Stab\,J_+} = \coh{\SL{2}}/\coh{\Stab\,J_+} = q^2-1$, a similar argument to Proposition \ref{prop:multiplicativity} implies that
		$$
			\coh{\XPh{\Xfp{n}{s}} \sslash \SL{2}} = \frac{\coh{\XPh{\Xfp{n}{s}}}}{q^2-1} = 2^{n}q^{n}(q-1)^{s-1}.
		$$
	\end{itemize}
From the stratification $\Xred{\Xfp{n}{s}} = \XPh{\Xfp{n}{s}} \sqcup \XTilde{\Xfp{n}{s}}$, where $\XTilde{\Xfp{n}{s}}$ is open orbitwise-closed, we get that
$$
    \coh{\Xred{\Xfp{n}{s}} \sslash \SL{2}} = \coh{\XPh{\Xfp{n}{s}} \sslash \SL{2}} + \coh{\XTilde{\Xfp{n}{s}} \sslash \SL{2}} =  
{\left(2^{n} + {\left(q - 1\right)}^{n - 1}\right)} {\left(q -
1\right)}^{s} q^{n - 1}.
$$
\end{itemize}
Summarizing, the analysis above shows that
\begin{empheq}{align*}
	\coh{\Xfp{n}{s} \sslash \SL{2}} = 2^{n} {\left(q - 1\right)}^{s} q^{n - 1} + {\left(q^{3} - q\right)}^{n -
1} {\left(q^{2} - 1\right)}^{s}.
\end{empheq}

\subsection{Surface groups}
Let $\Sigma$ be the closed orientable surface of genus $g \geq 1$ with $s \geq 1$ punctures and let us denote $\Xsp{g}{s} = \mathfrak{X}_{\SL{2}}(\pi_1(\Sigma), Q_s^+)$. As in Section \ref{section:stratification}, the decomposition into reducible and irreducible representations gives an equality
$$
	\coh{\Xsp{g}{s} \sslash \SL{2}} = \coh{\Xred{\Xsp{g}{s}} \sslash \SL{2}} + \coh{\Xirred{\Xsp{g}{s}} \sslash \SL{2}}.
$$

To understand this stratification, observe that $\XD{\Xsp{g}{s}} = \XDh{\Xsp{g}{s}} = \XI{\Xsp{g}{s}} = \emptyset$ since they are closed subvarieties of the ones of the free case. Set $\XU{\Xfp{n}{s}} = \XU{\Xf{n+s}} \cap \Xfp{n}{s}$ and $\XU{\Xsp{g}{s}} = \XU{\Xf{2g+s}} \cap \Xsp{g}{s}$ and take $A \in \XU{\Xfp{2g}{s}}$, namely
$$
	A = \left(
	\begin{pmatrix}
		\lambda_1 & \alpha_1 \\
		0 & \lambda_1^{-1}
	\end{pmatrix},
	\begin{pmatrix}
		\mu_1 & \beta_1 \\
		0 & \mu_1^{-1}
	\end{pmatrix}, \ldots,
	\begin{pmatrix}
		\lambda_{g} & \alpha_{g} \\
		0 & \lambda_{g}^{-1}
	\end{pmatrix},
	\begin{pmatrix}
		\mu_{g} & \beta_{g} \\
		0 & \mu_{g}^{-1}
	\end{pmatrix},
	\begin{pmatrix}
		1 & c_1 \\
		0 & 1
	\end{pmatrix}, \ldots,
	\begin{pmatrix}
		1 & c_s \\
		0 & 1
	\end{pmatrix}
	\right),
$$
with $\lambda_i, \mu_i \in \k^*$, $\alpha_i, \beta_i \in \k$ and $c_i \in \k^*$. Then we have that $A \in \XU{\Xsp{g}{s}}$ if and only if
\begin{equation}\label{eq:cond-upper:parabolic-jplus}
	{\sum_{i=1}^g \lambda_i\mu_i \left(\left(\lambda_i-\lambda_i^{-1}\right)\beta_i - \left(\mu_i-\mu_i^{-1}\right)\alpha_i\right) + \sum_{i=1}^s c_i= 0.}
\end{equation}
Let us analyze each stratum separately.

\begin{itemize}
	\item $\Xirred{\Xsp{g}{s}}$. By Proposition \ref{prop:multiplicativity}, the action of $\PGL{2}$ implies that
	$$
		\coh{\Xirred{\Xsp{g}{s}} \sslash \SL{2}} = \frac{\coh{\Xirred{\Xsp{g}{s}}}}{\coh{\PGL{2}}} = \frac{\coh{\Xsp{g}{s}} - \coh{\Xred{\Xsp{g}{s}}}}{q^3-q}.
	$$
The calculation of each of the strata of $\Xred{\Xsp{n}{s}}$ mimics the corresponding one for $\Xred{\Xfp{n}{s}}$ but now we have to take into account equation (\ref{eq:cond-upper:parabolic-jplus}).

\begin{itemize}
	\item $\XPh{\Xsp{g}{s}}$. In this case $\XP{\Xsp{g}{s}} = \left\{\pm 1\right\}^{2g} \times \left(\k^{2g} \times \pi_s\right)$, where
$$
    \pi_s = \left\{\left.\sum_{j = 1}^s c_j = 0\;\;\right|\, c_j \neq 0\right\}.
$$
To compute the virtual class of this space, observe that $\pi_s = (\k^*)^{s-1} - \pi_{s-1}$. Therefore, using the base case $\pi_1 = \emptyset$, we have
\begin{align*}
    \coh{\pi_s} &= (q-1)^{s-1} - \coh{\pi_{s-1}} = \sum_{k=1}^{s-1} (-1)^{k+1}(q-1)^{s-k} \\&= (-1)^s \left(\frac{(1-q)^s - 1}{q} + 1\right).
\end{align*}
As in the free case, the action of $\SL{2}$ on $\XP{\Xsp{n}{s}}$ is by rescaling on $\k^{n} \times \pi_s$, so we obtain a Zariski locally trivial fibration
$$
    \k^* \longrightarrow \left(\SL{2}/\Stab\,J_+\right) \times \left\{\pm 1\right\}^{2g} \times \left(\k^{2g} \times \pi_s\right) \longrightarrow \XPh{\Xsp{g}{s}}.
$$
Therefore, we have
$$
    \coh{\XPh{\Xsp{g}{s}}} = 2^{2g}(q^2-1)\frac{q^{2g}}{q-1}\left((-1)^s \left(\frac{(1-q)^s - 1}{q} + 1\right)\right).
$$
	
	\item $\XDh{\Xsp{g}{s}}$. This case makes no contribution since $\XD{\Xfp{n}{s}} = \emptyset$.

	\item $\XI{\Xsp{g}{s}}$. Again, this stratum is empty.
	\item $\XTilde{\Xsp{g}{s}}$. In this case, any element is conjugate to one of the form
$$
	\left(
	\begin{pmatrix}
		\lambda_1 & \alpha_1 \\
		0 & \lambda_1^{-1}
	\end{pmatrix}, \ldots,
	\begin{pmatrix}
		\mu_g & \beta_n \\
		0 & \mu_g^{-1}
	\end{pmatrix},
	\begin{pmatrix}
		1 & c_{1} \\
		0 & 1
	\end{pmatrix}, \ldots,
	\begin{pmatrix}
		1 & c_{s} \\
		0 & 1
	\end{pmatrix}
	\right)
$$
with $(\lambda_1, \ldots, \mu_g) \in (\k^*)^{2g} - \left\{(\pm 1, \ldots, \pm 1)\right\}$ and $(\alpha_1, \ldots, \beta_{g}, c_1, \ldots, c_s) \in  \Pi_s$, where we set
$$
    { \Pi_s = \left\{	\sum_{i=1}^g \lambda_i\mu_i \left(\left(\lambda_i-\lambda_i^{-1}\right)\beta_i - \left(\mu_i-\mu_i^{-1}\right)\alpha_i\right) + \sum_{i=1}^s c_i= 0
\right\},}
$$
for fixed $(\lambda_i, \mu_i)$.
To compute the virtual class of $\Pi_s$, observe that $\Pi_s = \k^{2g} \times (\k^*)^{s-1} - \Pi_{s-1}$. Using as base case that $\Pi_1$ is $\k^{2g}$ minus a hyperplane, we have
\begin{align*}
    \coh{\Pi_s} &= q^{2g}(q-1)^{s-1} - \coh{\Pi_{s-1}} \\
    &= q^{2g}\sum_{k=1}^{s} (-1)^{k+1}(q-1)^{s-k} + (-1)^s q^{2g-1} = q^{2g-1}(q-1)^s.
\end{align*}
There is a fibration
$$
	U \longrightarrow \PGL{2} \times \Omega \longrightarrow \XTilde{\Xsp{g}{s}},
$$
where $\Pi_s \to \Omega \to (\k^*)^{2g} - \left\{(\pm 1, \ldots, \pm 1)\right\}$ is a Zariski locally trivial fibration. Thus, we get
$$
	\coh{\XTilde{\Xsp{g}{s}}} = \frac{q^3-q}{(q-1)q} \left((q-1)^{2g} - 2^{2g}\right)\left(q^{2g-1}(q-1)^s\right).
$$
\end{itemize}
Summing these virtual classes, we obtain that
\begin{align*}
	\coh{\Xred{\Xsp{g}{s}}} &=          
2^{2g} \left(-1\right)^{s} {\left(q + 1\right)} q^{2g}
{\left(\frac{{\left(-q + 1\right)}^{s} - 1}{q} + 1\right)} \\&\;\;\;\;- {\left(2^{2
\, g} - {\left(q - 1\right)}^{2g}\right)} {\left(q - 1\right)}^{s}
{\left(q + 1\right)} q^{2g - 1}.
\end{align*}
In \cite[Theorem 5.10]{GP-2018} (see also \cite{Gonzalez-Prieto:Thesis}), it is proven that the virtual class of the whole representation variety is
\begin{align*}
\coh{\Xsp{g}{s}} =& \,{\left(q^2 - 1\right)}^{2g + s - 1} q^{2g - 1} +
\frac{1}{2} \, {\left(q -
1\right)}^{2g + s - 1}q^{2g -
1}(q+1){\left({2^{2g} + q - 3}\right)} 
\\ &+ \frac{\left(-1\right)^{s}}{2} \,
{\left(q + 1\right)}^{2g + s - 1} q^{2g - 1} (q-1){\left({2^{2g} +q -1}\right)}.
\end{align*}
Therefore, subtracting the contribution of the previous strata, we obtain
\begin{align*}
	\coh{\Xirred{\Xsp{g}{s}}} =&\, 2^{2g - 1} \left(-1\right)^{s} {\left(q + 1\right)}^{2g + s - 1}
{\left(q - 1\right)} q^{2g - 1} \\&- 2^{2g} \left(-1\right)^{s} {\left(q + 1\right)} q^{2g}
{\left(\frac{{\left(1-q\right)}^{s} - 1}{q} + 1\right)}\\&+ \frac{1}{2} \, \left(-1\right)^{s}
{\left(q + 1\right)}^{2g + s - 1} {\left(q - 1\right)}^{2} q^{2g
- 1}  \\&+ {\left(2^{2
\, g} - {\left(q - 1\right)}^{2g}\right)} {\left(q - 1\right)}^{s}
{\left(q + 1\right)} q^{2g - 1} \\
&+ \frac{1}{2} \, {\left(q -
1\right)}^{2g + s - 1} {\left(q + 1\right)} {\left(q - 3\right)}
q^{2g - 1}\\
	& + \frac{{\left(2^{2g} q^{2} + 2^{2g + 1} q + 2^{2
\, g} + 2 \, {\left(q + 1\right)}^{2g + s}\right)} {\left(q -
1\right)}^{2g + s - 1} q^{2g - 1}}{2 \, {\left(q + 1\right)}}.
\end{align*}

	\item For $\Xred{\Xsp{g}{s}}$ the situation is analogous to the case of $\Xred{\Xfp{2g}{s}}$. Since the action on $\Xred{\Xfp{2g}{s}}$ is closed, it is also so on $\Xred{\Xsp{g}{s}} = \Xred{\Xfp{2g}{s}} \cap \Xsp{g}{s}$. Therefore, stratifying in terms of the stabilizers for the action, we have:
	\begin{itemize}
		\item $\XTilde{\Xsp{g}{s}}$. By Proposition \ref{prop:multiplicativity} we get
		$$
			\coh{\XTilde{\Xsp{g}{s}} \sslash \SL{2}} = \frac{\coh{\XTilde{\Xsp{g}{s}}}}{q^3-q} = (-1)^s\left((q-1)^{2g} - 2^{2g}\right)q^{2g-2}(1-q)^{s-1}.
		$$
		\item $\XPh{\Xsp{g}{s}}$. Here the action of $\PGL{2}$ is not free but it has stabilizer isomorphic to $\Stab\,J_+ \cong \k$. As in the free case, analogously to Proposition \ref{prop:multiplicativity}, we get
		$$
			\coh{\XPh{\Xsp{g}{s}} \sslash \SL{2}} = \frac{\coh{\XPh{\Xsp{g}{s}}}}{q^2-1} = (-1)^s2^{2g}\frac{q^{2g}}{q-1}\left(\frac{(1-q)^s - 1}{q} + 1\right).
		$$
	\end{itemize}
Therefore, using the stratification $\Xred{\Xsp{g}{s}} = \XPh{\Xsp{g}{s}} \sqcup \XTilde{\Xsp{g}{s}}$ with $\XTilde{\Xsp{g}{s}}$ an open orbitwise-closed set, Theorem \ref{prop:decomposition-quotient} gives us
\begin{align*}
    \coh{\Xred{\Xsp{g}{s}} \sslash \SL{2}} &= \coh{\XPh{\Xsp{g}{s}} \sslash \SL{2}} + \coh{\XTilde{\Xsp{g}{s}} \sslash \SL{2}} \\&= 
\frac{2^{2g}q^{2g} \left(-1\right)^{s} {}
{\left(\frac{{\left(1-q\right)}^{s} - 1}{q} + 1\right)}}{q - 1} -
{{\left(2^{2g} - {\left(q - 1\right)}^{2g}\right)} {\left(q
- 1\right)}^{s-1} q^{2g - 2}}{}.
\end{align*}

\end{itemize}
Summarizing, the analysis above finally implies that
\begin{empheq}{align*}
	\coh{\Xsp{g}{s} \sslash \SL{2}} =&\,  {\left(q^2 - 1\right)}^{2g + s -
2} q^{2g - 2} +\left(-1\right)^{s} 2^{2g}  {\left(q - 1\right)} q^{2g - 2}
{\left({1-\left(1-q\right)}^{s - 1}\right)}\\
&+ \frac{1}{2}{\left(q - 1\right)}^{2g +s - 2} q^{2g - 2} \, {\left(2^{2g} + q - 3\right)}  \\
&+ \frac{1}{2} {\left(q + 1\right)}^{2g +
s - 2} q^{2g - 2}\,
\left(2^{2g} + q - 1\right).
\end{empheq}

\subsection{Parabolic data of Jordan type}
\label{subsec:parabolic-data-gen}

Let us denote by $\Gamma_{g,s}$ the fundamental group of the genus $g$ compact surface with $s$ removed points. Consider the parabolic structure $Q = \left\{(\gamma_1, [C_1]_{\SL{2}}), \ldots, (\gamma_s, [C_s]_{\SL{2}})\right\}$, where $C_i = J_+, J_-$ or $-\Id$. Let $r_+$ be the number of $J_+$, let $r_-$ be the number of $J_-$ and let $t$ be the number of $-\Id$ (so that $r_+ + r_- + t = s$). In addition, let us denote $\sigma = (-1)^{r_- + t}$. Observe that $J_+ \in [-J_-]_{\SL{2}}$ and $[-\Id]_{\SL{2}} = \left\{ - \Id\right\}$ so, depending on $\sigma$, we have:
\begin{itemize}
	\item If $\sigma = 1$, then we have $\mathfrak{X}_{\SL{2}}(\Gamma_{g,s}, Q) = \mathfrak{X}_{\SL{2}}(\Gamma_{g, r}, Q_{r}^+)$ where $r = r_+ + r_-$. Hence, we have $\coh{\mathfrak{X}_{\SL{2}}(\Gamma_{g,s}, Q) \sslash \SL{2}} = \coh{\Xsp{g}{r} \sslash \SL{2}}$ and the virtual class follows from the computation above.
	\item If $\sigma = -1$, then we have $\mathfrak{X}_{\SL{2}}(\Gamma_{g,s}, Q) = \mathfrak{X}_{\SL{2}}(\Gamma_{g, r + 1}, Q_{r}^-)$ where $r = r_+ + r_-$ and $Q_r^- = \left\{(\gamma_1, [J_+]_{\SL{2}}), \ldots, (\gamma_r, [J_+]_{\SL{2}}), (\gamma_{r+1}, \left\{-\Id\right\}) \right\}$. This is the so-called twisted representation variety.
This variety does not contain reducible representations. To check that,
let $A = (A_1, B_1, \ldots, C_1, \ldots, C_r, -\Id) \in \mathfrak{X}_{\SL{2}}(\Gamma_{g, r + 1}, Q_{r}^-)$ so that
$$
	\prod_{i=1}^{2g} [A_i, B_i] \prod_{j=1}^{s} C_i = -\Id.
$$
If $v \in \k^2 - \left\{0\right\}$ is a common eigenvector of $A$, then, since all the eigenvalues of the commutators $[A_i, B_i]$ and the {matrices $C_i$} are equal to $1$, the left hand side of the previous equation fixes $v$ but the right hand side does not. This proves that such $v$ cannot exist. Therefore, the action of $\PGL{2}$ on $\mathfrak{X}_{\SL{2}}(\Gamma_{g, r + 1}, Q_{r}^-)$ is closed and free.

As proven in \cite[Theorem 5.10]{GP-2018}, the virtual class of the representation variety is
\begin{align*}
	\coh{\mathfrak{X}_{\SL{2}}(\Gamma_{g, r + 1}, Q_{r}^-)} = &\, {\left(q - 1\right)}^{2g + r - 1} (q+1)q^{2g - 1}{{\left( {\left(q + 1\right)}^{2 \,
g + r-2}+2^{2g-1}-1\right)} } \\
&+ \left(-1\right)^{r + 1}2^{2g - 1}  {\left(q + 1\right)}^{2g + r
- 1} {\left(q - 1\right)} q^{2g - 1}.
\end{align*}
Therefore, Proposition \ref{prop:multiplicativity} gives us
\begin{align*}
	\coh{\mathfrak{X}_{\SL{2}}(\Gamma_{g,s}, Q) \sslash \SL{2}}  =&\, \left(-1\right)^{r-1}2^{2g - 1}  {\left(q + 1\right)}^{2g + r - 2} q^{2g - 2} \\ &+ {\left(q - 1\right)}^{2g + r - 2} q^{2g - 2}\left( {\left(q + 1\right)}^{2g + r - 2} + 2^{2g - 1} - 1\right).
\end{align*} 
\end{itemize}

%Hack to make arXiv LaTeX compiler work
%\includegraphics[scale=0]{./img/img1.png}

\bibliography{StratificationAlgebraicQuotients.bib}{}

\begin{thebibliography}{10}

\bibitem{Baraglia-Hekmati:2016}
D.~Baraglia and P.~Hekmati.
\newblock Arithmetic of singular character varieties and their
  {$E$}-polynomials.
\newblock {\em Proc. Lond. Math. Soc. (3)}, 114(2):293--332, 2017.

\bibitem{Bass:1968}
H.~Bass.
\newblock {\em Algebraic {$K$}-theory}.
\newblock W. A. Benjamin, Inc., New York-Amsterdam, 1968.

\bibitem{behrend512640motive}
K.~Behrend and A.~Dhillon.
\newblock On the motive of the stack of bundles.
\newblock {\em arXiv preprint math.AG/0512640}.

\bibitem{Beke2017Feb}
T.~Beke.
\newblock {The Grothendieck ring of varieties and of the theory of
  algebraically closed fields}.
\newblock {\em Journal of Pure and Applied Algebra}, 221(2):393--400, Feb 2017.

\bibitem{Berczi-Dolan-Hawes-Kirwan:2016}
G.~B\'erczi, B.~Doran, T.~Hawes, and F.~Kirwan.
\newblock Constructing quotients of algebraic varieties by linear algebraic
  group actions.
\newblock {\em Preprint arXiv:1512.02997}, 2016.

\bibitem{Biquard-Jardim}
O.~Biquard and M.~Jardim.
\newblock Asymptotic behaviour and the moduli space of doubly-periodic
  instantons.
\newblock {\em J. Eur. Math. Soc. (JEMS)}, 3(4):335--375, 2001.

\bibitem{Borel-1991}
A.~Borel.
\newblock {\em Linear algebraic groups}, volume 126 of {\em Graduate Texts in
  Mathematics}.
\newblock Springer-Verlag, New York, second edition, 1991.

\bibitem{Borisov}
L.~A. Borisov.
\newblock The class of the affine line is a zero divisor in the {G}rothendieck
  ring.
\newblock {\em J. Algebraic Geom.}, 27(2):203--209, 2018.

\bibitem{Cavazos-Lawton:2014}
S.~Cavazos and S.~Lawton.
\newblock {$E$}-polynomial of {${\rm SL}_2(\Bbb C)$}-character varieties of
  free groups.
\newblock {\em Internat. J. Math.}, 25(6):1450058, 27, 2014.

\bibitem{Corlette:1988}
K.~Corlette.
\newblock Flat {$G$}-bundles with canonical metrics.
\newblock {\em J. Differential Geom.}, 28(3):361--382, 1988.

\bibitem{Culler-Shalen}
M.~Culler and P.~B. Shalen.
\newblock Varieties of group representations and splittings of 3-manifolds.
\newblock {\em Annals of Mathematics}, 117(1):109--146, 1983.

\bibitem{DeligneII:1971}
P.~Deligne.
\newblock Th\'eorie de {H}odge. {II}.
\newblock {\em Inst. Hautes \'Etudes Sci. Publ. Math.}, (40):5--57, 1971.

\bibitem{DeligneIII:1971}
P.~Deligne.
\newblock Th\'eorie de {H}odge. {III}.
\newblock {\em Inst. Hautes \'Etudes Sci. Publ. Math.}, (44):5--77, 1974.

\bibitem{donaldson1986geometry}
S.~K. Donaldson and P.~Kronheimer.
\newblock The geometry of 4-manifolds.
\newblock In {\em Proceedings of the International Congress of Mathematicians
  (Berkeley 1986)(AM Gleason, ed.)}, volume~1, pages 43--54. Citeseer, 1986.

\bibitem{Dolan-Kirwan:2007}
B.~Doran and F.~Kirwan.
\newblock Towards non-reductive geometric invariant theory.
\newblock {\em Pure Appl. Math. Q.}, 3(1, Special Issue: In honor of Robert D.
  MacPherson. Part 3):61--105, 2007.

\bibitem{Drezet:2004}
J.-M. Dr\'ezet.
\newblock Luna's slice theorem and applications.
\newblock In {\em Algebraic group actions and quotients}, pages 39--89. Hindawi
  Publ. Corp., Cairo, 2004.

\bibitem{ekedahl2009grothendieck}
T.~Ekedahl.
\newblock The grothendieck group of algebraic stacks.
\newblock {\em arXiv preprint arXiv:0903.3143}, 2009.

\bibitem{Florentino-Silva}
C.~Florentino and J.~Silva.
\newblock Hodge-$\textrm{D}$eligne polynomials of abelian character varieties.
\newblock {\em Preprint arXiv:1711.07909}, 2017.

\bibitem{Gonzalez-Martinez2012Jan}
C.~Gonz{\ifmmode\acute{a}\else\'{a}\fi}lez-Mart{\ifmmode\acute{\imath}\else\'{\i}\fi}nez.
\newblock {The {H}odge{\textendash}{P}oincar{\ifmmode\acute{e}\else\'{e}\fi}
  polynomial of the moduli spaces of stable vector bundles over an algebraic
  curve}.
\newblock {\em Manuscripta math.}, 137(1):19--55, 2012.

\bibitem{GP-2019}
{\ifmmode\acute{A}\else\'{A}\fi}.~Gonz{\ifmmode\acute{a}\else\'{a}\fi}lez-Prieto.
\newblock {Virtual classes of parabolic $\mathrm{SL}_2(\mathbb{C})$-character
  varieties}.
\newblock {\em Adv. Math.}, 368:107148, Jul 2020.

\bibitem{GPLM-2017}
{\ifmmode\acute{A}\else\'{A}\fi}.~Gonz{\ifmmode\acute{a}\else\'{a}\fi}lez-Prieto,
  M.~Logares, and V.~Mu{\ifmmode\tilde{n}\else\~{n}\fi}oz.
\newblock {A lax monoidal topological quantum field theory for representation
  varieties}.
\newblock {\em Bulletin des Sciences
  Math{\ifmmode\acute{e}\else\'{e}\fi}matiques}, 161:102871, Jul 2020.

\bibitem{GPLM-2020}
{\ifmmode\acute{A}\else\'{A}\fi}.~Gonz{\ifmmode\acute{a}\else\'{a}\fi}lez-Prieto,
  M.~Logares, and V.~Mu{\ifmmode\tilde{n}\else\~{n}\fi}oz.
\newblock {Representation variety for the rank one affine group}.
\newblock {\em To appear in Mathematical Analysis in Interdisciplinary Research
  (I.N. Parasidis, E. Providas and Th.M. Rassias, eds.)}, 2020.

\bibitem{GP-2018}
{\'A}.~Gonz\'alez-Prieto.
\newblock Motivic theory of representation varieties via {T}opological
  {Q}uantum {F}ield {T}heories.
\newblock {\em Preprint arXiv:1810.09714v2}, 2018.

\bibitem{Gonzalez-Prieto:Thesis}
{\'A}.~Gonz\'alez-Prieto.
\newblock Topological {Q}uantum {F}ield {T}heories for character varieties.
\newblock {\em PhD Thesis. Universidad Complutense de Madrid}, 2018.

\bibitem{GPHV}
A.~Gonz\'alez-Prieto, M.~Hablicsek, and J.~Vogel.
\newblock Equivariant theory of virtual classes of character stacks.
\newblock In preparation.

\bibitem{EGAIV}
A.~Grothendieck.
\newblock {\'E}l\'ements de g\'eom\'etrie alg\'ebrique: {IV}. \'{E}tude locale
  des sch\'emas et des morphismes de sch\'emas, {T}roisi\`eme partie.
\newblock {\em Publications Math\'ematiques de l'IH\'ES}, 28:5--255, 1966.

\bibitem{Vogel:2020}
M.~Hablicsek and J.~Vogel.
\newblock Virtual classes of representation varieties of upper triangular
  matrices via {T}opological {Q}uantum {F}ield {T}heories, 2020.

\bibitem{Hausel-Letellier-Villegas}
T.~Hausel, E.~Letellier, and F.~Rodriguez-Villegas.
\newblock Arithmetic harmonic analysis on character and quiver varieties.
\newblock {\em Duke Math. J.}, 160(2):323--400, 2011.

\bibitem{Hausel-Rodriguez-Villegas:2008}
T.~Hausel and F.~Rodriguez-Villegas.
\newblock Mixed {H}odge polynomials of character varieties.
\newblock {\em Invent. Math.}, 174(3):555--624, 2008.
\newblock With an appendix by Nicholas M. Katz.

\bibitem{Hitchin}
N.~J. Hitchin.
\newblock The self-duality equations on a {R}iemann surface.
\newblock {\em Proc. London Math. Soc. (3)}, 55(1):59--126, 1987.

\bibitem{Hoskins}
V.~Hoskins.
\newblock Moduli problems and geometric invariant theory.
\newblock Lecture notes, 2015.

\bibitem{Jardim}
M.~Jardim.
\newblock Nahm transform and spectral curves for doubly-periodic instantons.
\newblock {\em Comm. Math. Phys.}, 225(3):639--668, 2002.

\bibitem{Kirwan:2009}
F.~Kirwan.
\newblock Quotients by non-reductive algebraic group actions.
\newblock In {\em Moduli spaces and vector bundles}, volume 359 of {\em London
  Math. Soc. Lecture Note Ser.}, pages 311--366. Cambridge Univ. Press,
  Cambridge, 2009.

\bibitem{Lawton-Munoz:2016}
S.~Lawton and V.~Mu\~noz.
\newblock {$E$}-polynomial of the {$SL(3,\Bbb{C})$}-character variety of free
  groups.
\newblock {\em Pacific J. Math.}, 282(1):173--202, 2016.

\bibitem{LM}
M.~Logares and V.~Mu\~noz.
\newblock Hodge polynomials of the {$\rm{SL}(2,\Bbb C)$}-character variety of
  an elliptic curve with two marked points.
\newblock {\em Internat. J. Math.}, 25(14):1450125, 22, 2014.

\bibitem{LMN}
M.~Logares, V.~Mu\~noz, and P.~E. Newstead.
\newblock Hodge polynomials of {${\rm SL}(2,\Bbb{C})$}-character varieties for
  curves of small genus.
\newblock {\em Rev. Mat. Complut.}, 26(2):635--703, 2013.

\bibitem{Luna:1972}
D.~Luna.
\newblock Sur les orbites ferm\'ees des groupes alg\'ebriques r\'eductifs.
\newblock {\em Invent. Math.}, 16:1--5, 1972.

\bibitem{Luna:1973}
D.~Luna.
\newblock Slices \'etales.
\newblock pages 81--105. Bull. Soc. Math. France, Paris, M\'emoire 33, 1973.

\bibitem{Martinez:2017}
J.~Mart\'inez.
\newblock E-polynomials of ${PGL}(2,\mathbb{C})$-character varieties of surface
  groups.
\newblock {\em Preprint arXiv:1705.04649}, 2017.

\bibitem{MM:2016}
J.~Mart\'inez and V.~Mu\~noz.
\newblock E-polynomials of {${SL}(2,\Bbb{C})$}-character varieties of complex
  curves of genus 3.
\newblock {\em Osaka J. Math.}, 53(3):645--681, 2016.

\bibitem{MM}
J.~Mart\'inez and V.~Mu\~noz.
\newblock E-polynomials of the {${\rm SL}(2,\Bbb C)$}-character varieties of
  surface groups.
\newblock {\em Int. Math. Res. Not. IMRN}, (3):926--961, 2016.

\bibitem{Munoz-Martinez:2015}
J.~Mart\'inez-Mart\'inez and V.~Mu\~noz.
\newblock The {$SU(2)$}-character varieties of torus knots.
\newblock {\em Rocky Mountain J. Math.}, 45(2):583--602, 2015.

\bibitem{Mereb}
M.~Mereb.
\newblock On the {$E$}-polynomials of a family of {$SL_n$}-character varieties.
\newblock {\em Math. Ann.}, 363(3-4):857--892, 2015.

\bibitem{milne:2017}
J.~S. Milne.
\newblock {\em Algebraic Groups: The Theory of Group Schemes of Finite Type
  over a Field}.
\newblock Cambridge Studies in Advanced Mathematics. Cambridge University
  Press, 2017.

\bibitem{MFK:1994}
D.~Mumford, J.~Fogarty, and F.~Kirwan.
\newblock {\em Geometric invariant theory}, volume~34 of {\em Ergebnisse der
  Mathematik und ihrer Grenzgebiete (2)}.
\newblock Springer-Verlag, Berlin, third edition, 1994.

\bibitem{Nagata:1963}
M.~Nagata.
\newblock Invariants of a group in an affine ring.
\newblock {\em J. Math. Kyoto Univ.}, 3:369--377, 1963/1964.

\bibitem{Newstead:1978}
P.~E. Newstead.
\newblock {\em Introduction to moduli problems and orbit spaces}, volume~51 of
  {\em Tata Institute of Fundamental Research Lectures on Mathematics and
  Physics}.
\newblock Tata Institute of Fundamental Research, Bombay; by the Narosa
  Publishing House, New Delhi, 1978.

\bibitem{Peters-Steenbrink:2008}
C.~A.~M. Peters and J.~H.~M. Steenbrink.
\newblock {\em Mixed {H}odge structures}, volume~52 of {\em Ergebnisse der
  Mathematik und ihrer Grenzgebiete. 3. Folge. A Series of Modern Surveys in
  Mathematics}.
\newblock Springer-Verlag, Berlin, 2008.

\bibitem{Schiffmann2016}
O.~Schiffmann.
\newblock {Indecomposable vector bundles and stable Higgs bundles over smooth
  projective curves}.
\newblock {\em Ann. Math.}, 183(1):297--362, 2016.

\bibitem{Serre}
J.-P. Serre.
\newblock Espaces fibr\'{e}s alg\'{e}briques (d'apr\`es {A}ndr\'{e} {W}eil).
\newblock In {\em S\'{e}minaire {B}ourbaki, {V}ol. 2}, pages Exp. No. 82,
  305--311. Soc. Math. France, Paris, 1995.

\bibitem{Sikora2012}
A.~Sikora.
\newblock {Character varieties}.
\newblock {\em Trans. Amer. Math. Soc.}, 364(10):5173--5208, 2012.

\bibitem{Simpson:parabolic}
C.~T. Simpson.
\newblock Harmonic bundles on noncompact curves.
\newblock {\em J. Amer. Math. Soc.}, 3(3):713--770, 1990.

\bibitem{Simpson:1992}
C.~T. Simpson.
\newblock Higgs bundles and local systems.
\newblock {\em Inst. Hautes \'Etudes Sci. Publ. Math.}, (75):5--95, 1992.

\bibitem{SimpsonI}
C.~T. Simpson.
\newblock Moduli of representations of the fundamental group of a smooth
  projective variety. {I}.
\newblock {\em Inst. Hautes \'Etudes Sci. Publ. Math.}, (79):47--129, 1994.

\bibitem{SimpsonII}
C.~T. Simpson.
\newblock Moduli of representations of the fundamental group of a smooth
  projective variety. {II}.
\newblock {\em Inst. Hautes \'Etudes Sci. Publ. Math.}, (80):5--79 (1995),
  1994.

\bibitem{Popov-Vinberg:1989}
E.~B. Vinberg and V.~L. Popov.
\newblock Invariant theory.
\newblock In {\em Algebraic geometry, 4 ({R}ussian)}, Itogi Nauki i Tekhniki,
  pages 137--314, 315. Akad. Nauk SSSR, Vsesoyuz. Inst. Nauchn. i Tekhn.
  Inform., Moscow, 1989.

\end{thebibliography}
\bibliographystyle{abbrv}

\end{document}